\documentclass[11pt, leqno, a4paper, fleqno]{amsart}
\usepackage{amssymb,amsthm,amsmath,amscd,amsfonts,bbm,mathrsfs}
\usepackage[all]{xy}
\def\noproof{{\unskip\nobreak\hfill\penalty50\hskip2em\hbox{}%
     \nobreak\hfill$\Box$\parfillskip=0pt%
     \finalhyphendemerits=0\par}}
\def\enddemo{\ifmmode\eqno\Box\else\noproof\vskip0.8truecm\fi}

\newtheorem{theo}{Theorem}[section]
\newtheorem{theorem}[theo]{Theorem}
\newtheorem{definition}[theo]{Definition}
\newtheorem{conjecture}[theo]{Conjecture}

\newtheorem{prop}[theo]{Proposition}
\newtheorem{corollary}[theo]{Corollary}

\newtheorem{remark}[theo]{Remark}
\newtheorem{remarks}[theo]{Remarks}

\newtheorem{lemma}[theo]{Lemma}

\newcommand{\lra}{\longrightarrow}

\DeclareMathOperator{\DR}{DR}

\DeclareMathOperator{\Eis}{Eis}

\DeclareMathOperator{\Coind}{Coind}
\DeclareMathOperator{\Ind}{Ind}

\DeclareMathOperator{\coinf}{coinf}
\DeclareMathOperator{\res}{res}

\DeclareMathOperator{\cusp}{cusp}
\DeclareMathOperator{\sign}{sign}

\DeclareMathOperator{\ord}{ord}
\DeclareMathOperator{\od}{o}

\DeclareMathOperator{\har}{har}
\DeclareMathOperator{\id}{id}
\DeclareMathOperator{\Hom}{Hom}

\DeclareMathOperator{\St}{St}
\DeclareMathOperator{\hol}{hol}
\DeclareMathOperator{\pr}{pr}
\DeclareMathOperator{\rank}{rank}
\DeclareMathOperator{\sing}{sing}

\DeclareMathOperator{\tor}{tor}
\DeclareMathOperator{\wcdot}{\, \cdot\, }
\DeclareMathOperator{\bcdot}{\bullet }

\DeclareMathOperator{\even}{even}
\DeclareMathOperator{\odd}{odd}

\DeclareMathOperator{\Ker}{Ker}
\DeclareMathOperator{\Image}{Im}

\DeclareMathOperator{\Coker}{Coker}
\DeclareMathOperator{\Dist}{Dist}
\DeclareMathOperator{\Meas}{\Dist^b}
\DeclareMathOperator{\Maps}{Maps}
\DeclareMathOperator{\Real}{Re}

\DeclareMathOperator{\Div}{Div}

\DeclareMathOperator{\Gal}{Gal}
\DeclareMathOperator{\Norm}{N}

\DeclareMathOperator{\GL}{GL}
\DeclareMathOperator{\PGL}{PGL}
\DeclareMathOperator{\SL}{SL}

\DeclareMathOperator{\Orth}{O}
\DeclareMathOperator{\SOrth}{SO}

\DeclareMathOperator{\disc}{disc}

\DeclareMathOperator{\einhalb}{\mbox{\small{$\frac{1}{2}$}}}

\newcommand{\Cd}{{C_{\diamond}}}
\newcommand{\sCd}{{{\sC}_{\diamond}}}

\newcommand{\Tw}{{\mathfrak{Tw}}}

\newcommand{\fa}{{\mathfrak a}}
\newcommand{\fA}{{\mathfrak A}}
\newcommand{\fc}{{\mathfrak c}}

\newcommand{\ff}{{\mathfrak f}}

\newcommand{\fm}{{\mathfrak m}}
\newcommand{\fn}{{\mathfrak n}}
\newcommand{\fp}{{\mathfrak p}}
\newcommand{\fq}{{\mathfrak q}}

\newcommand{\bC}{{\mathbb C}}
\newcommand{\bD}{{\mathbb D}}

\newcommand{\bG}{{\mathbb G}}
\newcommand{\bH}{{\mathbb H}}

\newcommand{\bN}{{\mathbb N}}
\newcommand{\bP}{{\mathbb P}}
\newcommand{\bQ}{{\mathbb Q}}
\newcommand{\bR}{{\mathbb R}}

\newcommand{\bZ}{{\mathbb Z}}
\newcommand{\barQ}{{\overline{\mathbb Q}}}

\newcommand{\bA}{{\mathbf A}}
\newcommand{\bI}{{\mathbf I}}

\newcommand{\barE}{{\overline{E}}}
\newcommand{\barH}{{\overline{H}}}

\newcommand{\barcO}{{\overline{\cO}}}

\newcommand{\bff}{{\bf f}}
\newcommand{\bfC}{{\bf C}}
\newcommand{\bfP}{{\bf P}}

\newcommand{\bGm}{{\bG_m}}

\newcommand{\cA}{{\mathcal A}}

\newcommand{\cD}{{\mathcal D}}
\newcommand{\cE}{{\mathcal E}}
\newcommand{\cF}{{\mathcal F}}
\newcommand{\cG}{{\mathcal G}}
\newcommand{\cH}{{\mathcal H}}

\newcommand{\cL}{{\mathcal L}}
\newcommand{\cM}{{\mathcal M}}
\newcommand{\cN}{{\mathcal N}}
\newcommand{\cO}{{\mathcal O}}

\newcommand{\cT}{{\mathcal T}}
\newcommand{\cV}{{\mathcal V}}
\newcommand{\cW}{{\mathcal W}}
\newcommand{\cX}{{\mathcal X}}

\newcommand{\sA}{{\mathscr A}}
\newcommand{\sB}{{\mathscr B}}
\newcommand{\sC}{{\mathscr C}}
\newcommand{\sD}{{\mathscr D}}
\newcommand{\sE}{{\mathscr E}}

\newcommand{\sQ}{{\mathscr Q}}

\DeclareMathOperator{\wsE}{\widetilde{\sE}}

\newcommand{\wkappa}{{\widetilde{\kappa}}}
\newcommand{\wdelta}{{\widetilde{\delta}}}

\newcommand{\oE}{{\vec{\mathcal E}}}

\newcommand{\ua}{{\underline{a}}}

\newcommand{\ui}{{\underline{i}}}

\newcommand{\un}{{\underline{n}}}

\newcommand{\uz}{{\underline{z}}}
\newcommand{\uzero}{{\underline{0}}}

\newcommand{\uzwei}{{\underline{2}}}
\newcommand{\ep}{{\epsilon}}

\newcommand{\ualpha}{{\underline{\alpha}}}
\newcommand{\uep}{{\underline{\epsilon}}}
\newcommand{\umu}{{\underline{\mu}}}
\newcommand{\unu}{{\underline{\nu}}}

\newcommand{\barho}{{\overline{\rho}}}

\DeclareMathOperator{\gr}{gr}
\newcommand{\bu}{{\bullet}}
\newcommand{\noi}{\noindent}

\begin{document}

\title{On special zeros of $p$-adic $L$-functions of Hilbert modular forms}
\author{By Michael Spie{\ss}}
\address{Fakult\"{a}t f\"{u}r Mathematik,
Universit\"{a}t Bielefeld, D-33501 Bielefeld, Germany}
\email{mspiess@math.uni-bielefeld.de}
\maketitle
\begin{abstract}

Let $E$ be a modular elliptic curve over a totally real number
field $F$. We prove the weak exceptional zero
conjecture which links a (higher) derivative of the $p$-adic
$L$-function attached to $E$ to certain $p$-adic periods attached 
to the corresponding Hilbert modular form at the places
above $p$ where $E$ has split multiplicative reduction.
Under some mild restrictions on $p$ and the conductor of $E$ we deduce
the exceptional zero conjecture in the strong form (i.e.\ where the
automorphic $p$-adic periods are replaced by the $\cL$-invariants of
$E$ defined in terms of Tate periods) from a special case proved
earlier by Mok. Crucial for our method is a new construction of the $p$-adic
$L$-function of $E$ in terms of local data.
\end{abstract}

\tableofcontents

\section*{Introduction}
 
Let $E$ be a modular elliptic curve over a totally real number field $F$ and
let $p$ be a prime number and such that $E$ has either good ordinary
or multiplicative reduction at all places $\fp$ above $p$. Attached to $E$
are the (Hasse-Weil) $L$-function $L(E, s)$ (a function in the
complex variable $s$) and a $p$-adic $L$-function $L_p(E,s)$ (here
$s\in \bZ_p$). Both are linked by the {\it interpolation property}
which expresses the $p$-adic measure associated to $L_p(E,s)$ in terms
of twisted special $L$-values $L(E,\chi, 1)$.  A
special case is the formula
\[
L_p(E, 0)  \,\,\, = \,\,\, \prod_{\fp|p} \, e(\alpha_{\fp},1) \, \cdot L(E,1).
\]
Here $e(\alpha_{\fp},1)$ is certain Euler factor defined in terms of the
reduction of $E$ at $\fp$ (see Prop.\ \ref{prop:interpol} for its
definition). It is $=0$ if and only if $E$ has split multiplicative reduction at
$\fp$. Let $S_1$ be the set of primes $\fp$ of $F$ above $p$ where $E$ has split
multiplicative reduction, let $S_p$ be the set of all primes above $p$
and let $S_2= S_p-S_1$. Thus we have $L_p(E, 0)=0$ if $S_1\ne
\emptyset$. In \cite{hida} it has been conjectured that 
\begin{eqnarray}
\label{eqn:vanorder}
&&\ord_{s=0} L_p(E,s) \ge r \,\colon  \!\! = \,\,\sharp(S_1);\\
\label{eqn:taylorcoef}
&&\frac{d^r}{ds^r} L_p(E, s)|_{s=0} \,\, =\,\, r! \,\,
  \prod_{\fp\in S_1}\cL_{\fp}(E) \,\cdot\,\prod_{\fp\in S_2}
  e(\alpha_{\fp},1) \,\cdot \, L(E, 1).
\end{eqnarray}
Here the {\it $\cL$-invariant $\cL_{\fp}(E)$} is defined as $\cL_{\fp}(E) =
\ell_{\fp}(q_{E/F_{\fp}}))/\od_{\fp}(q_{E/F_{\fp}})$
where $q_{E/F_{\fp}}$ is the Tate period of $E/F_{\fp}$, $\ell_{\fp} = \log_p\circ N_{F_{\fp}/\bQ_p}$ and $\od_{\fp} = \ord_{\fp}$ is the normalized additive valuation of $F_{\fp}$. 

In this paper we prove (\ref{eqn:vanorder}) unconditionally and (\ref{eqn:taylorcoef}) under the following
assumptions (see Theorem \ref{theorem:ezc}):
 (i) $p\ge 5$ is unramified in $F$; (ii) $E$ has multiplicative reduction at a prime $\fq\nmid p$, or $r$ is odd, or the
sign $w(E)$ of the functional equation for $L(E,s)$ (i.e.\ the root
number of $E$) is $= -1$.

The statements (\ref{eqn:vanorder}) and (\ref{eqn:taylorcoef}) are known as {\it exceptional zero conjecture}. In the case
$F= \bQ$ it was formulated by Mazur, Tate and Teitelbaum \cite{mtt}
and proved by Greenberg and Stevens \cite{grstevens} and independently
by Kato, Kurihara and Tsuji. In the case $r=1$, (\ref{eqn:taylorcoef}) was proved by Mok \cite{mok} under the assumption (i), by extending the method of \cite{grstevens}.

To explain our proof let $\pi$ be the automorphic
representation of $\GL_2(\bA_F)$ associated to $E$. Thus $\pi$ has trivial central character and the local factor $\pi_v$ is discrete of weight 2 at all archimedean places $v$. The Hasse-Weil $L$-function of $E$ is then equal to the
automorphic $L$-function $L(s-\einhalb,\pi)$. Moreover $L_p(E,s)$ is
solely defined in terms of $\pi$ (thus we write $L_p(s,\pi)$ for $L_p(E,s)$).

In section \ref{subsection:linv} we shall introduce a second type of
$\cL$-invariant $\cL_{\fp}(\pi)$. It is defined in terms of the cohomology of
($S_p$-)arithmetic groups. We show that $\cL_{\fp}(\pi)$ does not change under certain
quadratic twists, i.e.\ we have $\cL_{\fp}(\pi\otimes \chi) =
\cL_{\fp}(\pi)$ for any quadratic character $\chi$ of the idele class
group $\bI/F^*$ of $F$ such that the local components $\chi_v$ of
$\chi$ at infinite places and at $v=\fp$ are trivial. 
We prove an analogue of (\ref{eqn:taylorcoef}) above
(unconditionally) with the arithmetic $\cL$-invariants $\cL_{\fp}(E)$
replaced by the automorphic $\cL$-invariants $\cL_{\fp}(\pi)$, i.e.\
we show
\begin{eqnarray}
\label{eqn:weakezc}
&&\frac{d^r}{ds^r} L_p(s,\pi)|_{s=0} \,\, =\,\, r! \,\,
  \prod_{\fp\in S_1}\cL_{\fp}(\pi) \,\cdot\,\prod_{\fp\in S_2}
  e(\alpha_{\fp},1) \,\cdot \, L(\einhalb, \pi).
\end{eqnarray}
In the case $F=\bQ$ these $\cL$-invariants have been introduced by 
Darmon (\cite{darmon}, section 3.2). He showed that they are 
invariant under twists and also proved (\ref{eqn:weakezc}). Also if the narrow class
number of $F$ is $=1$ a different construction  of $\cL_{\fp}(\pi)$
has been given in \cite{greenberg}.

To deduce (\ref{eqn:taylorcoef}) from (\ref{eqn:weakezc}) it is therefore enough to show $\cL_{\fp}(\pi) = \cL_{\fp}(E)$ for all $\fp\in S_1$. In future work
\cite{iovitaspiess} we plan to give an unconditional proof of it
(and thus of (\ref{eqn:taylorcoef})) by comparing $\cL_{\fp}(\pi)$ to the (similarly defined)
$\cL$-invariant of an automorphic representation $\pi'$ of a
totally definite quaternion algebra -- which corresponds to $\pi$ under
Jacquet-Langlands functoriality -- and by using $p$-adic
uniformization of Shimura curves (compare also \cite{bdi} where a
similar proof has been given in the case $F=\bQ$ under certain assumptions on $\pi$). 

However if $p$ satisfies the conditions (i) above and $E$ satisfies
(ii) then we can deduce the equality $\cL_{\fp}(\pi) =\cL_{\fp}(E)$
for fixed $\fp\in S_1$ by comparing the formulas (\ref{eqn:taylorcoef}) and (\ref{eqn:weakezc}) in the case
$r=1$ for certain quadratic twists of $E$ and $\pi$. More precisely,
by a result of Waldspurger \cite{waldspurger}, we can choose a
quadratic character $\chi$ such that the arithmetic and
automorphic $\cL$-invariants at $\fp$ do not change under twisting with
$\chi$, $L(\einhalb,\pi\otimes \chi)$ does not vanish and $\fp$ is the only place above $p$
where the twisted elliptic curve $E_{\chi}$ has
split multiplicative reduction. Then
by Mok's result and (\ref{eqn:weakezc}) we can express both $\cL_{\fp}(E)$ and
$\cL_{\fp}(\pi)$ by the same formula.

The $p$-adic $L$-function attached to $\pi$ is the $\Gamma$-transform of a certain canonical measure
$\mu_{\pi}$ on the Galois group $\cG_p$ of the maximal abelian
extension of $F$ which is unramified outside $p$ and $\infty$, i.e.\
it is given by
\[
L_p(s, \pi) \,\, = \,\,\int_{\cG_p}\, \langle \gamma\rangle^s\mu_{\pi}(d\gamma)
\]
(for the definition of $\langle \gamma\rangle^s$ see section \ref{subsection:abstrezc}).

Crucial for the proof of (\ref{eqn:vanorder}) and (\ref{eqn:weakezc}) is a new construction of $\mu_{\pi}$\footnote{In principle our construction is related to Manin's \cite{manin}. However in our set-up the measure $\mu_{\pi}$ is build in a simple manner from local distributions $\mu_{\pi_v}$ at each place $v$ of $F$}. We shall briefly explain it (for details see
\ref{subsection:pmeashmf}). 
Heuristically, we define $\mu_{\pi}$ as the
direct image of the distribution $\mu_{\pi_p}\times W^p\left(\begin{matrix}
    x & 0\\ 0 & 1\end{matrix}\right) d^{\times} x$
under the reciprocity map $\bI = F_p^*\times \bI^p \to \cG_p$
of class field theory.
Here the first factor $\mu_{\pi_p}$ is the product distribution on $F_p^* =
\prod_{\fp\in S_p} F_{\fp}^*$ of certain canonical distributions
$\mu_{\fp}$ on $F_{\fp}^*$ attached to each local factors
$\pi_{\fp}$, $\fp\in S_p$. Moreover $d^{\times} x$ denotes the Haar
measure on the group of $S_p$-ideles $\bI^p = \prod'_{v\nmid p} \,
F_v^*$ and $W^p$ is a certain Whittaker function of $\pi^p =
\bigotimes_{v\nmid p} \pi_v$ (it is the product of local
Whittaker functions). 
 
To put this construction on a sound foundation consider the map $\phi_{\pi}$ given by 
\[
\phi_{\pi}(U,x^p) \,\,\, = \,\,\, \sum_{\zeta\in F^*}\,\, \mu_{\pi_p}(\zeta U)\,\, W^p\left(\begin{matrix} \zeta x^p & 0\\ 0 & 1\end{matrix}\right)
\]
where the first argument $U$ is a compact open subset of $F_p^*$
and the second an idele $x^p\in \bI^p$. Then $\phi_{\pi}(\zeta U, \zeta x^p) =
\phi_{\pi}(U,x)$ for all $\zeta\in F^*$. Thus if we set $\phi_U(x_p, x^p) \colon \!\! = \phi_{\pi}(x_p U, x^p)$
then $\phi_U$ can be viewed as a function on the idele class group
$\bI/F^*$ (so the map $U\to \phi_U$ is a distribution on $F_p^*$
with values in a certain space of functions on $\bI/F^*$). 

 For a locally constant map $f: \cG_p\to
\bC$ there exists a compact open subgroup $U \subset U_p = \prod_{\fp\in S_p} \cO_{\fp}^*\subset F_p^*$ such that
$f\circ \rho: \bI/F^* \to \bC$ factors through $\bI/F^*(U\times
U^p)$ (here $\rho: \bI/F^*\to \cG_p$ denotes the reciprocity map). 
Then $\int_{\cG_p} \,f(\gamma)\,\mu_{\pi}(d\gamma)$ is given by 
\[
\int_{\cG_p} \, \,f(\gamma)\,\mu_{\pi}(d\gamma) = [U_p:U]
\int_{\bI/F^*} \, f(\rho(x))\phi_U(x) \, d^{\times} x.
\]
By using properties of the cohomology groups of arithmetic subgroups of $\GL_2(F)$ we show that $\mu_{\pi}$ is {\it bounded} (i.e.\ it is a {\it $p$-adic measure} in the
sense of section \ref{subsection:padicmeas} below) and so any 
continuous map $\bZ_p\to \bC_p$ can be integrated against it. 

One way to describe the local distribution $\mu_{\fp}$ for $\fp\in
S_p$ is that it is the image of a certain 
Whittaker functional of $\pi_{\fp}$ under a canonical map -- denoted
by $\delta$ -- from the dual of $\pi_{\fp}$ to the space of
distributions on $F_{\fp}^*$. We will give the definition of $\delta$ in
the case $\fp\in S_1$, or equivalently, when $\pi_{\fp}$ is the Steinberg representation $\St$
(i.e.\ $\pi_{\fp}$ is isomorphic to the space of locally constant functions $\bP^1(F_{\fp}) \to
\bC$ modulo constants). For $c\in \Hom(\St, \bC)$ we
define $\delta(c)$ by $\int_F \, f(x) \delta(c)(dx)=
c(\tilde{f})$. Here for a locally constant map with
compact support $f: F_{\fp}\to \bC$ we define $\tilde{f}: \bP^1(F_{\fp}) \to
\bC$ by $\tilde{f}(\infty) =0$ and  $\tilde{f}([x:1]) = f(x)$.
Thus in the case $\pi_{\fp} = \St$, the target of
$\delta$ is the space of distributions on $F_{\fp}$. 

In particular the local contribution $\mu_{\fp}$ of $\mu_{\pi}$ at $\fp\in S_1$ is actually a distribution on
$F_{\fp}$ (and not only on $F_{\fp}^*$). Therefore, allowed as first argument in
$\phi_{\pi}(U, x^p)$ are not only compact open subsets $U$ of
$F_p^*$ but also of the larger space $\prod_{\fp\in S_1} F_{\fp} \times \prod_{\fp\in S_2}
F_{\fp}^*$. This fact is crucial for our proof that the vanishing order
$L_p(s, \pi)$ at $s=0$ is $\ge r$. The map $\delta$ and distributions $\mu_{\fp}$ will be introduced in sections \ref{subsection:diststeinberg} and \ref{subsection:localdistr} respectively.

Chapter \ref{section:abstractezc} is the technical heart of this
paper. It provides an axiomatic approach to study trivial zeros of
$p$-adic $L$-function which can be applied in other situations
as well (e.g.\ to the case of $p$-adic $L$-functions of totally real
number fields \cite{spiess}, \cite{daschar}). We consider 
arbitrary two-variable function $\phi: (U,x^p) \mapsto \phi(U, x^p)$
($U\subset \prod_{\fp\in S_1} F_{\fp} \times \prod_{\fp\in S_2}F_{\fp}^*$
compact open and $x^p\in \bI^p$) satisfying certain axioms and attach 
a $p$-adic distribution $\mu$ on $\cG_p$ as above. By "integrating away" the infinite places we obtain a certain cohomology class $\kappa\in H^d(F_+^*,\cD)$ associated to $\phi$ (where $d= [F:\bQ]-1$, $F^*_+$ denotes the group of totally positive elements of $F$ 
and $\cD$ is a certain space of distributions on the adelic space $\prod_{\fp\in S_1} F_{\fp} \times \prod_{\fp\in S_2}F_{\fp}^*\times \prod'_{v\nmid p\infty} \,
F_v^*$) and the distribution $\mu$ can be defined solely in terms of $\kappa$. The space $\cD$ contains a canonically subspace $\cD^b$ (consisting -- in a certain sense -- of $p$-adic measures) and $\mu$ is a $p$-adic measure provided that $\kappa$ lies in the image of $H^d(F_+^*,\cD^b)\to H^d(F_+^*,\cD)$ (see section \ref{subsection:intcohom}). 

In this case we define $L_p(s, \phi)$ as the $\Gamma$-transform of $\mu$ and show that $L_p(s, \phi)$ has a zero of order $\ge r$ at $s=0$. Furthermore we give a
description of the $r$-th derivative $\frac{d^r}{ds^r} L_p(s, \phi)|_{s=0}$ as a certain cap-product. More precisely, we associate to
any continuous homomorphism $\ell: F_{\fp}^*\to \bC_p$ a 
cohomology class $c_{\ell}\in H^1(F^*_+, C_c(F_{\fp}, \bC_p))$ (for its definition and the notation see \ref{definition:logclass}). If $S_1 =
\{\fp_1, \ldots, \fp_r\}$ we will show
\begin{eqnarray}
\label{eqn:taylorcup}
&&\frac{d^r}{ds^r} L_p(s, \phi)|_{s=0} \,\,\, = \,\,\, (-1)^{\binom{r}{2}}\,\,r! \,\, (\kappa \cup c_{\ell_{\fp_1}} \cup \ldots \cup c_{\ell_{\fp_{r}}}) \cap \vartheta.
\end{eqnarray}
Here $\vartheta$ is essentially the fundamental class of the quotient $M/F^*_+$ where 
$M$ is a certain $d+r$-dimensional manifold on which $F^*_+$ acts freely (see section \ref{subsection:capprod}). If $U_0 = \prod_{\fp\in S_1} \cO_{\fp} \times
\prod_{\fp\in S_2} \cO_{\fp}^*$ and $\phi_0(x) \colon \!\!= \phi(x_p U_0,
x^p)$ for $x= (x_p, x^p) \in F_p^*\times \bI^p = \bI$, we will also
prove 
\begin{eqnarray}
\label{eqn:lwertcup}
&&
\int_{\bI/F^*} \, \phi_0(x) \, d^{\times}x\,\,\, =\,\,\, (-1)^{\binom{r}{2}}\,\,r! \,\, (\kappa \cup c_{\od_{\fp_1}} \cup \ldots \cup c_{\od_{\fp_r}}) \cap \vartheta.
\end{eqnarray}

In chapter \ref{section:pmeashmf} we will verify that the theory
developed in the previous chapter can be applied in the case $\phi
=\phi_{\pi}$. The difficult part is to show that the cohomology
class $\kappa_{\pi}$ attached to $\phi_{\pi}$ comes from a class in $H^d(F_+^*,\cD^b)$. 
This is achieved by showing that it lies in the image of a specific cohomology class $\widehat{\kappa}_{\pi}\in H^d(\PGL_2(F), \cA)$ under a canonical map $\Delta_*: H^d(\PGL_2(F), \cA)\to H^d(F_+^*,\cD)$ (for the definition of the coefficients $\cA$ and the map $\Delta_*$ we refer to section \ref{subsection:finiteautomorph} and \ref{subsection:eichshim}). The fact that any arithmetic subgroup of $\PGL_2(F)$ has the finiteness property (VFL) (introduced by Serre in \cite{serre}) implies that $\Delta_*$ factors through $H^d(F_+^*,\cD^b)$.

In the last chapter \ref{section:ezc} we will introduce the
automorphic $\cL$-invariant $\cL_{\fp}(\pi)$ and deduce (\ref{eqn:weakezc}) from (\ref{eqn:taylorcup})
and (\ref{eqn:lwertcup}). The cohomology group $H^d(\PGL_2(F), \cA)$ carries an action of a Hecke algebra and $\widehat{\kappa}_{\pi}$ lies in the $\pi$-isotypic component $H^d(\PGL_2(F), \cA)_{\pi}$.
 Using the fact that the classes $c_{\ell}$ ``come'' from certain $\PGL_2$ cohomology classes as well (they will be introduced in section \ref{subsection:steinbergext}) and the fact that $H^d(\PGL_2(F), \cA)_{\pi}$ is onedimensional (a results due to Harder \cite{harder}) we show that the cup products $\kappa \cup c_{\ell_{\fp}}$ and $\kappa \cup c_{\fp}$ differ by a factor $\cL_{\fp}(\pi)$ which is defined in terms of the cohomology of $\PGL_2(F)$. 
\medskip

\paragraph{{\it Acknowledgement}} I thank Vytautas Paskunas for several helpful conversations and Kumar Murty for providing me with the reference \cite{friedhoff}. Also I am grateful to H.\ Deppe, L.\ Gehrmann, S.\ Molina and M.\ Seveso for useful comments on an earlier draft.
\medskip

\paragraph{{\it Notation}}

The following notations are valid throughout this paper. A list with further notations will be given at the beginning of chapters \ref{section:localmass} and \ref{section:abstractezc}.

Unless otherwise stated all rings are commutative with unit.

We fix a prime number $p$ and embeddings 
\[
\iota_{\infty}: \barQ\hookrightarrow \bC, \qquad \iota_p: \barQ\hookrightarrow \bC_p.
\]
We let $\ord_p$ denote the valuation on $\bC_p$ and $\barQ$ (via
$\iota_p$) normalized so that $\ord_p(p) =1$. The valuation ring of
$\barQ$ with respect to $\ord_p$ will be denoted by $\barcO$.

If $X$ and $Y$ are topological spaces then $C(X,Y)$ denotes the set of continuous maps $X\to Y$.  If we consider $Y$ with the discrete topology then we shall also write $C^0(X,Y)$ instead of $C(X,Y)$. If $Y=R$ is a topological ring then $C_c(X,R)$ is the submodule of $C(X,R)$ of continuous maps with compact support. If we consider $R$ with the discrete topology then we shall also write $C^0_c(X,Y)$ instead of  $C_c(X,R)$.

Put $G \colon \!\! = \PGL_2$, and let $B$ be the subgroup of upper triangular matrices (modulo the center $Z$ of $\GL_2$), $T = \left\{ \left(\begin{matrix} *& 0\\
0 & *\end{matrix}\right)\right\}/Z$ be the maximal torus of $G$ in $B$. We write elements of $G$ often simply as matrices $\left(\begin{matrix} a & b\\
c & d\end{matrix}\right)$ (and neglect the fact that we consider them only modulo the center of $\GL_2$). We identift $\bG_m$ with $T$ via the isomorphism $t\mapsto \left(\begin{matrix} t& 0\\ 0 & 1\end{matrix}\right)$. If $R$ is a ring the determinant induces a homomorphism $\det: G(R)\to R^*/(R^*)^2$. 

\section{Generalities on distributions and measures}
\label{section:distmeas}

\subsection{Distributions and measures} Let $\cX$ be a totally disconnected $\sigma$-locally compact topological space (in practice $\cX$ will be a e.g.\ profinite set like an infinite Galois group or a certain space of adeles). For a topological Hausdorff ring $R$ we denote by $\Cd(\cX, R)$ the subring of $C(\cX, R)$ consisting of maps $f:\cX \to R$ with $f(x) \to 0$ as $x\to \infty$ (equivalently by setting $f(\infty) = 0$ the map $f$ extends continuously to the one-point compactification of $\cX$). We have $C^0_c(\cX, R) \subseteq C_c(\cX, R) \subseteq \Cd(\cX, R) \subseteq C(\cX, R)$. Note that if $\cX = \cX_1\times \cX_2$ where $\cX_1$ and $\cX_2$ are $\sigma$-locally compact and if $f_1\in \Cd(\cX_1, R)$, $f_2\in \Cd(\cX_1, R)$ then the map $(f_1\otimes f_2)(g_1, g_2) \colon \!\! = f_1(g_1) \cdot f_2(g_2)$ lies in $\Cd(\cX, R)$. 

Let $M$ be an $R$-module. Recall that an $M$-valued {\it distribution} on $\cX$ is a homomorphism $\mu: C^0_c(\cX, \bZ) \to M$. It extends to an $R$-linear map 
\begin{equation}
\label{distr1}
C^0_c(\cX, R) \lra M,\,\,\, f\mapsto \int_{\cX}\, f\,\, d\mu.
\end{equation}
We shall denote the $R$-module of $M$-valued distributions on $\cX$ by $\Dist(\cX, M)$. If $\cX = \cX_1\times \cX_2$, $\mu\in \Dist(\cX, M)$ and $f_1\in C^0_c(\cX_1, R)$ then $f_2\mapsto \int\, f_1\otimes f_2\,\, d\mu$ is an $M$-valued distribution on $\cX_2$ which will be denoted by $\mapsto \int_{\cX_1}\, f_1\,d\mu$ i.e.\ we have a pairing
\begin{equation}
\label{fubini}
\Dist(\cX, M)\times C^0_c(\cX_1, R) \lra \Dist(\cX_2, M),\, (\mu, f_1)\mapsto \int_{\cX_1}\, f_1\,d\mu.
\end{equation}

Next we introduce the notion of a measure on $\cX$ with values in a $p$-adic Banach space. Assume that $R = K$ is a $p$-adic field. By that we mean that $K$ is a field of characteristic 0 which is equipped with a $p$-adic value, i.e.\ a nonarchimedian absolute value $|\,\,\,|: K \to \bR$ whose restriction to $\bQ$ is the usual $p$-adic value and such $K$ is complete with respect to $|\,\,\,|$. We denote a $p$-adic value often as $|\,\,\,|_p$ and the corresponding valuation ring by $\cO_K$. 

A norm on a $K$-vector space $V$ is a function $\|\,\,\,\|: V \to \bR$ such that (i) $\|av\| = |a|_p\|v\|$, (ii) $\|v+w\| \le \max(\|v\|, \|w\|)$ and (iii) $\|v \|\ge 0$ with equality iff $v=0$ for all $a\in K$, $v, w\in V$. Two norms 
$\|\,\,\,\|_1$, $\|\,\,\,\|_2$ are {\it equivalent} if there exists $C_1, C_2\in \bR_+$ with $C_1\|v\|_2\le\|v\|_1 \le C_2\|v\|_2$ for all $v\in V$. A normed $K$-vector space $(V, \|\,\,\,\|)$ is a ($K$-) {\it Banach space} if $V$ is complete with respect to $\|\,\,\,\|$. Recall that any finite-dimensional $K$-vector space admits a norm, any two norms are equivalent and it is complete. The $K$-vector space $\Cd(\cX, K)$ with the supremums norm $\|f\|_{\infty} = \sup_{\gamma\in \cX}\, |f(\gamma)|_p$ is a $K$-Banach space.
 
Let $V$ be a $K$-vector space. Recall that an $\cO_K$-submodule $L\subseteq V$ is a {\it lattice} if $\bigcup_{a\in K^*}\, aL = V$ and $\bigcap_{a\in K^*}\, aL = \{0\}$. For a given lattice $L\subseteq V$ the function $p_L(v) \colon \!\! = \inf_{v\in aL} |a|_p$ is a norm on $V$. If $\|\,\,\,\|$ is another norm then $p_L$ is equivalent to $\|\,\,\,\|$ if and only if $L$ is open and bounded in $(V,\|\,\,\,\|)$. A lattice $L\subseteq V$ is {\it complete} if $V$ is complete with respect to $p_L$. Finally a torsion free $\cO_K$-module $L$ is said to be {\it complete} if $L$ is a complete lattice in $L\otimes_{\cO_K} K$. For example the $\cO_K$-dual of a free module is a complete torsionfree $\cO_K$-module.

Let $(V, \|\,\,\,\|)$ be a Banach space. An element $\mu\in \Dist(\cX, V)$ is a {\it measure (or bounded distribution)} if $\mu$ is continuous with respect to the supremums norm, i.e.\ if there exists $C\in \bR$, $C>0$ such that $\|\int_{\cX}\, f\,\, d\mu\|\le  C \|f\|_{\infty}$ for all $f\in C^0(\cX, K)$. We will denote the space of $V$-valued measures on $\cX$ by $\Meas(\cX, V)$. If $L\subseteq V$ is an open and bounded lattice then $\Meas(\cX, V)$ is the image of the canonical inclusion $\Dist(\cX, L)\otimes_{\cO_K} K \to \Dist(\cX, V)$. An element $\mu\in \Meas(\cX, V)$ can be integrated not only against locally constant functions but against any $f\in \Cd(\cX, K)$. In fact since $C^0_c(\cX, K)$ is dense in the Banach space $(\Cd(\cX, K), \|\,\,\,\|_{\infty})$ the functional \eqref{distr1} extends uniquely to a continuous functional 
\begin{equation}
\label{meas1}
\Cd(\cX, K) \lra V,\,\,\, f\mapsto \int\, f\,\, d\mu. 
\end{equation}
If $\cX = \cX_1\times \cX_2$ then we obtain as a refinement of the bilinear map \eqref{fubini}  a pairing
\begin{equation}
\label{fubini2}
\Meas(\cX, V)\times \Cd(\cX_1, K) \lra \Meas(\cX_2, V),\, (\mu, f_1)\mapsto \int_{\cX_1}\, f_1\,d\mu.
\end{equation}

\subsection{$p$-adic measures}
\label{subsection:padicmeas} 
Given $\mu\in \Dist(\cX, \bC)$ we want to clarify what do we mean by saying that $\mu$ is a $p$-adic measure. For simplicity assume that $\cX$ is compact. The distribution $\mu$ extends to a $\bC_p$-linear map
\begin{equation}
\label{distr1a}
C^0(\cX, \bC_p) \lra \bC_p\otimes_{\barQ} \bC,\,\, f\mapsto \int\, f\,\, d\mu
\end{equation}
and we denote its image by $V_{\mu}$ so that we can view $\mu$ as an element of $\Dist(\cX, V_{\mu})$. It is called a {\it $p$-adic measure} if $V_{\mu}$ is a finitely generated $\bC_p$-vector space and if $\mu\in \Dist^b(\cX, V_{\mu})$. Equivalently, the image of $\mu$ (considered as a map $C^0(\cX, \bZ) \to \bC$) is contained in a finitely generated $\overline{\cO}$-module. So if $\mu\in \Dist(\cX, \bC)$ is a $p$-adic measure \eqref{distr1a} extends to continuous functional $C(\cX, \bC_p) \lra V_{\mu},\, f\mapsto \int\, f\,\, d\mu$. 

\section{Local distributions attached to ordinary representations}
\label{section:localmass}

\subsection{Gauss sums}

Throughout this chapter $F$ denotes a finite extension of $\bQ_p$, $\cO = \cO_F$ its ring of integers and $\fp$ the maximal ideal of $\cO$. We denote by  $U$ the group of units of $\cO$ and put $U^{(n)} = \{x\in U| \,\, x\equiv 1 \mod \fp^n\}$. Let $q$ denote the number of elements of $\cO/\fp$. We fix an (additive) character $\psi: F \to \barQ^*$ such that $\Ker(\psi) = \cO$ and a generator $\varpi$ of $\fp$. We denote by $|x|$ the modulus of $x\in F^*$ (i.e.\ $|\varpi| = q^{-1}$) and by $\ord = \ord_F$ the additive valuation (normalized by $\ord(\varpi) =1$). The normalized Haar measure on $F$ will be denoted by $dx$ (normalized by $\int_{\cO} \, dx = 1$). We put $d^{\times} x = (1-\frac{1}{q})^{-1}\frac{dx}{|x|}$ so that $\int_U \, d^{\times} x = 1$.
 
\begin{lemma}
\label{lemma:addint}
Let $X\subseteq \{x\in F^*\mid\, \ord(x) \le -2\}$ be a compact open subset such for all $a\in X$ there exists $n\in \bZ$, $1 \le n\le -\ord(a)-1$ such that $a U^{(n)} \subseteq X$. Then, 
\[
\int_{X}\, \psi(x)d^{\times} x = 0.
\]
\end{lemma}

{\em Proof.} It is enough to consider the case $X = a U^{(n)}$ with $1 \le n\le -\ord(a)-1$. Choose $b\in F^*$ with $\ord(b) + \ord(a) = -1$. Hence $\psi(ab) \ne 1$ and $\ord(b) \ge n$ and therefore
\begin{eqnarray*}
\int_{X}\, \psi(x)d^{\times} x & = & \int_{U^{(n)}}\, \psi(ax)d^{\times} x = \int_{U^{(n)}}\, \psi(a(1+b)x)d^{\times} x\\
& = &\int_{U^{(n)}}\, \psi(ax)\psi(abx)d^{\times} x.
\end{eqnarray*}
Since $\ord(abx-ab) = -1 + \ord(x-1) \ge n -1 \ge 0$, we have $\psi(abx) = \psi(ab)$  for all $x\in U^{(n)}$. It follows 
\[
\int_{X}\, \psi(x)d^{\times} x = \psi(ab)\int_{U^{(n)}}\, \psi(ax)d^{\times} x =  \psi(ab)\int_{X}\, \psi(x)d^{\times} x,
\]
hence $\int_{X}\, \psi(x)d^{\times} x= 0$.
\enddemo

Recall that the {\it conductor $\fc(\chi)$} of a quasicharacter $\chi: F^* \to \bC^*$ is the largest ideal $\fp^n$ of $\cO$ such that $U^{(n)}\subseteq \Ker(\chi)$.

\begin{lemma}
\label{lemma:nongauss2}
Let $\chi: F^* \to \bC^*$ be a quasicharacter of conductor $\fp^n, n\ge 1$ and let $a\in F^*$ with $\ord(a) \ne -n$. Then we have 
\[
 \int_{U}\, \psi(ax) \chi(x)d^{\times} x = 0.
\]
\end{lemma}

{\em Proof.}  1.~case $\ord(a)>-n$: Choose $b\in F^*$ with $\max(-\ord(a),0)\le\ord(b) <n$, $1+b\in U$ and $\chi(1+b) \ne 1$. Then,
\begin{eqnarray*}
\int_{U}\, \psi(ax) \chi(x)d^{\times} x & = & \int_{U}\, \psi(ax(1+b)) \chi(x(1+b))d^{\times} x = \\
& = & \chi(1+b)\int_{U}\, \psi(ax) \psi(abx) \chi(x)d^{\times} x\\
& = & \chi(1+b)\int_{U}\, \psi(ax) \chi(x)d^{\times} x
\end{eqnarray*}
hence $\int_{U}\, \psi(ax) \chi(x)d^{\times} x=0$.

\noi 2.~case $\ord(a)<-n$: By \ref{lemma:addint} above we have
\[
 \int_{U}\, \psi(ax) \chi(x)d^{\times} x = \sum_{bU^{(n)}\in U/U^{(n)}} \chi(b) \int_{abU^{(n)}}\, \psi(x)d^{\times} x = 0.
\]
\enddemo

We recall the definition of the Gauss sum of a quasicharacter (with respect to the fix choice of $\psi$).

\begin{definition}
\label{definition:gauss}
Let $\chi: F^* \to \bC^*$ be a quasicharacter with conductor $\fp^n$, $n\ge 0$ and $a\in F^*$ with $\ord(a) = -n$. We define the Gauss sum of $\chi$ by
\[
\tau(\chi) = \tau(\chi,\psi) = [U:U^{(n)}]\int_{aU}\, \psi(x) \chi(x)d^{\times} x.
\]
\end{definition}

For a quasicharacter $\chi:F^* \to \bC^*$ we define 
\begin{eqnarray}
\label{intsteinberg}
\int_{F^*} \chi(x) \psi(x) dx \colon \!\! & = & \lim_{n\to+ \infty} \int_{x\in F^*, -n \le \ord(x)\le n} \chi(x)\psi(x) dx.
\end{eqnarray}

\begin{lemma}
\label{lemma:localmeasure}
Let $\chi:F^* \to \bC^*$ be a quasicharacter with conductor $\fp^f$. Assume that $|\chi(\varpi)| < q$. Then the integral \eqref{intsteinberg} converges and we have
\[
\int_{F^*} \chi(x) \psi(x) dx \quad  = \quad \left \{\begin{array}{ll} 
\frac{1-\chi(\varpi)^{-1}}{1-\chi(\varpi) q^{-1}} & \mbox{if $f=0$;}\\
\tau(\chi) & \mbox{if $f>0$.} 
\end{array}\right.
\]
\end{lemma}

{\em Proof.} Firstly, we remark 
\begin{equation}
\label{nongauss1}
 \int_{U}\, \psi(ax) d^{\times} x = \left\{\begin{array}{ll} 1 & \mbox{if $\ord(a)\ge 0$;}\\
                            - \frac{1}{q-1} & \mbox{if $\ord(a)= -1$;}\\
                           0 & \mbox{if $\ord(a)\le -2$;}\\
\end{array}\right.
\end{equation} 
for all $a\in F^*$. Since $(1-1/q)d^{\times} x = \frac{dx}{|x|}$, we obtain
\begin{eqnarray*}
\int_{F^*} \chi(x) \psi(x) dx & = & \sum_{n=-\infty}^{\infty} (1- 1/q)q^{-n} \int_{\varpi^n U}\, \chi(x) \psi(x) d^{\times} x.
\end{eqnarray*}
If $f >0$ then by \ref{lemma:nongauss2} we have
\[
\int_{F^*} \chi(x) \psi(x) dx = (1- 1/q)q^{f} \int_{\varpi^{-f} U}\, \chi(x) \psi(x) d^{\times} x = \tau(\chi).
\]
On the other hand if $f=0$ then by \eqref{nongauss1} we get
\begin{eqnarray*}
\int_{F^*} \chi(x) \psi(x) dx & = & (1- 1/q)\left(-\frac{q}{(q-1)\chi(\varpi)} + \sum_{n=0}^{\infty} \, (\chi(\varpi) q^{-1})^n\right)\\
& = & \frac{1-\chi(\varpi)^{-1}}{1-\chi(\varpi) q^{-1}}.
\end{eqnarray*}
\enddemo

\subsection{Ordinary representations of $\PGL_2(F)$}
\label{subsection:ordinary}

We introduce more notation. Let $K = G(\cO)$. For an ideal $\fc\subset \cO$ let $K_0(\fc)\subseteq K$ denote the subgroup of matrices $A$ (modulo $Z$) which are upper triangular modulo $\fc$. 

Let $\pi: G(F) \to \GL(V)$ be an irreducible admissible infinite-dimensional representation (where $V$ is a $\bC$-vector space). Recall \cite{casselman} that there exists a largest ideal $\fc(\pi)$ -- the {\it conductor} of $\pi$ -- such that $V^{K_0(\fc)} = \{v\in V\mid \, \pi(k) v = v \,\,\forall\,k \in K_0(\fc)  \}\ne 0$. In this case $V^{K_0(\fc)}$ is one-dimensional. 

The representation $\pi$ is called {\it tamely ramified} if the conductor divides $\fp$. This holds if and only if $\pi= \pi(\chi^{-1}, \chi)$ for an unramified quasicharacter $\chi: F^* \to \bC^*$ (see e.g.\ \cite{bump}, Ch.\ IV). More precisely if the conductor is $\cO_F$, then $\pi$ is spherical hence a principal series representation $\pi(\chi^{-1}, \chi)$ where $\chi: F^* \to \bC^*$ is an unramified quasicharacter with $\chi^2 \ne |\cdot |$. If $\fc(\pi)=\fp$, then $\pi$ is a special representation $\pi(\chi^{-1}, \chi)$ where $\chi$ is unramified with $\chi^2 = |\cdot |$. 

\begin{definition}
\label{definition:ordinary} 
Assume that $\pi = \pi(\chi^{-1}, \chi)$ is tamely ramified. Then $\pi$ is called ordinary if either $\chi^2 = |\cdot |$ or if $\pi$ is spherical and tempered and if $\chi(\varpi)q^{1/2}$ is a $p$-adic unit (i.e.\ it lies in $\barcO^*$).
\end{definition}

Thus if $\pi = \pi(\chi^{-1}, \chi)$ is tamely ramified and if we put $\alpha \colon \!\! = \chi(\varpi)q^{1/2}\in \bC$ then $\pi$ is ordinary if either $\alpha = \pm 1$ or if $\alpha\in \barcO^*$ and $|\alpha| = q^{1/2}$. Note that $\alpha$ determines $\pi$ uniquely, i.e.\ there exists a one-to-one correspondence between the set (of isomorphism classes) of ordinary representations of $G(F)$ and the set $\{\alpha\in \barcO^*|\,\,\, \alpha = \pm 1\,\,\mbox{or}\,\,\, |\alpha | = q^{1/2}\}$. We will call an element of the latter set an {\it ordinary parameter}. We will denote the class corresponding to $\alpha$ by $\pi_{\alpha}$ and define $\chi_{\alpha}(x) \colon \!\! = \alpha^{\ord(x)}$ (thus $\pi_{\alpha} = \pi(\chi_{\alpha}^{-1}|\cdot|^{-1/2}, \chi_{\alpha}|\cdot|^{1/2})$). If $\alpha = \pm 1$ (resp.\ $\alpha \ne \pm 1$) then $\pi_{\alpha}$ is special (resp.\ spherical). If $\alpha = 1$ then $\pi_{\alpha}= \St$ is the Steinberg representation.

\subsection{Bruhat-Tits tree} 

In the next section we recall the well-known construction of models of the spherical principal series representations and special representations of $G(F)$ it terms of the Bruhat-Tits tree $\cT$ of $G(F)$ (see e.g.\ \cite{barthellivne}). In fact we will work over an arbitrary ring $R$ rather than $\bC$. Here we recall a few facts regarding the tree $\cT$ (see e.g.\ \cite{serre-tree}). The set of vertices $\cV= \cV(\cT)$ of $\cT$ is the set of homothety classes of lattices in $F^2$. For any two vertices $v_1,v_2$ of $\cT$ denote $d(v_1,v_2)$ the distance between $v_1$ and $v_2$. A vertex $v$ is even or odd if its distance to to the standard vertex $v_0 = [\cO^2]$ is even or odd. The set of even (resp.\ odd) vertices will be denoted by $\cV_{\even}$ (resp.\ $\cV_{\odd}$). The group $G(F)$ operates on $\cT$ by $g [L] = [\widetilde{g}L]$ (where $\widetilde{g}\in \GL_2(F)$ is a lift of $g\in G(F)$). Let $=\oE= \oE(\cT)$ (resp.\ $\cE= \cE(\cT)$) denote the set of oriented (resp.\ unoriented) edges of $\cT$. For $e\in \oE$ let $o(e)$ (resp.\ $t(e)$) denote the origin (resp.\ target) of $e$ and let $\bar{e}$ be the same edge as $e$ but with opposite orientation as $e$. Given $e\in \oE$ the set of ends $U(e)$ of $e$ is an open compact subset of $\bP^1(F)$. We recall its definition. For $(x,y)\in F^2$ let $\ell_{(x,y)}$ denote the linear form $F^2\to F,(x',y') \mapsto \det\left(\begin{matrix} x & x'\\
y & y'\end{matrix}\right)$. Given $P= [x:y] \in \bP^1(F)$ and representatives $L_1$ and $L_2$ of  $o(e)$ and $t(e)$ respectively with $\varpi L_1\subsetneq L_2\subsetneq L_1$ we have $P\in U(e)$ iff $\ell_{(x,y)}(L_2) \subsetneq \ell_{(x,y)}(L_1)$. For $e\in \oE(\cT)$ and $g\in G(F)$ we get $U(\bar{e}) = \bP^1(F)- U(e)$ and $g U(e) = U(g e)$ for all. For $n\in \bZ$ put $v_n = [\cO\oplus \fp^n]$. The set $\{v_n\mid n\in \bZ\}$ determines the standard apartment $A$ of $\cT$. The edge of $A$ with origin $v_{n+1}$ and target $v_n$ will be denoted by $e_n$. One easily checks that $U(e_n) = \fp^{-n}\subset \bP^1(F)$, so $\bigcap_n U(e_n)= \{0\}$ and $\bigcap_n U(\bar{e}_n)= \{\infty\}$ so the sequence $\{e_n\}_{n\in \bZ}$ is the geodesic from $\infty$ to $0$.

Following \cite{barthellivne} we define the height $h(v)\in \bZ$ of $v\in \cV$ as follows. The geodesic ray from $v$ to $\infty$ has a non-empty intersection with $A$. If $v_n$ is any point in the intersection we define $h(v) = n - d(v, v_n)$. It is independent of the choice of $v_n$ and satisfies $h(v_n) = n$. We need the following simple 

\begin{lemma}
\label{lemma:Phohe}
(a) For all $e\in \cE(\cT)$ we have
\begin{equation*}
\label{Phohe}
h(t(e)) = \left\{\begin{array}{ll} h(o(e)) + 1 & \mbox{if $\infty\in U(e)$,}\\
                            h(o(e)) - 1 & \mbox{otherwise.}
\end{array}\right.
\end{equation*}
\noi (b) For all $a, b\in F, a\ne 0$ and $v\in \cV(\cT)$ we have 
\[
h\left(\left(\begin{matrix} a & b\\
0 & 1\end{matrix}\right) v\right) = -\ord(a) + h(v).
\]
\end{lemma}

{\em Proof.} (a) follows immediately from the definition. For (b) it suffices to consider the case $v=v_0$ since the group $B(F)$ acts transitively on $\cV(\cT)$. Put $e = \left(\begin{matrix} a & b\\
0 & 1\end{matrix}\right) e_0$. Since $\infty\not\in b+a\cO = U(e)$ we obtain
\[
h\left(\left(\begin{matrix} a & b\\
0 & 1\end{matrix}\right) v_0\right) = h(t(e)) = h(o(e)) -1 = h\left(\left(\begin{matrix} a\varpi^{-1} & b\\
0 & 1\end{matrix}\right) v_0\right) -1.
\]
If $b\ne 0$ we have for $n\in \bN$ with $n \ge \ord(a)-\ord(b)$ and $m=\ord(a)$
\begin{eqnarray*}
h\left(\left(\begin{matrix} a & b\\
0 & 1\end{matrix}\right) v_0\right) \,\,= \,\,h\left(\left(\begin{matrix} a\varpi^{-n} & b\\
0 & 1\end{matrix}\right) v_0\right) - n
\,\, = \,\, h\left(\left(\begin{matrix} a\varpi^{-n} & 0\\
0 & 1\end{matrix}\right) v_0\right) - n\\
= \,\, h\left(\left(\begin{matrix} a & 0\\
0 & 1\end{matrix}\right) v_0\right) \,\,= \,\, h(v_{-m}) \,\, =\,\, -m.\hspace{2.8cm}
\end{eqnarray*}
\enddemo

\subsection{Representations of $\PGL_2(F)$ attached to the Bruhat-Tits tree}
\label{subsection:principalseries}

Let $R$ be an arbitrary ring. For an $R$-module $M$ let $C(\cV, M)$ denote the $R$-module of maps $\phi: \cV(\cT)\to M$ and $C(\oE, M)$ the $R$-module of maps $ c: \oE(\cT)\to M$. Moreover we denote by $C^{\pm}(\cE, M)\subseteq C(\oE, M)$ the submodule of $c\in C(\oE, M)$ with $c(\bar{e}) = \mp c(e)$ for all $e\in \oE(\cT)$. Both $C(\cV, M)$ and $C(\oE, M)$ are left $G(F)$-modules via $(g\phi)(v) \colon = \phi(g^{-1}v)$ and $(gc)(e) = c(g^{-1}e)$ and $C^+(\cE, M), C^-(\cE, M)\subseteq C(\oE, M)$ are $G(F)$-stable submodules. We let $C_c(\cV, M)\subseteq C(\cV, M)$ (resp.\ $C_c(\oE, M)\subseteq C^(\oE, M)$, $C_c^{\pm}(\cE, M)\subseteq C^{\pm}(\cE, M)$) be the submodule of $\phi\in C(\cV, M)$ with $\phi(v) = 0$ for almost all $v$ (resp.\ $c\in C^{\pm}(\cE, M)$ with  $c(e) = 0$ for almost all $e$). We define pairings
\begin{eqnarray}
& \langle\,\,\,, \,\,\,\rangle: C_c(\cV, R)\times C(\cV, M) \  M, \,\,\, \langle\phi_1, \phi_2\rangle \colon = \sum_{v\in \cV} \,\phi_1(v) \phi_2(v)\hspace{0.5cm}
\label{pairtree1}
\\
& \langle\,\,\,, \,\,\,\rangle: C_c^{\pm}(\cE, R)\times C^{\pm}(\cE, M) \to M, \,\,\, \langle c_1, c_2\rangle \colon = \sum_{e\in \cE} \, c_1(e) c_2(e)
\label{pairtree2}
\end{eqnarray}
(note that in the second pairing the summand $c_1(e) c_2(e)$ does not depend on the choice of  orientation of $e$). Define maps
\begin{eqnarray*}
\label{eqn:delta}
& \delta: C(\oE, M)\lra C(\cV, M), \,\,\, \delta(c)(v) \colon= \sum_{t(e) = v} c(e),\hspace{2.4cm}\\
& \delta_{\pm}^*: C(\cV, M) \lra C^{\pm}(\cE, M), \,\,\,\delta_{\pm}^*(\phi)(e) \colon = \phi(t(e))\mp\phi(o(e)).
\hspace{1cm}\label{delta*}
\end{eqnarray*}
They are adjoint with respect to \eqref{pairtree1}, \eqref{pairtree2}, i.e.\ we have $\langle\delta(c), \phi\rangle = \langle c, \delta_{\pm}^*(\phi)\rangle$ for all $c\in C_c^{\pm}(\cE, R), \phi\in C(\cV, M)$ (by abuse of notation we denote the restriction of $\delta$ to any submodule of $C(\oE, M)$ also by $\delta$ and similarly for $\delta_{\pm}^*$).

For the function $\tau_{\pm}\in C(\cV,R)$ defined by $\tau_{\pm}(v) = (\pm 1)^{d(v_0,v)}= (\pm 1)^{h(v)}$ for $v\in \cV$ we consider the map 
\begin{equation*}
\label{taupm}
\langle \wcdot, \tau_{\pm}\rangle: C_c(\cV, R) \,\lra\, R,\,\, \phi\,\mapsto\, \, \langle \phi, \tau_{\pm}\rangle= \sum_{v\in \cV_{\even}}\, \phi(v) \pm \sum_{v\in \cV_{\odd}}\, \phi(v).
\end{equation*}
One readily verifies that the sequence of $R$-modules
\begin{equation*}
\label{cokerdelta}
\begin{CD}
0 @>>> C_c^{\pm}(\cE,R)@> \delta >> C_c(\cV,R) @> \langle \wcdot, \tau_{\pm}\rangle >> R @>>> 0
\end{CD}
\end{equation*}
is exact. Dually we have an exact sequence 
\begin{equation*}
\label{kerdelta}
\begin{CD}
0 @>>> M @> m\mapsto \tau_{\pm} m >> C(\cV,M)@> \delta_{\pm}^* >> C^{\pm}(\cE,M) @>>> 0.
\end{CD}
\end{equation*}
In particular the restriction $\delta_{\pm}^* = (\delta_{\pm}^*)|_{C_c(\cV,R)}: C_c(\cV,R)\to C_c^{\pm}(\cE,R)$ is injective. 

The kernel of $\delta:C^+(\cE, M)\to C(\cV, M)$ is the set of harmonic cocyles $C_{\har}(\cT, M)$, i.e.\ the set of maps \(c: \oE(\cT) \to M\) such that \(c(\bar{e}) = -c(e)\) for all $e\in \oE(\cT)$ and \(\sum_{o(e) = v} c(e)  = 0\) for all $v\in \cV(\cT)$. We recall the relation between harmonic cocyles and boundary distributions on $\bP^1(F)$. The map 
\begin{equation*}
\label{steinbarg}
 \Coker(\delta_+^*:C_c(\cV, R)\to C_c^+(\cE, R)) \, \lra \, C^0(\bP^1(F), R)/R
\end{equation*}
given by $c\mapsto \sum_{e\in \cE(\cT)} \, c(e) 1_{U(e)}$ is an isomorphism (note that $c(e) 1_{U(e)} \equiv c(\bar{e}) 1_{U(\bar{e})}$ mod $R$ -- i.e.\ modulo constant functions). Thus \eqref{pairtree2} induces a pairing
\begin{equation}
\label{pairtree4}
C^0(\bP^1(F), R)/R \times C_{\har}(\cT, M) \to M, (f,c) \mapsto \int_{\bP^1(F)} f(P) \mu_c(dP)
\end{equation}
i.e.\ a map $C_{\har}(\cT, M)\to \Dist(\bP^1(F), M), c\mapsto \mu_c$ so that $\mu_c$ has total mass $=0$. For $f= 1_{U(e)}$, $e\in \oE(\cT)$ we have \(\int_{\bP^1(F)} 1_{U(e)}(P) \mu_c(dP) = c(e)\).

The Hecke operator $T: C(\cV,M)\to C(\cV,M)$ is defined by \((T\phi)(v) = \sum_{o(e) = v} \phi(t(e))\). By (\cite{barthellivne}, Thm.\ 10) the $R[T]$-module $C_c(\cV,R)$ is free. Thus for $a\in R$ the map $T-a: C_c(\cV,R)\to C_c(\cV,R)$ is injective. If $a\ne \pm (q+1)$ we define $\sB_a(F,R)$ to be the cokernel so that there exists a short exact sequence 
\begin{equation}
\label{heckeinduced}
\begin{CD}
0 @>>> C_c(\cV,R)@> T- a\id>> C_c(\cV,R) @>>> \sB_a(F,R) @>>> 0
\end{CD}
\end{equation}
of $R[G(F)]$-modules. Note that $\sB_a(F,R)$ is free as an $R$-module.

For $a=\pm (q+1)$ we have $\langle \wcdot, \tau_{\pm}\rangle \circ (T-a) =0$ since $\langle T\phi, \tau_{\pm}\rangle = \langle \phi, T \tau_{\pm}\rangle = a\langle \phi,  \tau_{\pm}\rangle$ for all $\phi\in  C_c(\cV,R)$. We put $\sB_a(F,R) = \Ker(\langle \wcdot, \tau_{\pm}\rangle)/\Image(T- a)$, so that the sequence 
\begin{equation}
\label{heckeharmonic2}
0 \to \sB_a(F,R) \to \Coker(T-a: C_c(\cV, R)\to  C_c(\cV, R))\stackrel{\langle \wcdot, \tau_{\pm}\rangle}{\lra}  R \to 0
\end{equation} 
is exact. It is easy to see that $\sB_a(F,R)$ is again an $R[G(F)]$-module which is free as an $R$-module. Since $\delta \circ \delta_{\pm}^* = (q+1)\id \pm T$ we see that $\delta$ induces a map $C_c^{\pm}(\cE,R)\to \sB_{\pm(q+1)}(F,R)$ such that
\begin{equation}
\label{heckeharmonic}
\begin{CD}
0 @>>> C_c(\cV,R) @> \delta_{\pm}^* >> C_c^{\pm}(\cE,R)@>>> \sB_{\pm(q+1)}(F,R) @>>> 0
\end{CD}
\end{equation} 
is an exact sequence of $R[G(F)]$-modules. 

Dually, for $a\in R$ we define $\sB^a(F, M)$ as follows. If $a\ne \pm(q+1)$ we let $\sB^a(F, M)$ be the kernel of $T-a: C(\cV,M)\to C(\cV,M)$. If $a=\pm (q+1)$ then $T-a$ maps the submodule $\tau_{\pm} M = \{\tau_{\pm} \cdot m\mid m\in M\}$ to zero so it induces an endomorphism of the quotient $C(\cV,M)/\tau_{\pm} M$ and we define $\sB^a(F, M)$ to be its kernel. 

Since $T: C_c(\cV,R)\to C_c(\cV,R)$ and $T: C(\cV,M)\to C(\cV,M)$ are adjoint with respect to \eqref{pairtree1} we obtain a pairing 
\begin{equation}
\label{pairtree5}
\langle\,\,\,, \,\,\,\rangle: \sB_a(F, R)\times \sB^a(F, M) \to  M
\end{equation}
which induces an isomorphism $\sB^a(F, M)\to \Hom_R(\sB_a(F, R), M)$. 

For $\ep = \pm 1\}$ let as before $\chi_{\ep}(x) = \ep^{\ord(x)}$. We can view $\chi_{\ep}$ as a character of $F^*/(F^*)^2$ and so $\chi_{\ep}(\det(g))$ is defined for $g\in G(F)$. 
Since $T(\tau_{\pm}\cdot \phi) = \pm T(\phi)\cdot \tau_{\pm}$ the isomorphism $C_c(\cV,R) \to C_c(\cV, R), \phi\mapsto \phi\cdot \tau_{\pm}$ induces an isomorphism
\begin{equation}
\label{twist1}
\Tw_{\ep}: \sB_a(F,R) \lra \sB_{\ep a}(F,R)
\end{equation}
which satisfies $\chi_{\ep}(\det(g)) \Tw_{\ep}(g\phi) = g\Tw_{\ep}(\phi)$ for all $\phi\in \sB_a(F,R)$ and $g\in G(F)$. If $\ep = + 1$ then $\Tw_{\ep}$ is of course the identity. In general we have $\Tw_{\ep}\circ \Tw_{\ep} = \id$. The operators \eqref{twist1} will be in section \ref{subsection:linv} in order to show that certain $\cL$-invariants do not change  under quadratic twists.

We want to reinterpret the sequences \eqref{heckeinduced} and \eqref{heckeharmonic} 
in terms of the induced representation $\Ind^{G(F)}_{K} R$ and $\Ind^{G(F)}_{K_0(\fp)} R$. Since $G(F)$ acts transitively on $\cV(\cT)$ and $\oE(\cT)$ and the stabilizer of $v_0$ and $e_0$ is $K$ and $K_0(\fp)$ respectively we have $C_c(\cV,R)\cong \Ind^{G(F)}_{K} R$ and $C(\oE,R)\cong \Ind^{G(F)}_{K_0(\fp)} R$. The element $W=\left(\begin{matrix} 0 & 1\\
\varpi & 0 \end{matrix}\right)\in G(F)$ normalizes $K_0(\fp)$ hence induces an involution
$W : \Ind^{G(F)}_{K_0(\fp)} R\lra \Ind^{G(F)}_{K_0(\fp)} R,\, W(\phi)(g) \, \colon\!\! =\phi(Wg)$
and under the above isomorphism $C^{\pm}(\cE,R)$ is mapped onto $(\Ind^{G(F)}_{K_0(\fp)} R)^{W=\mp 1}$. Hence we have exact sequences of $R[G(F)]$-modules 
\begin{equation}
\label{heckeinduced3}
\begin{CD}
0 \lra \Ind^{G(F)}_{K} R @> T- a\id>> \Ind^{G(F)}_{K} R @>>> \sB_a(F,R) \lra 0
\end{CD}
\end{equation}
for $a\in R$, $a\ne \pm (q+1)$ and 
\begin{equation}
\label{heckeinduced4}
0 \lra \Ind^{G(F)}_{K} R \lra (\Ind^{G(F)}_{K_0(\fp)} R)^{W=\mp 1} \lra \sB_{\pm (q+1)}(F,R) \lra 0.
\end{equation}

If $R=\bC$ and $a = \alpha +q/\alpha$ for some $\alpha\in \bC^*$, $\alpha \ne \pm 1$ (resp.\ $\alpha=\pm 1$) then it is well-known that $\sB_a(F,\bC)$ is a model of the principal series representation (resp.\ special representation) $\pi(\chi_{\alpha}^{-1}|\wcdot|^{-1/2}, \chi_{\alpha}|\wcdot|^{1/2})$ (see e.g.\ \cite{serre}, or \cite{barthellivne}). In particular $\sB_a(F,\bC)$ admits a (up to scalar) unique Whittaker functional, i.e.\ a nontrivial linear map $\lambda: \sB_a(F,\bC) \to \bC$ such that 
\[
\lambda\left(\left(\begin{matrix} 1 & x\\
0 & 1\end{matrix}\right)\phi\right) \,\, =\,\,\psi(x) \lambda(\phi)
\]
for all $x\in F$ and $\phi\in  \sB_a(F,\bC)$. This fact will be used in section \ref{subsection:localdistr}. 

\subsection{Distributions attached to elements of $\sB^a(F,M)$}
\label{subsection:diststeinberg}

Given $\rho\in C(\cV, R)$ define $R$-linear maps
\begin{eqnarray*}
\label{wdelta}
& \wdelta_{\rho}: C(\cE, M)\lra C(\cV, M), \, \, \wdelta_{\rho}(c)(v) \colon= \sum_{t(e) = v} \rho(o(e)) c(e),\hspace{2.2cm}\\
& \wdelta^{\rho}: C(\cV, M) \lra C^+(\cE, M), \,\wdelta^{\rho}(\phi)(e) \colon = \rho(o(e))\phi(t(e))- \rho(t(e))\phi(o(e)).
\end{eqnarray*}
They are adjoint with respect to \eqref{pairtree1}, \eqref{pairtree2}, i.e.\ we have $\langle \wdelta_{\rho}(c), \phi\rangle = \langle c, \wdelta^{\rho}(\phi)\rangle$ for all $c\in C_c^+(\cE, R), \phi\in C(\cV, M)$. For the constant function $\rho\equiv 1$ we have $\wdelta^1 = \delta_+^*, \wdelta_1 = \delta$. Note that for $\rho_1, \rho_2\in C(\cV,R)$ and $\phi\in C(\cV,M)$ we have 
\begin{equation}
\label{delta*delta}
(\wdelta_{\rho_1} \circ \wdelta^{\rho_2})(\phi) = T(\rho_1 \cdot \rho_2)\cdot \phi - \rho_2 \cdot T(\rho_1\cdot \phi).
\end{equation}
Hence for $a\in R$ and $\rho\in \sB^a(F,R)$ the maps $\wdelta_{\rho}$ and $\wdelta^{\rho}$ induce $R$-linear maps
\begin{eqnarray*}
&\wdelta_{\rho}: \,C^0(\bP^1(F), R)/R\cong \Coker(\delta_+^*:C_c(\cV, R)\to C_c^+(\cE, R)) \lra \sB_a(F,R),\\
&\wdelta^{\rho}: \,\sB^a(F,M) \lra \Ker(\delta: C^+(\cE, M)\to C(\cV, M)) = C_{\har}(\cT, M).\hspace{1cm}
\end{eqnarray*}
In fact by applying \eqref{delta*delta} to $\rho_1 =1$ and $\rho_2 = \rho$ (resp.\ $\rho_1 =\rho$ and $\rho_2 = 1$) we see that $\wdelta^{\rho}$ maps $\Ker(T-a)$ into $\Ker(\delta)$ (resp.\ $\wdelta_{\rho}$ induces a map $\Coker(\delta) \to \Coker(T-a)$). 

Let $\alpha\in R^*$ and put $a= \alpha + q/\alpha$ and $\chi_{\alpha}: F^* \to R^*, x\mapsto \chi_{\alpha}(x)=\alpha^{\ord(x)}$. In the following we assume that $q$ is not a zero-divisor in $R$ so that $\alpha = \pm 1$ if and only if $a = \pm (q+1)$. One easily checks that the function $\rho(v) \colon = \alpha^{h(v)}$ lies in $\sB^a(F,R)$. Put $\wdelta_{\alpha}: = \wdelta_{\rho}$ and $\wdelta^{\alpha}: = \wdelta^{\rho}$, so 
\begin{equation*}
\label{deltaalpha}
\wdelta_{\alpha}: C^0(\bP^1(F), R)/R \lra \sB_a(F, R), \quad
\wdelta^{\alpha}:  \sB^a(F, M)\lra C_{\har}(\cT, M)
\end{equation*}
are adjoint with respect to  \eqref{pairtree4} and \eqref{pairtree5}. 

\begin{lemma}
\label{lemma:wdelta}
(a) We have $\wdelta_{\alpha}(g f) \, =\, \chi_{\alpha}(a)^{-1} \,g\,\wdelta_{\alpha}(f)$
for all $g = \left(\begin{matrix} a & b\\
0 & 1\end{matrix}\right)\in B(F)$ and $f\in C^0(\bP^1(F), R)/R$.

\noi (b) If $\alpha=\pm 1$ then $\wdelta_{\alpha}: C^0(\bP^1(F), R)/R \lra \sB_a(F, R)$ is an isomorphism. 

\noi (c) If $\alpha \ne \pm 1$ then 
\[
\begin{CD}
0 \lra C^0(\bP^1(F), R)/R @> \wdelta_{\alpha} >> \sB_a(F, R) @> \phi\mapsto\langle \phi, \rho\rangle >> R \lra 0
\end{CD}
\]
is exact. 
\end{lemma} 

{\em Proof.} (a) follows immediately from $\rho\left(\left(\begin{matrix} a & b\\
0 & 1\end{matrix}\right) v\right) = \chi_{\alpha}(a)^{-1} \rho(v)$ by Lemma \ref{lemma:Phohe} and the simple proof of (b) will be left to the reader.

For (c) consider the commutative diagram
\[
\begin{CD}
0 @>>> C_c(\cV, R) @>\delta_+^*>> C_c^+(\cE, R) @>>> C^0(\bP^1(F), R)/R @>>> 0\\
@. @VV \eqref{rhovert} V @VV \wdelta_{\alpha} V@VV \wdelta_{\alpha} V\\
0 @>>> C_c(\cV,R)@> a\id - T >> C_c(\cV,R) @>>> \sB_a(F,R) @>>> 0
\end{CD}
\]
where the first vertical map is the isomorphism
\begin{equation}
C_c(\cV, R)\lra C_c(\cV, R),  \,\, \phi\mapsto (v\mapsto \phi(v) \rho(v)).
\label{rhovert}
\end{equation}
So it remains to prove that the upper row of the diagram 
\[
\begin{CD}
0 @>>> C^+_c(\cE, R) @> \wdelta_{\alpha} >> C_c(\cV, R) @> \phi\mapsto\langle \phi, \rho\rangle >> R @>>> 0\\
@. @VV \eqref{rhoedge} V @VV \eqref{rhovert} V@VV \id V\\
0 @>>> C^+_c(\cE, R)@> \delta >> C_c(\cV,R) @>  \phi\mapsto\langle \phi, \tau_+\rangle  >> R @>>> 0
\end{CD}
\]
is exact where \eqref{rhoedge} is the isomorphism
\begin{equation}
C^+_c(\cE, R)\lra C^+_c(\cE, R),  \,\, c\mapsto (e\mapsto c(e) \phi(t(e))\phi(o(e))). \label{rhoedge}
\end{equation}
However the lower row is exact.
\enddemo

We define $R$-linear maps
\begin{eqnarray}
\label{delta3}
& \delta_{\alpha}: C_c(F^*, R) \,\,\,  \lra \,\,\, \sB_a(F,R) & \mbox{if $\alpha \ne 1$,}\\
& \delta_{\alpha}: C_c(F, R) \,\,\,  \lra \,\,\, \sB_a(F,R) & \mbox{if $\alpha = 1$.}\nonumber
\end{eqnarray}
as follows. If $\alpha \ne 1$ and $f\in C_c(F^*, R)$ we define $\delta_{\alpha}(f)$ by extending $\chi_{\alpha}(x) f(x)$ by zero to $\bP^1(F)$ and then applying $\wdelta_{\alpha}$. If $\alpha = 1$ and $f\in C_c(F, R)$ we extend $f$ by zero to $\bP^1(F)$ and then apply $\wdelta_{\alpha}$. We let $F^*$ act on $C_c(F^*, R)$ and $C_c(F, R)$ by $(a\cdot f)(x) = f(a^{-1}x)$. It induces a $T(F)$-action via the isomorphism $T\cong\bG_m$. Also we define a $F^*$ operation on $\Dist(F, M)$ and $\Dist(F^*, M)$ by
\[
\int f(x) (t\mu)(dx) = t (\int (t^{-1}f)(x)\mu(dx)).
\]
for all $f\in C_c(F^*, R)$ resp.\ in $f\in C_c(F, R)$. The following result is an immediate consequence of Lemma \ref{lemma:wdelta} (a).

\begin{lemma}
\label{lemma:spechar}
The map $\delta_{\alpha}$ is $T(F)$-equivariant.
\end{lemma} 

Let $H$ be a subgroup of $G(F)$ and $M$ a $R[H]$-module (in the applications in chapter \ref{section:pmeashmf} both $H$ and $M$ will be of "global nature"). We define a $H$-action on $\sB^a(F, M)$ by requiring that $\langle \phi, h\cdot \lambda\rangle = h \cdot\langle h^{-1}\phi, \lambda\rangle$ for all $h\in H$, $\phi\in \sB_a(F, R)$ and $c\in \sB^a(F, M)$. By passing to duals we get $(T(F)\cap H)$-equivariant homomorphisms
\begin{eqnarray*}
\delta^{\alpha}:  \sB^a(F,M) \,\,\,  \lra \,\,\, \left \{\begin{array}{ll} \Dist(F^*, M) & \mbox{if $\alpha \ne 1$,}\\
\Dist(F, M) & \mbox{if $\alpha = 1$}
\end{array}\right.
\end{eqnarray*}
characterized by 
\begin{equation}
\label{delta4}
\langle \delta_{\alpha}(f), \lambda\rangle  \,\,\, =  \,\,\, \left \{\begin{array}{ll} \,\, \int_{F^*} \, f(x) \,\delta^{\alpha}(\lambda)(dx) & \mbox{if $\alpha \ne 1$,}\\
\,\, \int_{F} \, f(x) \, \delta^{\alpha}(\lambda)(dx) & \mbox{if $\alpha = 1$.}
\end{array}\right.
\end{equation}

\subsection{Local distributions}
\label{subsection:localdistr}

In this section we assume $R= \bC$. Let $\alpha\in \barcO^*$ be an ordinary parameter, i.e.\ $\alpha = \pm 1$ or $|\alpha | = q^{1/2}$. Define $\mu_{\alpha} \colon \!\! = \psi(x) \chi_{\alpha}(x) dx\in \Dist(F^*, \bC)$ (resp.\ $\in \Dist(F, \bC)$ if $\alpha =1$). We call it the {\it local distribution associated to $\pi_{\alpha}$} (the justification for this terminology will become apparent in section \ref{subsection:pmeashmf}). $\mu_{\alpha}$ is the image of a Whittaker functional under \eqref{delta4} (see Prop.\ \ref{prop:whittdelta} below).

\begin{prop}
\label{prop:localmeasurespherical}
Let $\chi:F^* \to \bC^*$ be a quasicharacter with conductor $\fp^f$. Assume that $|\chi(\varpi)| < q^{1/2}$. Then the integral $\int_{F^*} \chi(x) \mu_{\alpha}(dx)$ converges and we have
\[
\int_{F^*} \chi(x) \mu_{\alpha}(dx) \,\,= \,\, \tau(\chi) e(\alpha, \chi) L(\einhalb,\pi_{\alpha}\otimes\chi)
\]
where 
\[
e(\alpha, \chi) \,\, = \,\, \left \{\begin{array}{ll} 
(1- \alpha \chi(\varpi)^{-1})  & \mbox{if $f=0, \alpha = \pm 1$;}\\
\left(1-\frac{\chi(\varpi)}{\alpha}\right)\left(1-\frac{1}{\alpha\chi(\varpi)}\right) & \mbox{if $f=0, \alpha \ne \pm 1$;}
\\
\alpha^{-f} & \mbox{if $f>0$.} 
\end{array}\right.
\]
\end{prop} 

{\em Proof.} Recall that for the local $L$-factors we have $L(s, \pi_{\alpha}\otimes \chi) = 1$ if $f>0$ and
\[
L(s, \pi_{\alpha}\otimes \chi) = (1- \chi(\varpi)\alpha q^{-(s+1/2)})^{-1} 
\]
if $\alpha = \pm 1, f=0$ and 
\begin{eqnarray*}
L(s, \pi_{\alpha}\otimes \chi) & = & L(s, \chi\chi_{\alpha}^{-1} |\cdot|^{-1/2}) L(s, \chi\chi_{\alpha} |\cdot|^{1/2}) \\
& = & (1- \chi(\varpi)\alpha^{-1} q^{-(s-1/2)})^{-1}(1- \chi(\varpi)\alpha q^{-(s+1/2)})^{-1}
\end{eqnarray*}
if $\alpha \ne \pm 1, f=0$. Thus the assertion follows from Lemma \ref{lemma:localmeasure}. 
\enddemo

\begin{prop}
\label{prop:whittdelta} 
\noi (a) There exists a unique Whittaker functional $\lambda= \lambda_a$ for $\sB_a(F,\bC)$ such that $\delta^{\alpha}(\lambda_a) = \mu_{\alpha}$.

\noi (b) Let $\cW_{\alpha} = \cW(\pi_{\alpha})$ denote the Whittaker model of $\pi_{\alpha}$. If $\alpha \ne 1$ (resp.\  $\alpha = 1$) then for any $f\in C_c(F^*, \bC)$ (resp.\ $f\in C_c(F, \bC)$) there exists $W = W_f\in \cW_{\alpha}$ such that 
\[
\int_{F^*} \, (af)(x) \, \mu_{\alpha}(dx) = W\left(\begin{matrix} a & 0\\
0 & 1\end{matrix}\right)
\]
for all $a\in F^*$ (resp.\  $\int_{F} \, (af)(x) \, \mu_{\alpha}(dx) = W\left(\begin{matrix} a & 0\\
0 & 1\end{matrix}\right)$). 

\noi (c) Let $H$ be an open subgroup of $U$ and put $W_H = W_{1_H}$. Then, for any $f\in C^0_c(F^*, \bC)^H$ we have
\[
\int_{F^*} f(x) \mu_{\alpha}(dx) = [U:H] \int_{F^*} f(x) \, W_H\left(\begin{matrix} x & 0\\
0 & 1\end{matrix}\right) d^{\times} x.
\]
\end{prop}

{\em Proof.} (a) We let the (aditive) group $F$ act on the Schwartz space $C_c(F, \bC)$ as usual by $(x\cdot f)(y) \colon = f(y-x)$. Thus the functional 
\[
\Lambda: C_c(F, \bC)\,\lra \, \bC, \,\, f\mapsto \int_F f(x)\, \psi(x) \, dx
\]
satisfies $\Lambda(xf) = \psi(x) \Lambda(f)$ for all $x\in F$ and $f\in C_c(F, \bC)$.
Also we let $F$ act on $C^0(\bP^1(F), \bC)/\bC$ as $x \phi \colon \!\! = \left(\begin{matrix} 1 & x\\
0 & 1\end{matrix}\right) \phi$ so that
\begin{equation}
\label{steinbint}
C^0(\bP^1(F), \bC)/\bC\lra  C_c(F, \bC),\,\, \phi\mapsto f(x)\colon = \phi([x:1])-\phi(\infty)
\end{equation}
is an $F$-equivariant isomorphism. Thus the composite 
\begin{equation}
\label{steinwhitt}
\St(F,\bC)= C^0(\bP^1(F), \bC)/\bC\stackrel{\eqref{steinbint}}{\lra}  C_c(F, \bC)\stackrel{\Lambda}{\lra} \bC
\end{equation}
is a Whittaker functional of the Steinberg representation. It follows from Lemma \ref{lemma:wdelta} (a), (b) that for $\alpha = \pm 1$ the composition 
\[
\lambda: \sB_a(F,\bC) \stackrel{\delta_{\alpha}^{-1}}{\lra} C^0(\bP^1(F), \bC)/\bC\stackrel{\eqref{steinbint}}{\lra}C_c(F, \bC) \stackrel{\Lambda}{\lra}\bC
\]
is a Whittaker functional.

Assume now $\alpha \ne \pm 1$ and let $\lambda: \sB_a(F, \bC)\to \bC$ be a Whittaker functional. Since $\langle u \phi, \rho\rangle = \langle \phi, u \rho\rangle = \langle  \phi, \rho\rangle$ for all $\phi\in \sB_a(F, \bC)$ and $u=\left(\begin{matrix} 1 & x\\
0 & 1\end{matrix}\right)$ the map $\langle \wcdot, \rho\rangle: \sB_a(F, \bC)\lra \bC, \,\,\phi\mapsto\langle \phi, \rho\rangle$ is not a Whittaker functional. Therefore by Lemma \ref{lemma:wdelta} (a), (c) the map $\lambda \circ \wdelta_{\alpha}: \St(F,\bC)=C^0(\bP^1(F), \bC)/\bC\to \sB_a(F, \bC)\to \bC$ is a Whittaker functional of $\St(F,\bC)$ so -- after replacing $\lambda$ by a scalar multiple -- we may assume that $\lambda \circ \wdelta_{\alpha}$ is equal to the Whittaker functional \eqref{steinwhitt}.
Then 
$\delta^{\alpha}(\lambda)(f) = (\lambda \circ \wdelta_{\alpha})(\chi_{\alpha} \cdot f) = \Lambda(\chi_{\alpha} \cdot f) =  \int_{F^*} f(x)\, \chi_{\alpha}(x)\, \psi(x) \, dx =\mu_{\alpha}(f)$ for all $f\in C_c(F^*, \bC)$.

\noi (b) By (a) the function $W(g) \colon \!\! = \lambda( g \cdot\delta_{\alpha}(f))$ lies in $\cW_{\alpha}$ and we have 
\[
\int_{F^*} \, (af)(x) \, \mu_{\alpha}(dx) = \lambda( \delta_{\alpha}(af)) = \lambda\left(\left(\begin{matrix} a & 0\\
0 & 1\end{matrix}\right)\delta_{\alpha}(f)\right) = W\left(\begin{matrix} a & 0\\
0 & 1\end{matrix}\right).
\]
\noi (c) It is enough to consider the case $f = 1_{aH}$ for $a\in F^*$. Then
\begin{eqnarray*}
\int_{F^*} f(x) \mu_{\alpha}(dx) = \int_{F^*} (a1_H)(x) \mu_{\alpha}(dx) \,\, = \,\, W_H\left(\begin{matrix} a & 0\\
0 & 1\end{matrix}\right)\hspace{1cm}\\
=  m\int_{F^*} 1_H(x) \, W_H\left(\begin{matrix} ax & 0\\
0 & 1\end{matrix}\right) d^{\times} x 
= m \int_{F^*} f(x) \, W_H\left(\begin{matrix} x & 0\\
0 & 1\end{matrix}\right) d^{\times} x
\end{eqnarray*}
with $m= [U:H]$.\enddemo

\subsection{Extensions of the Steinberg representation}
\label{subsection:steinbergext}

In this section we assume that $R$ is a topological Hausdorff ring. We consider certain extensions of the $R[G(F)]$-module $\St(F, R) = C(\bP^1(F), R)/R$ associated to a continuous homomorphism $\ell$ from $F^*$ to the additive group of $R$ (for a related construction see \cite{breuil}, 2.1). Let 
\begin{equation*}
\label{projection}
\pi: G(F) \, \lra \, \bP^1(F), \,\,\, g= \left(\begin{matrix} a & b\\
c & d\end{matrix}\right) \mapsto g\infty = [a:c]
\end{equation*}
be the canonical $G(F)$-equivariant projection. Note that 
\begin{equation*}
\label{deltacont}
\delta:  \Cd(F, R)\,\lra\, \St(R), \,\, f\mapsto \delta(f)(P) = \left \{\begin{array}{ll} \, f(x) & \mbox{if $P=[x:1]$,}\\
\, 0 & \mbox{if $P=\infty$}
\end{array}\right.
\end{equation*}
is an isomorphism of $R[T(F)]$-modules (its inverse $\delta^{-1}$ is given by $\delta^{-1}(\phi)(x) = \break\phi([x:1])-\phi(\infty)$).
Define $\wsE(\ell)$ as the $R$-module of pairs $(f, y)\in C(G(F),R)\times\break R$ with
\begin{equation*}
\label{equation:model1a}
\phi\!\left( g \cdot \left(\begin{matrix} t_1 & u\\
0 & t_2\end{matrix}\right) \right) = \phi(g) + \ell(t_1/t_2) y
\end{equation*}
for all $t_1, t_2\in F^*$, $u\in F$ and $g\in G(F)$. We denote by $\wsE(\ell)_0\cong R$ the submodule consisting of pairs $(\phi,0)$ with $\phi: G(F) \to R$ constant and put $\sE(\ell) = \wsE(\ell)/\wsE(\ell)_0$. The left $G(F)$-action on $\wsE(\ell)$ given by $g\cdot (\phi(h), y) = (\phi(g^{-1}h), y)$ induces a $G(F)$-action on $\sE(\ell)$.

\begin{lemma}
\label{lemma:extsteinberg}
(a) Let $\epsilon : \sE(\ell) \to R$ be given by $\epsilon(\phi,y) = y$. Then the sequence of $R[G(F)]$-modules 
\begin{equation}
\label{steinext}
\begin{CD}
0 @>>> \St(F, R) @>\phi \mapsto (\phi\circ \pi, 0) >> \sE(\ell) @>\epsilon >> R @>>> 0
\end{CD}
\end{equation}
is exact.

\noi (b) Let $\delta^*: H^1(G(F), \St(F, R)) \to H^1(F^*, \Cd(F, R))$ be the homomorphism induced by the maps $\bG_m \cong T \subseteq G, x \mapsto \left(\begin{matrix} x & 0\\0 & 1\end{matrix}\right)$ and $\delta^{-1}$, let $[\sE(\ell)]$ denote the cohomology class of the extension $\eqref{steinext}$ and let $c_{\ell}\in H^1(F^*, \Cd(F, R))$ be the class of the cocyle
\begin{equation}
\label{cycle}
z_{\ell}(a) \colon \!\! = (1-a) (\ell\cdot 1_{\cO}) \colon \!\!= \left\{\begin{array}{ll} \ell(a) 1_{a\cO} + \ell \cdot 1_{\cO- a\cO} & \mbox{if $\ord(a) \ge 0$;}\\
\ell(a) 1_{a\cO} - \ell \cdot 1_{a\cO- \cO} & \mbox{if $\ord(a) < 0$.}
\end{array}\right.
\end{equation}
Then, $\delta^*([\sE(\ell)]) = 2 c_{\ell}$.

\noi (c) If $\ell = \ord_F: F \to \bZ \to R$ and the topology on $R$ is discrete then \eqref{steinext} is isomorphic to \eqref{heckeharmonic2}.
\end{lemma}

{\em Proof.} (a) It suffices to show that $\sE(\ell)\to R,(\phi, y) \mapsto y$ is surjective. Define 
\begin{equation*}
\label{equation:cocycle}
\phi_0\!\left(\left(\begin{matrix} a & b\\
c & d\end{matrix}\right)\right) \,\,\, = \,\,\, \left \{\begin{array}{ll} \ell\left(\frac{a^2}{ad-bc}\right) &  \mbox{if $\ord(a) < \ord(c)$,}\\
\ell\left(\frac{c^2}{ad-bc}\right) &  \mbox{$\ord(a) \ge \ord(c)$}\\
\end{array}\right.
\end{equation*}
so $\phi_0\!\left(\left(\begin{matrix} a & b\\
c & d\end{matrix}\right)\right) = \ell\left(\frac{a^2}{ad-bc}\right) + 2\ell(c/a)1_{c\cO}(a)$ if $a, c\ne 0$. One easily checks that $(\phi_0, 1)\in \sE(\ell)$.

\noi (b) Note that $\pi\!\left(\begin{matrix} x & -1\\
1 & 0\end{matrix}\right) = [x:1]$ and $\pi\left(\begin{matrix} 1 & 0\\
0 & 1\end{matrix}\right) = \infty$. Thus for $a\in F^*$ let $\phi\in \St(R)$ be given by $\phi\circ \pi \, =\, \left(\begin{matrix} a & 0\\
0 & 1\end{matrix}\right)\phi_0 - \phi_0$. Then for $x\in F$ we have
\begin{eqnarray*}
\delta^{-1}(\phi)(x) & = & \phi_0\!\left(\left(\begin{matrix} a & 0\\
0 & 1\end{matrix}\right)^{-1}\left(\begin{matrix} x & -1\\
1 & 0\end{matrix}\right)\right) - \phi_0\!\left(\left(\begin{matrix} a & 0\\
0 & 1\end{matrix}\right)^{-1}\right)\\
&& - \phi_0\!\left(\left(\begin{matrix} x & -1\\
1 & 0\end{matrix}\right)\right) + \phi_0\!\left(\left(\begin{matrix} 1 & 0\\
0 & 1\end{matrix}\right)\right)\\
& = & \ell(x^2/a) + 2\ell(a/x) 1_{a\cO}(x) - \ell(1/a) - (\ell(x^2) + 2\ell(1/x) 1_{\cO}(x))\\
& = & 2 \ell(x)(1_\cO(x)- 1_{a\cO}(x)) + 2\ell(a) 1_{a\cO}(x) \, =\, 2z_{\ell}(a)(x).
\end{eqnarray*}
\noi (c) Note that if $\ell = \ord$ then $\phi_0(k) =0$ for all $k\in K$. Hence for $g = k\cdot \left(\begin{matrix} t_1 & *\\
0 & t_2\end{matrix}\right)\in G(F) = KB(F)$ and $h \in K$ we have $(h \cdot \phi_0)(g) = \phi_0(h^{-1}k) + \ord(t_1/t_2) = \ord(t_1/t_2) = \phi_0(g)$, i.e.\ $(\phi_0, 1)$ is $K$-invariant. By Frobenius reciprocity we obtain a homomorphism $\Psi: C_c(\cV, R)\cong \Ind^{G(F)}_{K} R \to \sE(\ell)$. One can easily verify that $\Psi$ induces an isomorphism $\Coker(\delta\circ \delta^*_+) = \Coker(T-(q+1)\id) \cong \sE(\ell)$ and that the sequence $C^+_c(\cE, R) \stackrel{\delta}{\lra} C_c(\cV, R) \stackrel{\ep \circ \Psi}{\lra} R \to 0$ is exact, so the assertion follows.\enddemo

\subsection{Semi-local theory}
\label{subsection:semiloc}

We briefly discuss how to generalize some of the previous constructions to the semi-local case. Let $F_1, \ldots, F_m$ be finite extensions of $\bQ_p$ and let $q_i$ be the number of elements of the residue field of $F_i$ for $i=1, \ldots, m$. We put $F= F_1\times \ldots\times F_m$ and $F_S = \prod_{i\in S}\, F_i$ for a subset $S\subseteq \{1, \ldots, m\}$. Let $R$ be a ring and $a_1, \ldots, a_m\in R$ put $\ua = (a_1, \ldots, a_m)$ and define the $R[G(F)]$-module $\sB_{\ua}(F,R)$ as the tensor product of $\sB_{a_1}(F_1,R), \ldots, \sB_{a_m}(F_m,R)$
\begin{equation*}
\label{seimlocaltensor1}
\sB_{\ua}(F,R)\,\, =\,\, \bigotimes_R\, \sB_{a_i}(F_i,R).
\end{equation*}
To define the semi-local analogues of the maps \eqref{delta3} let $\alpha_1, \ldots, \alpha_m\in R^*$ and assume $a_i = \alpha_i + q_i/\alpha_i$ for $i=1, \ldots, m$. Let $S_1 = \{i\in \{1, \ldots, m\}|\, \alpha_i = 1\}$ and $S_2 = S_1^c \colon \!\! = \{1, \ldots, m\}-S_1$. It is easy to see that 
\begin{eqnarray}
\label{seimlocaltensor2}
\bigotimes_{i\in S_1}\, C_c^0(F_i, R) \otimes \bigotimes_{i\in S_2}\, C_c^0(F_i^*, R) & \lra & C^0_c(F_{S_1}\times F_{S_2}^*, R),\\
\bigotimes_{i\in S_1} f_i \otimes \bigotimes_{i\in S_2} f_i \hspace{1.5cm}& \mapsto &((g_i)_{i=1, \ldots, m} \mapsto \prod_{i=1}^m f_i(g_i))\nonumber
\end{eqnarray}
is an isomorphism.
We define the $R[T(F)]$-linear map 
\begin{equation}
\label{delta7}
\delta_{\ualpha}: C^0_c(F_{S_1}\times F_{S_2}^*, R) \,\,\,  \lra \,\,\, \sB_{\ua}(F,R)
\end{equation}
as the composite of the inverse of \eqref{seimlocaltensor2} and $\bigotimes_{i=1, \ldots, m} \delta_{\alpha_i}$.

For a $R$-module $M$ we define $\sB^{\ua}(F,M) \, = \, \Hom_R(\sB_{\ua}(F, R), M)$ and let 
\begin{equation}
\label{upairdual}
\langle\,\,\,, \,\,\,\rangle: \sB_{\ua}(F, R) \times \sB^{\ua}(F,M)\lra M
\end{equation}
be the evaluation pairing. If $H$ is a subgroup of $G(F)$ and $M$ a $H$-module then we define an $H$-action on $\sB^{\ua}(F,M)$ as before by $\langle \phi, h\cdot c\rangle = h \cdot\langle h^{-1}\phi, c\rangle$ for $h\in H$, $\phi\in \sB_{\ua}(F, R)$ and $c\in \sB^{\ua}(F, M)$.
By passing in \eqref{delta7} to duals we get a $(T(F)\cap H)$-linear map
\begin{equation*}
\label{udelta}
\delta^{\ualpha}:  \sB^{\ua}(F,M) \,\,\,  \lra \,\,\, \Dist(F_{S_1}\times F_{S_2}^*, M).
\end{equation*}
Note that $\sB^{\ua}(F,M) \cong \sB^{\ua}(F_S,\sB^{\ua}(F_{S^c},M))$ for any subset $S$ of $\{1, \ldots, m\}$. Note also that we have a canonical map $\sB^{\ua}(F,R) \otimes_R M\to \sB^{\ua}(F,M)$. In particular we get a map $\bigotimes_{i=1}^m \,  \sB^{a_i}(F_i,R) \to \sB^{a_m}(F_m,\bigotimes_{i=1}^{m-1} \,  \sB^{a_i}(F_i,R))$ and by iterating this construction we get a homomorphism
\begin{equation}
\label{prodsteinbergdual}
\bigotimes_{i=1}^m \,  \sB^{a_i}(F_i,R)   \,\,\lra    \,\,    \sB^{\ua}(F,R).
\end{equation}

Finally, we introduce the semi-local analogues of the maps \eqref{twist1}. Let $\chi: F^* \to R^*$ be an unramified quadratic homomorphism, i.e.\ for each $i=1, \ldots, m$ the restriction $\chi_i$ of $\chi$ to the factor $F_i$ is unramified and $\chi^2 =1$. Put $\uep = (\ep_1, \ldots, \ep_m) = (\chi_1(\varpi_1), \ldots, \chi_m(\varpi_m))\in \{\pm 1\}^m$ so that $\chi_i = \chi_{\ep_i}$ in the notation of section \ref{subsection:principalseries}. Define 
\begin{equation}
\label{twist2}
\Tw_{\chi} = \bigotimes_{i=1, \ldots, m} \Tw_{\ep_i}: \sB_{\ua}(F,R) \lra \sB_{\uep\ua}(F,R).
\end{equation}
Again $\chi(\det(g))\,\Tw_{\chi}(g \phi) \, =\, g\Tw_{\chi}(\phi)$ holds for all $g\in G(F)$ and $\phi\in \sB_{\ua}(F,R)$. Also for a subgroup $H$ of $G(F)$ and a $R[H]$-module $M$ the isomorphism \eqref{twist2} induces an isomorphism of $R$-modules
\begin{equation*}
\label{twist3}
\quad \Tw_{\chi}: \sB^{\uep\ua}(F,M) \lra \sB^{\ua}(F,M)
\end{equation*}
which is adjoint to \eqref{twist2} with respect to the pairing \eqref{upairdual}
 and satisfies $h\Tw_{\chi}(c)\break = \chi(\det(h))\Tw_{\chi}(h c)$ for all $h\in H$ and $c\in \sB^{\ua}(F,M)$.

\section{Special zeros of $p$-adic $L$-functions}
\label{section:abstractezc}

\paragraph{{\it Notation}} We introduce the following notation which will be used throughout the rest of this paper. $F$ denotes a totally real number field of degree $d+1$ over $\bQ$ with ring of integers $\cO_F$. For a non-zero ideal $\fa\subseteq \cO_F$ we set $N(\fa) = \sharp(\cO_F/\fa)$. We denote by $\bfP_F$ the set of all places of $F$ and by $\bfP_F^{\infty}$ (resp.\ $S_{\infty}$) the subset of finite (resp.\ infinite) places. For a prime number $\ell$, we shall write $S_{\ell}$ for the set of places above $\ell$. We denote by $\sigma_0, \ldots, \sigma_d$ the different embeddings of $F$ into $\bR$ and let $\infty_0, \ldots, \infty_d$ be the corresponding archimedian places of $F$.  Elements of $\bfP_F$ will be denoted by $v,w$ or also by $\fp, \fq$ if they are finite. If $\fp\in \bfP_F^{\infty}$, we denote the corresponding prime ideal of $\cO_F$ also by $\fp$. For $v\in \bfP_F$,  we denote by $F_v$ the completion of $F$ at $v$. If $v$ is finite then $\cO_v$ denotes the valuation ring of $F_v$ and $\ord_v$ the corresponding normalized (additive) valuation on $F_v$ (so $\ord_v(\varpi) =1$ if $\varpi\in \cO_v$ is a local uniformizer at $v$). Also for $v\in \bfP_F$ we let $|\wcdot|_v$ be the associated normalize multiplicative valuation on $F_v$. Thus if $v\in S_{\infty}$ corresponds to the embedding $\sigma:F \to \bR$ then $|x|_v = |\sigma(x)|$ and if $v=\fq$ is finite then $|x|_{\fq} = N(\fq)^{-\ord_{\fq}(x)}$. For $v\in \bfP_F$ we put $U_v = \bR_+^*$ if $v$ is infinite and $U_v = \cO_v^*$ if $v$ is finite. Moreover if $v =\fp$ is finite and $n\ge 0$, then we also put $U_v^{(n)} = \{ x\in U_v|\, \ord_v(x-1)\ge n\}$. 

Let $\bA= \bA_F$ be the adele ring of $F$ and $\bI= \bI_F$ the group of ideles. Let $|\wcdot |: \bI_F \lra \bR^*$ be the absolute modulus, i.e.\ $|(x_v)_v | = \prod_v\, |x_v|_v$ for $(x_v)_v\in \bI_F$. For a finite subset $S\subseteq \bfP_F$ we let $\bA^S$ (resp.\ $\bI^S$) denote the $S$-adeles (resp.\ $S$-ideles) and put $F_S = \prod_{v\in S} F_v$. We also define $U^S =\prod_{v\not\in S} U_v$ and $U_S = \prod_{v\in S} U_v$. For $T\subseteq \bfP_{\bQ} = \{ 2,3,5, \ldots, \infty\}$ and $S=\{v\in \bfP_F\mid\, v|_{\bQ}\in T\}$ we often write $F_T$, $\bA^T$, $\bI^T$ etc.\ for $F_S$, $\bA^S$, $\bI^S$ etc. We also write $U^p$, $U_p$, $U^{p, S}$, $U^{p,\infty}$ etc.\ for $U^{\{p\}}$, $U_{\{p\}}$, $U^{S_p\cup S}$, $U^{S_p\cup S_{\infty}}$ etc. and use a similar notation for adeles and ideles. Thus for example for a finite subset $S$ of $\bfP_F^{\infty}$, $\bI^{S,\infty}$ denotes the set of $S\cup S_{\infty}$-ideles and for $\ell\in \bfP_{\bQ}$ we have $F_{\ell} = F\otimes \bQ_{\ell} = \prod_{v\in S_{\ell}} F_v$.

We fix an (additive) character $\psi: \bA\to \bC^*$ which is trivial on $F$. For $v\in \bfP_F$ let $\psi_v$ denote the restriction of $\psi$ to $F_v\hookrightarrow \bA$. For convenience we choose $\psi$ so that $\Ker(\psi_{\fp}) = \cO_{\fp}$ for all $\fp\in S_p$. Let $dx$ (resp.\ $dx_v$) denote the associated self-dual Haar measure on $\bA$ (resp.\ on $F_v$). Thus $dx = \prod_v dx_v$. For $v\in\bfP_F$ we define a normalized Haar measure $dx_v^{\times}$ on $F_v^*$ by $dx_v^{\times} = m_v\frac{dx_v}{|x_v|_v}$ where $m_v = (1-\frac{1}{N(v)})^{-1}$ if $v\in\bfP_F^{\infty}$ and $m_v=1$ if $v\in S_{\infty}$. For a character $\chi: \bI/F^* \to \bC^*$ and $v\in \bfP_F$ we denote by $\chi_v$ its $v$-component, i.e.\ $\chi_v : F_v^*\hookrightarrow \bI \stackrel{\chi}{\lra} \bC^*$. The {\it Gauss sum} $\tau(\chi)= \tau(\chi, \psi)$ of $\chi$ is then defined as $\tau(\chi) = \prod_{\fp\mid \ff(\chi)}\, \tau(\chi_{\fp})$.

We denote by $F^*_+$ the totally positive elements of $F$ and by $G(F)^+$ (resp.\ $G(F_{\infty})^+$) the subgroup of $G(F)$ (resp.\ $G(F_{\infty})$) of elements with totally positive determinant (note that $G(F)^+\cong \SL_2(F)/\{\pm 1\}$). The subgroups $B(F)^+\subseteq B(F)$ and $T(F)^+\subseteq T(F)$ are defined similarly. Furthermore we define subgroups $K_{\infty}^+\subseteq K_{\infty}\subseteq G(F_{\infty})$ as the image of $\Orth(2)^{\Hom(F,\bR)}\subseteq \GL_2(F_{\infty})$ and $\SOrth(2)^{\Hom(F,\bR)}$ under the projection $\GL_2(F_{\infty}) \to G(F_{\infty})$ (thus $K_{\infty}^+= K_{\infty}^+\cap G(F_{\infty})^+$).

There is a canonical  $G(F_{\infty})^+$-action on $\bH^{d+1}$ where $\bH \colon \!\! = \{z\in \bC\mid \, \Image(z) >0\}$; the embeddings $\sigma_0, \ldots, \sigma_d$ allow us to identify $G(F_{\infty})^+$ with $(G(\bR)^+)^{d+1}$ and the latter group acts on $\bH^{d+1}$ through linear transformations factor-by-factor. For $g = (g_0, \ldots, g_d)\in G(F_{\infty})^{+}$ and $\uz = (z_0, \ldots, z_d)\in \bH^{d+1}$ we define $j(g, \uz) = \prod_{\nu=0}^d j(g_\nu, z_\nu)$ where $j(\gamma, z) = \det(\gamma)^{-1/2} (cz+d)$ if $\gamma = \left(\begin{matrix} a & b\\ c & d\end{matrix}\right)\in G(\bR)^+$, $z\in \bH$.

Let $\fn$ be a non-zero ideal of $\cO_F$. For $v\in \bfP_F^{\infty}$ we put $K_0(\fn)_v = \{A \in G(\cO_v)|\,\, A\equiv \left(\begin{matrix} * & *\\ 0 & *\end{matrix}\right)\mod \fn\cO_v\}$ and set $K_0(\fn) = \prod_{v\in \bfP_F^{\infty}} \, K_0(\fn)_v$. If $S \subseteq \bfP_F^{\infty}$ we also put $K_0(\fn)^S = \prod_{v\in \bfP_F^{\infty}-S} \, K_0(\fn)_v$.

\subsection{Rings of functions on ideles and adeles}

\paragraph{\bf The module $C_c^{\flat}(F_v,  K)$}

Let $v$ be a finite place of $F$ and let $K$ be a Hausdorff topological field (in the application $v$ will be a place above $p$ and $K$ a $p$-adic field). We identify $C_c(F_v^*, K)$ with the submodule $\{f\in C_c(F_v, K)\mid\, \mbox{$f\equiv 0$ near $0$}\}$ of $C_c(F_v, K)$ and define
\begin{equation*}
\label{flat}
C^{\flat}_c(F_v, K)\,\,\, = \,\,\, C^0_c(F_v, K) + C_c(F^*_v, K).
\end{equation*}
Both $C_c(F_v^*, K)$ and  $C^{\flat}_c(F_v, K)$ are $F^*_v$-submodules of $C_c(F_v, K)$. For $f\in C_c(F_v^*, K)^{U_v}$ and $x\in F_v^*$ the infinite sum 
\begin{equation*}
\label{unitint1}
(\sum_{n=0}^{\infty}\, \varpi^n f)(x) \colon \!\! = \sum_{n=0}^{\infty}\, f(\varpi^{-n} x)
\end{equation*}
is finite and one easily checks that $F_v^* \to K, \, x\mapsto (\sum_{n=0}^{\infty}\, \varpi^n f)(x)$ extends to a function in $C^0_c(F_v, K)$ which will be denoted by $(1-\varpi)^{-1}f$. For example if $f= 1_{U_v}$ then $(1-\varpi)^{-1}f = 1_{\cO_v}$. Thus we obtain a $F_v^*$-equivariant $K$-linear monomorphism 
\begin{equation}
\label{unitint2}
C_c(F_v^*, K)^{U_v} \lra C^0_c(F_v, K), \,\,\, f\mapsto (1-\varpi)^{-1}f
\end{equation}
Its image is $C^0_c(F_v, K)^{U_v}$. Hence if we consider the following two-step filtration $\cF_v^{\bu}$ on $C^{\flat}_c(F_v, K)$
\begin{equation}
\label{unitint3}
\cF_v^0 = C^{\flat}_c(F_v, K), \quad \cF_v^1 =  C^0_c(F_v, K)^{U_v}, \quad \cF_v^{2} =0
\end{equation}
then we have for the associated graded $F_v^*$-modules $\gr_{\cF_v}^n = \cF_v^n/\cF_v^{n+1}$ 
\begin{equation}
\label{unitint4}
\gr_{\cF_v}^n\,\,\, \cong \,\,\, \left \{\begin{array}{ll} 
C_c(F_v^*,K)/C_c(F_v^*,K)^{U_v} & \mbox{if $n=0$,}\\
C_c(F_v^*, K)^{U_v} & \mbox{if $n=1$,}\\
0 & \mbox{otherwise.}
\end{array}\right.
\end{equation}
Note also that 
\begin{eqnarray}
\label{induced}
C_c(F_v^*,K)^{U_v} \,\, \cong \,\, C_c(F_v^*/U_v,K) \,\, \cong \,\, \Ind^{F_v^*}_{U_v} K,\\ \nonumber
C_c(F_v^*, K)/C_c(F_v^*, K)^{U_v} \,\, \cong\,\,  \Ind^{F_v^*}_{U_v} (C_c(U_v^*, K)/K).
\end{eqnarray}

\paragraph{\bf The module $\sC_c^{\flat}(S_1,S_2,  K)$}

Consider now two (possibly empty) disjoint subsets $S_1, S_2$ of $S_p$ and let $R$ be a topological Hausdorff ring. We define
\begin{eqnarray*}
\sC(S_1, S_2, R)  & = & C(F_{S_1}\times F_{S_2}^* \times \bI^{p,\infty}/ U^{p,\infty},R),\\
\sCd(S_1, S_2, R)   & = &\Cd(F_{S_1}\times F_{S_2}^* \times \bI^{p,\infty}/ U^{p,\infty},R),\\
\sC_c(S_1, S_2, R)   & = & C_c(F_{S_1}\times F_{S_2}^* \times \bI^{p,\infty}/ U^{p,\infty},R),\\
\sC^0(S_1, S_2, R)  & = & C^0(F_{S_1}\times F_{S_2}^* \times \bI^{p,\infty}/ U^{p,\infty},R).
\end{eqnarray*}
We have 
\[
\sC_c^0(S_1, S_2, R)\subseteq \sC_c(S_1, S_2, R)\subseteq \sCd(S_1, S_2, R) \subseteq \sC(S_1, S_2, R)
\]
(for the first inclusion see Prop.\ \ref{prop:freermod2} below). Note that if $R$ carries the discrete topology the first three rings are all equal.

Assume now that $S= S_p$ so $S_1\overset{\cdot}{\cup} S_2 = S_p$ and that $K=R$ is a field. We define the submodule $\sC_c^{\flat}(S_1,S_2, K)$ of $\sC_c(S_1,S_2, K)$ as the image of the embedding
\begin{equation*}
\label{nicefunc2}
\bigotimes_{v\in S_1} C^{\flat}_c(F_v, K) \otimes C_c(\bI^{S_1,\infty}/ U^{S,\infty}, K) \lra \sC_c(S_1, S_2, K).
\end{equation*}
We have 
\[
\sC_c^0(S_1, S_2, K) \subseteq \sC_c^{\flat}(S_1, S_2, K) \subseteq \sC_c(S_1, S_2, K)\subseteq \sC(S_1, S_2, K).
\]
The filtrations \eqref{unitint3} on $C^{\flat}_c(F_v, K)$ for all $v\in S_1$ induce a filtration $\cF^{\bu}$ on $\sC_c^{\flat}(S_1,S_2,  K)$. For $\un= (n_v)_{v\in S_1}\in \bZ^{S_1}$ put $|\un| = \sum_v n_v$. Then $\cF^m C^{\flat}_c(F_v, K)$ is defined as the image of 
\begin{equation*}
\label{unitint5}
\bigoplus_{\un \in \bZ^{S_1}, |\un| =m}\, \bigotimes_{v\in S_1} \cF_v^{n_v} C^{\flat}_c(F_v, K) \otimes C_c(\bI^{S_1,\infty}/ U^{S,\infty}, K) \lra \sC_c(S_1, S_2, K).
\end{equation*}
We get for the associated graded quotients $\gr_{\cF}^m = \cF^m/\cF^{m+1}$ 
\begin{equation}
\label{unitint6}
\gr_{\cF}^m \,=\, \bigoplus_{|\un| =m}\, \bigotimes_{v\in S_1} \gr_{\cF_v}^{n_v} \otimes C_c(\bI^{S_1,\infty}/ U^{S,\infty}, K).
\end{equation}
We fix a splitting of the exact sequence $1\to E_+ \to F^*_+\to \Gamma\colon \!\! =F^*_+/E_+ \to 1$, i.e.\ we fix a subgroup $\cT\subseteq F^*_+$ such that $F^*_+ = E_+ \times \cT$. 

\begin{prop}
\label{prop:freermod2} 
$\sC_c^{\flat}(S_1, S_2, K)$ is a free $K[\cT]$-module.
\end{prop}

{\em Proof.} It is enough to prove that each graded quotient  $\gr_{\cF}^m$ and therefore each summand in \eqref{unitint6} is a free $K[\cT]$-module. Since
\begin{equation*}
\label{induced2}
C_c(\bI^{S_1,\infty}/ U^{S,\infty}, K) \,\, \cong \,\, \Ind^{\bI^{S_1,\infty}}_{U^{S_1,\infty}} \, C_c(U^{S_1,\infty}/ U^{S,\infty}, K)
\end{equation*}
we deduce using \eqref{unitint4} and \eqref{induced} that each summand in \eqref{unitint6} is isomorphic to a $K[\bI^{\infty}]$-module of the form $\Ind^{\bI^{\infty}}_{U^{\infty}} V$ for some $K[U^{\infty}]$-module $V$. Hence it is free a $K[\cT]$-module. Indeed, since by assumption $\cT\cap U^{\infty} = 1$ we have $\Ind^{\bI^{\infty}}_{U^{\infty}} V \,\, = \,\, \bigoplus_{i=1}^h\, \Ind^{\cT}_1 x_i V$ (as $\cT$-modules) where $\{x_1,\ldots, x_h\}$ is a system of representatives of $\bI^{\infty}/U^{\infty} \cT= \bI^{\infty}/U^{\infty} F^*_+$ (cf.\ \cite{brown}, Prop.\ 5.6, p.\ 69).\enddemo

\subsection{Computation of $\partial((\log_p \circ \cN)^k)$ for $k=0,\ldots, r$}
\label{subsection:capprod}

\paragraph{\bf Definition of $\partial$} Assume again that $S_1\overset{\cdot}{\cup} S_2 = S_p$ and that $R$ is a topological Hausdorff ring. Let $\cG_p= \Gal(M/F)$ is the Galois group of the maximal abelian extension $M/F$ which is unramified outside $p$ and $\infty$. We shall now construct a canonical homomorphism
\begin{equation}
\label{recipr2}
\partial: C(\cG_p, R) \lra  H_d(F^*_+, \sC_c(S_1, S_2, R)).
\end{equation}
Let $E_+$ be the group of totally positive units of $\cO_F$. Firstly, there exists an isomorphism
\begin{equation}
\label{recipr1}
C(\cG_p, R) \lra H_0(F^*_+/E_+, H^0(E_+, \sC_c(S_p,R)))
\end{equation}
defined as follows. Let $\barE_+$ be the closure of $E_+$ in $U_p$ and let $\pr: \bI^{\infty}/ U^{p,\infty} \to \bI^{\infty}/ (\barE_+ \times U^{p,\infty})$ denote the projection. The map
\begin{equation*}
C_c(\bI^{\infty}/ (\barE_+ \times U^{p,\infty}), R)\lra H^0(E_+, C_c(\bI^{\infty}/ U^{p,\infty},R)),\,\, f\mapsto f\circ \pr
\end{equation*}
is an isomorphism. Hence its inverse induces an isomorphism
\begin{equation}
\label{recipr2a}
H_0(F^*_+/E_+,\! H^0(E_+, C_c(\bI^{\infty}/ U^{p,\infty},\! R)))\!\cong \! H_0(F^*_+/E_+,\!C_c(\bI^{\infty}/ (\barE_+ \times U^{p,\infty}), \!R)).
\end{equation}
The reciprocity map of class field theory $\rho: \bI/F^* \to \cG_p$ induces a surjection $\barho:\bI^{\infty}/(\barE_+ \times U^{p,\infty})\to \cG_p$ whose kernel is discrete and $\cong F^*_+/E_+$. It follows that the map
\begin{equation}
\label{recipr2b}
\rho^{\sharp}: H_0(F^*_+/E_+, C_c(\bI^{\infty}/ (\barE_+ \times U^{p,\infty}), R))\lra C(\cG_p, R)
\end{equation}
defined by $\rho^{\sharp}([f])(\barho(x)) = \sum_{\zeta\in F^*_+/E_+} \, f(\zeta x)$ is an isomorphism as well. The map \eqref{recipr1} is the composite of \eqref{recipr2a} with \eqref{recipr2b}. 

Let $A$ be any $F^*_+$-module. Next we construct a homomorphism
\begin{equation}
\label{edge1}
H_0(F^*_+/E_+, H^0(E_+, A))\lra  H_d(F^*_+, A)
\end{equation}
Since $E_+ \cong \bZ^d$ we have $H_d(E_+,\bZ)\cong \bZ$. Choose a generator $\eta$ of $H_d(E_+,\bZ)$. Since the action of $F^*_+/E_+$ on $H_d(E_+,\bZ)$ is trivial, taking the cap product with $\eta$ yields an $F^*_+/E_+$-equivariant map $H^0(E_+, A) \to H_d(E_+, A)$ hence
\begin{equation}
\label{capeta}
H_0(F^*_+/E_+, H^0(E_+, A)) \lra  H_0(F^*_+/E_+, H_d(E_+, A))
\end{equation}
We define \eqref{edge1} as the composite of \eqref{capeta} with the edge morphism 
\begin{equation}
\label{edge2}
H_0(F^*_+/E_+, H_d(E_+, A))\to H_d(F^*_+, A)
\end{equation}
of the Hochschild-Serre spectral sequence. 

There is in fact a canonical choice for $\eta$. Consider the action of $E_+$ on $\bR^{d+1}_0= \{(x_0, \ldots, x_d) \in \bR^{d+1} |\, \sum_{i=0}^d \, x_i =0\}$ given by $a\cdot (x_0, \ldots, x_d) = (\log(\sigma_0(a)) +x_0, \ldots,\log(\sigma_d(a)) +x_d)$. The $d$-dimensional manifold $\bR^{d+1}_0/E_+$ is oriented and compact. We chose $\eta\in H_d(E_+,\bZ)$ so that it corresponds  to the fundamental class under the canonical isomorphism $H_d(E_+,\bZ)\cong H_d(\bR^{d+1}_0/E_+,\bZ)$ (thus $\eta$ depends on our chosen ordering of the real places of $F$).

Finally we define \eqref{recipr2} is the composite of \eqref{recipr1}, \eqref{edge1} (for $M = \sC_c(S_p, R)$) and the map $H_d(F^*_+, \sC_c(S_1, S_2, R))\lra H_d(F^*_+, \sC_c(S_1, S_2, R))$ induced by the inclusion $\sC_c(S_p, R)\subseteq \sC_c(S_1, S_2, R)$.

\paragraph{\bf Fundamental homology classes.}

We put $r = \sharp(S_1)$, $m=\sharp(S_p)$ and order the places above $p$, so that $S_1 = \{\fp_1, \ldots, \fp_r\}$ and $S_2=\{\fp_{r+1}, \ldots, \fp_{m}\}$. Beside $\eta\in H_{d-1}(E_+, \bZ)$ we consider two more canonical homology classes $\vartheta$ and $\varrho$. To begin with we introduce the following $F^*_+$-action on $\bR^{d+1}$, $\bR^r$ and  $\bI^{S_1,\infty}/U^{S_1,\infty}$ 
\begin{eqnarray*}
a\cdot (x_0, \ldots, x_d) & \colon \!\! = &   (\log(\sigma_1(a)) +x_1, \ldots,\log(\sigma_d(a)) +x_d) ,\\
a\cdot (y_1, \ldots, y_r)  & \colon \!\! = & (\ord_{\fp_1}(a) + y_1, \ldots,\ord_{\fp_r}(a) + y_r), \\
a\cdot (x_v)_{v\not\in S_1 \cup S_{\infty}} & \colon \!\! = &  (a x_v)_{v\not\in S_1 \cup S_{\infty}}.
\end{eqnarray*}
Let  $M$ be the submanifold of $\bR^{d+1}\times \bR^r\times \bI^{S_1,\infty}/U^{S_1,\infty}$ defined by the equation 
\begin{eqnarray*}
 \sum_{i=1}^d \, x_i \,- \, \left(\sum_{j=1}^r\, \log(N(\fp_j)) y_j\right) + \left(\sum_{v\not\in S_1 \cup S_{\infty}} \, \log(|x_v|_v)\right) \,= \, 0.
\end{eqnarray*}
and put $M_1 \colon \!\! = \bR^r\times \bI^{S_1,\infty}/U^{S_1,\infty}$. We have 
\[
H_0(M, \bZ)\cong H_0(M_1, \bZ) \cong C_c(\bI^{S_1,\infty}/ U^{S_1,\infty}, \bZ)
\] 
where the first isomorphism is induced by the projection $M\to M_1$. The group $F^*_+$ (resp.\ $\Gamma \colon \!\! =F^*_+/E_+$) acts properly discontinuously on $M$ (resp.\ on $M_1$) and the projection $\pi: M/F^*_+\to M_1/\Gamma$ is a fiber bundle with fiber $\cong \bR^{d+1}/E_+$ (in fact it is easy to see that it is trivial i.e.\ it is homeomorphic to the trivial bundle $M_1/\Gamma \times \bR^{d+1}/E_+$ over $M_1/\Gamma$). The base $M_1/\Gamma$ is a compact oriented $r$-dimensional manifold. 

\noi {\bf Definition of $\vartheta$.} Define $\vartheta$ as the image of the fundamental class under the composition
\begin{eqnarray*}
& H_{d+r}(M/F^*_+, \bZ)\! \cong \! H_{d+r}(F^*_+,H_0(M, \bZ))\!\cong \! H_{d+r}(F^*_+,C_c(\bI^{S_1,\infty}/ U^{S_1,\infty}, \bZ))\nonumber\\
\label{fundclass1}
&\lra H_{d+r}(F^*_+, C_c(\bI^{S_1,\infty}/ U^{p,\infty}, \bZ)) = H_{d+r}(F^*_+, \sC_c(\emptyset, S_2, \bZ))\end{eqnarray*}
where the last map is induced by the projection $\bI^{S_1,\infty}/ U^{p,\infty}\to \bI^{S_1,\infty}/ U^{S_1,\infty}$. If $R$ is arbitrary topological Hausdorff ring, then -- by abuse of notation -- we denote the image of $\vartheta$ under the canonical map $H_{d+r}(F^*_+, \sC_c(\emptyset, S_2, \bZ)) \to H_{d+r}(F^*_+, \sC_c(\emptyset, S_2, \bZ))$ also by $\vartheta$.

\noi {\bf Definition of $\varrho$.} Let $\cT$ be any subgroup of $F^*_+$ such that $\cT \cap E_+ = \{1\}$ and $\cT E_+$ has finite index in $F^*_+$ (we are mainly interested in the case $F^*_+ = E_+ \times \cT$) so that the group $\cT$ acts properly discontinuously on $M_1$. Let $\varrho_{\cT}\in H_r(\cT,C_c(\bI^{S_1,\infty}/ U^{p,\infty}, \bZ))$ be the image of the fundamental class of the oriented $r$-dimensional compact manifold $M_1/\cT$ under the canonical the map 
\begin{eqnarray}
\label{fundclass}
& \quad H_r(M_1/\cT, \bZ)) \cong H_r(\cT,H_0(M_1, \bZ))\cong H_{d+r}(\cT,C_c(\bI^{S_1,\infty}/ U^{S_1,\infty}, \bZ)\\
&\lra H_r(\cT, C_c(\bI^{S_1,\infty}/ U^{p,\infty}, \bZ))\nonumber
\end{eqnarray}

\begin{remarks} 
\label{remarks:varrho}
\rm (a) If $\cT$ and $\cT'$ are subgroups as above with $\cT' \subseteq \cT$ then we have $\res(\varrho_{\cT}) = \varrho_{\cT'}$. 

\noi (b) Let $\cT_1 = \{x\in \cT|\, \ord_{\fq}(x) = 0 \,\, \forall \, \fq\not\in S_1\}$, let $\cT_2$ be a subgroup of $\cT$ with $\cT = \cT_1\times \cT_2$ and let $\cF_2$ be a fundamental domain for the action of $\cT_2$ on $\bI^{S_1,\infty}/ U^{p,\infty}$ such that $U_{S_2}\cF = \cF$. Then $C_c(\bI^{S_1,\infty}/ U^{p,\infty}, \bZ)\cong \Ind^{\cT}_{\cT_1} C(\cF_2, \bZ)$ hence by Shapiro's Lemma
\begin{equation}
\label{shapiro}
H_{\bcdot}(\cT,C_c(\bI^{S_1,\infty}/ U^{p,\infty}, \bZ)) \,\, \cong \,\, H_{\bcdot}(\cT_1, C(\cF_2, \bZ))
\end{equation}
Let $\varrho_1\in H_r(\cT_1, \bZ)\cong H_r(\bR^r/\cT_1, \bZ)$ be the fundamental class of $\bR^r/\cT_1$. Then $\varrho_1\otimes 1_{\cF_2}$ is mapped to $\varrho_{\cT}$ under 
\[
H_r(\cT_1, \bZ)\otimes C(\cF_2, \bZ)^{\cT_1} \stackrel{\cap}{\lra} H_r(\cT_1, C(\cF_2, \bZ)) \stackrel{\eqref{shapiro}}{\cong} H_r(\cT,C_c(\bI^{S_1,\infty}/ U^{p,\infty}, \bZ)).
\]
\end{remarks}

For a subgroup $\cT\subseteq F^*_+$ such that $F^*_+ = E_+ \times \cT$ we shall explain the relation between the homology classes $\varrho_{\cT}$, $\eta$ and $\vartheta$. Consider the Hochschild-Serre spectral sequence
\begin{equation*}
\label{ss2}
E_{pq}^2 = H_p(E_+,H_q(\cT,C_c(\bI^{S_1,\infty}/ U^{p,\infty}, \bZ)))\Rightarrow H_{p+q}(F^*_+, C_c(\bI^{S_1,\infty}/ U^{p,\infty}, \bZ)).
\end{equation*}
Here we have $E_{pq}^2=0$ if $p>d$ or $q>r$ (the latter follows from \eqref{shapiro} above since $\cT_1$ is free-abelian of rank $r$). Thus we get an isomorphism $E_{d-1, r}^2 \cong E_{d+r-1}$ i.e.\
\begin{equation}
\label{hsedge2}
H_d(E_+,H_r(\cT,C_c(\bI^{S_1,\infty}/ U^{S_1,\infty}, \bZ))) \cong H_{d+r}(F^*_+, C_c(\bI^{S_1,\infty}/ U^{S_1,\infty}, \bZ)).
\end{equation}

\begin{lemma}
\label{lemma:fundclass} 
$\varrho_{\cT}$ is mapped to $\vartheta$ under the composite
\begin{eqnarray}
\label{coreta} 
& H_r(\cT,C_c(\bI^{S_1,\infty}/ U^{p,\infty}, \bZ))^{E_+} \stackrel{\cap \eta}{\lra} H_d(E_+,H_r(\cT,C_c(\bI^{S_1,\infty}/ U^{p,\infty}, \bZ)))\nonumber\\ & \stackrel{\eqref{hsedge2})}{\cong}\, H_{d+r}(F^*_+, C_c(\bI^{S_1,\infty}/ U^{p,\infty}, \bZ)) \nonumber
\end{eqnarray}
\end{lemma}

Recall the definition of cohomology classes defined in \ref{lemma:extsteinberg} (b).

\begin{definition}
\label{definition:logclass}
Let $\fp\in S_1$, let $R$ be a topological Hausdorff ring and let $\ell: F_{\fp}^* \to R$ be a continuous homomorphism. We denote by $c_{\ell}\in H^1(F_{\fp}^*, C_c(F_{\fp}, R))$ the cohomology class of the 1-cocyle \eqref{cycle} (i.e.\ of the cocycle $z_{\ell}(a) \colon \!\! = (1-a) (\ell\cdot 1_{\cO_{\fp}})$ for $a\in F^*_+$).
\end{definition}

By abuse of notation we shall write $c_{\ell}$ instead of $\res(c_{\ell})\in H^1(H, C_c(F_{\fp}, \break R))$ for any subgroup $H$ of $F_{\fp}^*$. We are interested in the case $H= F_+^*$, $R = \bC_p$ and either $\ell = \ord_{\fp}$ or $\ell = \log_p\circ \Norm_{F_{\fp}/\bQ_p}$ and will derive a formula for $(c_{\ell_{\fp_1}} \cup \ldots \cup c_{\ell_{\fp_r}}) \cap \vartheta$ in both cases.

We begin with the first case. Let $H \colon \!\! = \{x\in F^*_+|\, \ord_{\fp}(x) = 0 \,\,\forall\,\,\fp\in S_1\}$, $H_1\colon \!\! = \{x\in F^*_+|\, \ord_{\fp}(x) = 0 \,\,\forall\,\,\fp\in S_p\}$ and let $\cF_1$ denote a fundamental domain for the action of $H_1/E_+$ on $\bI^{p,\infty}/U^{p,\infty}$. Put $\cX \colon \!\! = \prod_{\fp\in S_1} \cO_{\fp} \times  \prod_{\fp\in S_2} \cO_{\fp}^* \times \cF_1\subseteq F_{S_1}\times F_{S_2}^* \times \bI^{p,\infty}/ U^{p,\infty}$. The characteristic function of $\cX$ clearly lies in $H^0(E_+, \sC^0_c(S_1,S_2, \bZ))$, hence defines an element $[1_{\cX}]\in H_0(F^*_+/E_+, H^0(E_+, \sC^0_c(S_1,S_2, \bZ)))$. 

\begin{prop}
\label{prop:thetaord}
For $\fp\in S_p$ put $c_{\fp} = c_{\ord_{\fp}}\in H^1(F^*_+, C^0_c(F_{\fp}, \bZ))$. We have
\[
\ep([1_{\cX}]) \,\, =\,\, (-1)^{\binom{r}{2}}\, (c_{\fp_1} \cup \ldots \cup c_{\fp_r}) \cap \vartheta
\]
Here $\ep$ denotes the map \eqref{edge1} for $A= \sC^0_c(S_1,S_2,\bZ)$.
\end{prop}

{\em Proof.} Similar as above we denote by $\varrho_1\in H_r(F^*_+/H, \bZ)$
the homology class which corresponds to the fundamental class of
$\bR^r/F^*_+$ under the natural isomorphism $H_r(F^*_+/H, \bZ) \cong H_r(\bR^r/F^*_+,\bZ)$. By taking the cap product with $\eta$ we can identify $C_c(\bI^{S_1,\infty}/  U^{S_1,\infty}, \bZ)^{E_+}$ with $H_d(E_+, C_c(\bI^{S_1,\infty}/ U^{S_1,\infty},\break \bZ))$. Note that $\cF\colon \!\! = \{1\}\times \cF_1\subseteq F_{S_2}^*/U_{S_2}\times \bI^{p,\infty}/ U^{p,\infty} = \bI^{S_1,\infty}/ U^{S_1,\infty}$ is a fundamental domain for the action of $\barH = H/E_+$. Hence if $\bD \colon \!\! = C_c(\bI^{S_1,\infty}/ U^{S_1,\infty}, \bZ))$ then $H_0(\barH, \bD) = H_0(\barH, \Ind^{\barH} \,C(\cF, \bZ)) \cong C(\cF, \bZ)$ and $H_q(\barH, \bD) =0$ for $q>0$. Consider the Hochschild-Serre spectral sequence
\begin{eqnarray*}
\label{ss1}
E_{pq}^2 = H_p(F^*_+/E_+,H_q(E^+,\bD)) \Rightarrow E_{p+q} = H_{p+q}(F^*_+, \bD).
\end{eqnarray*}
We have $E_{pq}^2=0$ if $q>d$ and $E_{pd}^2\cong H_p(F^*_+/E_+, \bD) \cong H_p(F^*_+/H, H_0(\barH, \bD))\break =0$ if $p>r$. It follows $E_{d+r} \cong E_{rd}^2$. Define
\begin{eqnarray}
\label{hsedge}
&& H_r(F^*_+/E_+, \bD) \stackrel{\cap\eta}{\lra}  H_r(F^*_+/E_+,H_d(E^+,\bD)) \cong H_{d+r}(F^*_+,\bD) \\
&& \hspace{2cm}\stackrel{\eqref{fundclass}}{\lra} H_{d+r}(F^*_+,\sC_c(\emptyset, S_2, \bZ)).\nonumber
\end{eqnarray}
It is easy to see that $\varrho$ is mapped to $\vartheta$ under \eqref{hsedge} and that $[1_{\cF}]\otimes \varrho_1$ is mapped to $\varrho$ under 
\begin{eqnarray*}
&& H^0(F^*_+/H, H_0(\barH, \bD))\otimes H_r(F^*_+/H, \bZ) \stackrel{\cap}{\lra}H_r(F^*_+/H, H_0(\barH, \bD))  \\
&& \hspace{2cm}\cong H_r(F^*_+/E_+, \bD).
\end{eqnarray*}
Note that we can view $c_{\fp} = c_{\ord_{\fp}}$ as an element of $H^1(F^*_+/H, C^0_c(F_{\fp}, \bZ)^H)$. The assertion thus follows from 
\begin{equation*}
\label{keyequalord1}
(-1)^{\binom{r}{2}}\, (c_{\fp_1} \cup \ldots \cup c_{\fp_r}) \cap \varrho_1 = [1_{\cO_{S_1}}]\in H_0(F^*_+/H, C^0_c(F_{S_1}, \bZ)^H)
\end{equation*}
where $\cO_{S_1} = \prod_{\fp\in S_1} \cO_{\fp}$. For that put $z_{\fp}=  z_{\ord_{\fp}}$ for $\fp \in S_1$ and choose generators $t_1, \ldots, t_r\in F^*_+/H$ such that $\ord_{\fp_i}(t_i) = 1$ and $\ord_{\fp_j}(t_i) = 0$ for all $j\ne i$, $1\le j \le r$. Note that $z_{\ord_{\fp_i}}(t_i) = t_i 1_{\cO_{\fp_i}}$ and  $z_{\ord_{\fp_i}}(t_j) = 0$ for $j\ne i$. Since the fundamental class of $\bR^r/\langle t_1,\ldots, t_r\rangle$ is the cross product of the fundamental classes of $\bR/\langle t_1 \rangle, \ldots, \bR/\langle t_r \rangle$ the $r$-cycle $\sum_{\sigma\in S_r} \,\sign(\sigma)\, [t_{\sigma(1)}| \ldots| t_{\sigma(r)}]$ is a representative of $\varrho_1$ (see \cite{maclane}, Ch.\ VIII, 8.8). Hence 
\begin{eqnarray*}
\label{keyequalord2}
&&\sum_{\sigma\in S_r} \,\sign(\sigma)\, z_{\fp_1}(t_{\sigma(1)}) \otimes t_{\sigma(1)} z_{\fp_2}(t_{\sigma(2)})\otimes \ldots \otimes t_{\sigma(1)}\ldots t_{\sigma(r-1)}z_{\fp_r}(t_{\sigma(r)})\\
&& = \,\,z_{\fp_1}(t_1) \otimes \ldots \otimes z_{\fp_r}(t_r) \,\,=\,\,\prod_{i=1}^r t_i \wcdot  1_{\cO_{S_1}}
\end{eqnarray*}
is a representative of $(-1)^{\binom{r}{2}} (c_{\fp_1} \cup \ldots \cup c_{\fp_r}) \cap \varrho_1$. 
\enddemo

We consider \eqref{recipr2} for $R=\bC_p$, i.e.\ $\partial: C(\cG_p, \bC_p) \lra  H_d(F^*_+, \sC_c(S_1, S_2, \bC_p))$. Let $\cN: \cG_p \to \bZ_p^*$ be defined by $\gamma \zeta = \zeta^{\cN(\gamma)}$ for all $p$-power roots of unity $\zeta$.

\begin{prop}
\label{prop:thetalog}
For $\fp\in S_1$ put $\ell_{\fp} \colon \!\! = \log_p \circ \Norm_{F_{\fp}/\bQ_p}: F_{\fp}^* \to \bC_p$. We have

\noi (a) $\partial((\log_p \circ \cN)^k) = 0$ for all $k=0,1, \ldots, r-1$.

\noi (b) $\partial((\log_p \circ \cN)^r) \,\, =\,\,  (-1)^{\binom{r}{2}}(c_{\ell_1} \cup \ldots \cup c_{\ell_r}) \cap \vartheta $.
\end{prop}

{\em Proof.} We choose again a subgroup $\cT\subseteq F^*_+$ such that $F^*_+ = E_+ \times \cT$. We denote by $\ell: \bI \to \bQ_p$ the composite 
\begin{equation*}
\label{log}
\ell:\bI \stackrel{\rho}{\lra} \cG_p \stackrel{\cN}{\lra} \bZ_p^* \stackrel{\log_p}{\lra} \bQ_p
\end{equation*}
and for a place $v$ of $F$ let $\ell_v: F_v\hookrightarrow \bI_F \stackrel{\ell}{\lra} \bQ_p$ be the $v$-component of $\ell$. Note that for $x = (x_v)\in \bI_F$ we have $\ell_v(x_v) = 0$ for almost all $v$ and 
 \begin{equation*}
 \ell(x) \,\, = \sum_{\fq} \, \ell_{\fq}(x_v).
 \end{equation*}
Let $\cF \subseteq \bI^{\infty}/U^{p,\infty}$ be a compact open fundamental domain for the action of $\cT$ such that $U_p\cF = \cF$. The function $(\log_p \circ \cN)^k$ is mapped under the inverse of \eqref{recipr1} to the class of $\ell^k 1_{\cF}$. 

Since $z_{\ell_{\fp}}(a) = (1-a) (\ell\cdot 1_{\cO_{\fp}})$ is a 1-cocyle with values in $C^{\flat}_c(F_{\fq}, \bC_p)$ we can view $c_{\ell_{\fp}}$ as an element of $H^1(F^*_+, C^{\flat}_c(F_{\fq}, \bC_p))$. Therefore the right hand side of (b) can be viewed as an element of $H_d(F^*_+, \sC^{\flat}_c(S_1,S_2, \bC_p))$. Note also that $\ell^k 1_{\cF}$ lies in $H^0(E_+,\sC_c^{\flat}(S_1, S_2, \bC_p))$. Therefore it suffices to show that the class $[\ell^k 1_{\cF}]\in H_0(\Gamma, H^0(E_+,\sC_c^{\flat}(S_1, S_2, \bC_p)))$ is mapped to $0$ (resp.\ to $(-1)^{\binom{r}{2}}(c_{\ell_1} \cup \ldots \cup c_{\ell_r})$) und
\begin{equation*}
\label{edge7}
H_0(\Gamma, H^0(E_+, \sC_c^{\flat}(S_1, S_2, \bC_p)))\stackrel{\eqref{edge1}}{\lra} H_{d-1}(F^*_+, \sC_c^{\flat}(S_1, S_2, \bC_p))
\end{equation*}
for $k=0, \ldots, r-1$ (resp.\ for $k=r$). 

After this preliminary remark we prove (a). Consider the commutative diagram
\begin{equation*}
\label{cap2}
\xymatrix@-0.5pc{H_0(\cT, H^0(E_+, \sC_c^{\flat}(S_1, S_2, \bC_p)))\ar[r]\ar[d]^{\eqref{edge1}}& H^0(E_+, H_0(\cT, \sC_c^{\flat}(S_1, S_2, \bC_p)))\ar[d]^{\cap\eta}\\
H_{d-1}(F^*_+, \sC_c^{\flat}(S_1, S_2, \bC_p))\ar[r]^{\coinf\hspace{0.5cm}} & H_{d-1}(E_+, H_0(\cT, \sC_c^{\flat}(S_1, S_2, \bC_p)))}
\end{equation*}
where the upper horizontal arrow is the canonical map induced by the inclusion $
H^0(E_+, \sC_c^{\flat}(S_1, S_2, \bC_p))\hookrightarrow \sC_c^{\flat}(S_1,
S_2, \bC_p)$. By Prop.\ \ref{prop:freermod2} the coinflation $H_{\bcdot}(F^*_+, 
\sC_c^{\flat}(S_1, S_2, \bC_p))\to H_{\bcdot}(E_+, H_0(\cT, \sC_c^{\flat}(S_1, S_2, \bC_p)))$ is an isomorphism. Hence it remains to prove that the image of $[\ell^k 1_{\cF}]$ under the upper horizontal map vanishes, i.e.~we have 
\begin{equation}
\label{ellvanish}
\ell^k 1_{\cF} \in I(\cT)\sC_c^{\flat}(S_1, S_2, \bC_p) \qquad \mbox{for all $k= 0,1,\ldots,r-1$}
\end{equation}
where $I(\cT)\subseteq \bC_p[\cT]$ denote the augmentation ideal.

We may shrink $\cT$. In fact if $\cT'\subseteq \cT$ is a subgroup of finite index then 
it follows from Prop.\ \ref{prop:freermod2} that $\res: H_0(\cT, \sC_c^{\flat}(S_1, S_2, \bC_p))\to H_0(\cT', \sC_c^{\flat}(S_1, S_2, \bC_p))$ is injective and if $\cF' \subseteq \bI^{\infty}/U^{p,\infty}$ is a fundamental domain for the action of $\cT'$ then we have $\res([\ell^k 1_{\cF}]) = [\ell^k 1_{\cF'}]$.

Hence we may assume that 
\begin{equation}
\label{tdecomposition}
\cT = \cT_p \times \cT^p \qquad \mbox{and} \qquad \cT_p= \langle t_1,\ldots, t_m\rangle
\end{equation}
with $\ord_{\fp_i}(t_i)>0$ and $\ord_{\fq}(t_i) = 0 = \ord_{\fp_i}(t)$ for all $i\in \{1, \ldots, m\}$, $t\in \cT^p$ and all finite places $\fq\ne \fp_i$. Put $F_i\colon \!\! = F_{\fp_i}$, $\cO_i\colon \!\! = \cO_{\fp_i}$ and $\cF_i \colon \!\! = \cO_i - t_i\cO_i$ for $i=1, \ldots, m$. Let $\cF^p\subseteq \bI^{p,\infty}/ U^{p,\infty}$ be a fundamental domain for the action of $\cT^p$. Then $\cF \colon \!\! = \prod_{i=1}^m \cF_i \times \cF^p\subseteq \bI^{\infty}/U^{p,\infty}$ is a fundamental domain for the $\cT$-action. 

We also denote by $\ell_p$ resp.\ $\ell^p$ the map $\ell_p: \bI_F \stackrel{\pr}{\lra} F_p \hookrightarrow \bI_F \stackrel{\ell}{\lra} \bQ_p$ resp.\ $\ell^p: \bI_F \stackrel{\pr}{\lra} \bI_F^p \hookrightarrow \bI_F \stackrel{\ell}{\lra} \bQ_p$ so that $\ell = \ell_p + \ell^p$ and $\ell_p = \sum_{i=1}^m \ell_i$  where $\ell_i = \ell|_{F_{\fp_i}}$ for $i=1, \ldots, m$.
We will show 
\begin{equation}
\label{ellvanish2}
\ell_p^k 1_{\cF} \in I(\cT_p)\sC_c^{\flat}(S_1, S_2, \bC_p) \qquad \mbox{for all $k= 0,1,\ldots,r-1$}
\end{equation}
where $I(\cT_p)\subseteq \bC_p[\cT_p]$ denote the augmentation ideal. Since $t \,\ell^p = \ell^p$ for all $t\in \cT_p$ this implies 
\begin{equation}
\label{ellvanish3}\ell_p^k \,(\ell^p)^j \,1_{\cF} \in I(\cT_p)\sC_c^{\flat}(S_1,  S_2, \bC_p)\subseteq I(\cT)\sC_c^{\flat}(S_1,  S_2, \bC_p)
\end{equation}
for all $k,j\ge 0$ with $k\le r-1$ and therefore \eqref{ellvanish}.

For $\Xi \subseteq \{1, \ldots, r\}$ we set 
\[
\cF_{\Xi} \colon \!\! = \prod_{i \in \Xi} \cO_i\times \prod_{i\in\Xi^c} \cF_i\times \prod_{i=r+1}^m \cF_i \times \cF_0
\]
where $\Xi^c$ denotes the complement of $\Xi$ in $\{1, \ldots, r\}$. For $\un= (n_1,\ldots, n_m)\break \in \bN_0^m$ with $n_i=0$ for all $i\in \Xi$ we let $\lambda(\Xi, \un) \colon \!\! = (\prod_{i=1}^m\, \ell_i^{n_i})\cdot 1_{\cF_{\Xi}}\in \sC_c^{\flat}(S_1, S_2, \bC_p)$ i.e.\ $\lambda(\Xi, \un)$ is given by
\[
\lambda(\Xi, \un)(x_1,\ldots, x_m, x^{p,\infty}) = \left \{\begin{array}{ll} 
\prod_{i=1}^m\,\ell_i(x_i)^{n_i} & \mbox{if $(x_1,\ldots, x_m, x^{p,\infty})\in \cF_{\Xi}$;}\\
0  & \mbox{otherwise.}
\end{array}\right.
\]
Put $|\un| \,\colon \!\! = \sum_{i=1}^m \, n_i$ and $\un! \,\colon \!\! = \prod_{i=1}^m n_i!$. Then,
\begin{equation}
\label{lambda}
\ell_p^k 1_{\cF} = \sum_{|\un| = k}\, \frac{k!}{\un!}\, \lambda(\emptyset, \un).
\end{equation}
Thus \eqref{ellvanish2} follows from 

\begin{lemma}
\label{lemma:freermod4} 
If $\sharp(\Xi) + |\un| < r$ then $\lambda(\Xi, \un)\in I(\cT_p)\sC_c^{\flat}(S_1, S_2, \bC_p)$.
\end{lemma}

{\em Proof.} We first remark that for two functions $f,g: \bI^{\infty}/U^{p,\infty}\to \bC_p$ and $t\in F^*$ we have
\begin{equation*}
\label{prodmodaug}
(1-t)(f\cdot g) = ((1-t)f)\cdot g + f\cdot ((1-t)g) - ((1-t)f)\cdot ((1-t)g).
\end{equation*}
For $\un, \un'\in \bN_0^m$ write $\un'< \un$ if $n'_i \le n_i$ for all $i$ and $\un'\ne \un$. By using $((1-t) \ell_i)(x) = \ell_i(x) - \ell_i(t^{-1} x) = \ell_i(t)$ one can easily show that we have
\begin{equation}
\label{lambdared}
(1-t) \prod_{i=1}^m\, \ell_i^{n_i} = \sum_{\un'< \un}\, a_{\un'} \,\,\prod_{i=1}^m\, \ell_i^{n'_i}
\end{equation}
for some $a_{\un'} \in \bC_p$. In fact if $\un'\in \bN_0^m$ with $\un'< \un$ and $|\un'| = |\un| -1$ and if $i\in \{1, \ldots, m\}$ with $n'_i = n_i - 1$ then 
\begin{equation}
\label{lambdared2}
a_{\un'} = n_i \ell_i(t).
\end{equation}

We prove the assertion by induction on $|\un|$. Assume first that $\un = \uzero = (0, \ldots, 0)$. Let $i\in \Xi^c$ ($\Xi^c\ne \emptyset$ since $\sharp(\Xi)<r$). Then 
\[
\lambda(\Xi, \uzero) = (1-t_i) \lambda(\Xi\cup\{i\}, \uzero) \in I(\cT_p)\sC_c^{\flat}(S_1, S_2, \bC_p).
\]
Now assume that $|\un|>0$. Since $\sharp(\Xi) + \sum_{i\in \Xi^c} \, n_i \le \sharp(\Xi) + |\un| < r$, there exists $j\in \Xi^c$ with $n_j =0$. Put $\Xi' \colon \!\! = \Xi\cup \{j\}$. Modulo $I(\cT_p)\sC_c^{\flat}(S_1, S_2, \bC_p)$ we obtain
\begin{eqnarray}
\label{lambdared3}
& &\lambda(\Xi,\un) \,\,=\,\, \prod_{i=1}^m\, \ell_i^{n_i}\cdot 1_{\cF_{\Xi}} = \prod_{i=1}^m\, \ell_i^{n_i}\cdot (1-t_j)1_{\cF_{\Xi'}}\\
&&= \,\, (1-t_j) \lambda(\Xi',\un) - ((1-t_j)\prod_{i=1}^m\, \ell_i^{n_i})\cdot 1_{\cF_{\Xi'}} + ((1-t_j)\prod_{i=1}^m\, \ell_i^{n_i})\cdot 1_{\cF_{\Xi}}\nonumber\\
&& \equiv \,\, - ((1-t_j)\prod_{i=1}^m\, \ell_i^{n_i})\cdot 1_{\cF_{\Xi'}} + ((1-t_j)\prod_{i=1}^m\, \ell_i^{n_i})\cdot 1_{\cF_{\Xi}}.\nonumber
\end{eqnarray}
By \eqref{lambdared} and the induction hypothesis we have
\[ ((1-t_j)\prod_{i=1}^m\, \ell_i^{n_i})\cdot 1_{\cF_{\Xi'}} \in \sum_{\un'< \un}\, \bC_p \lambda(\Xi', \un') \subseteq I(\cT_p)\sC_c^{\flat}(S_1, S_2, \bC_p)
\]
and 
\[
((1-t_j)\prod_{i=1}^m\, \ell_i^{n_i})\cdot 1_{\cF_{\Xi}} \in \sum_{\un'< \un}\, \bC_p \lambda(\Xi, \un') \subseteq I(\cT_p)\sC_c^{\flat}(S_1, S_2, \bC_p)
\]
hence $\lambda(\Xi,\un)\in I(\cT_p)\sC_c^{\flat}(S_1, S_2, \bC_p)$.
\enddemo

Next we want to write $(\sum_{i=1}^m \ell_i)^r 1_{\cF}$ -- modulo $I(\cT)\sC_c^{\flat}(S_1, S_2, \bC_p)$ -- as a linear combination of a particular subset of $\{ \lambda(\Xi, \un)| \,\, \sharp(\Xi) + |\un| = r\}$ (this will be used in the proof of Prop.\ \ref{prop:thetalog} (b) below). For that we need to introduce more notation. 

For $\Xi\subseteq \{1, \ldots, r\}$ and a map $f: \Xi\to \{1, \ldots, m\}$ we let $\un(f) \colon \!\! = (\sharp(f^{-1}(1)), \ldots, \sharp(f^{-1}(m)))\in \bN_0^m$. If in particular $f: \Xi\to \{1, \ldots, m\}$ is the inclusion we write $\un(\Xi)$ rather than $\un(f)$. We define
\begin{equation*}
\label{lambdabase}
\Lambda_{\Xi} \colon \!\!= \lambda(\Xi^c, \un(\Xi)) = (\prod_{i\in \Xi}\, \ell_i)\cdot 1_{\cF_{\Xi^c}} .
\end{equation*}
Note that $\Lambda_{\Xi} = \prod_{i\in \Xi} \Lambda_i$ where for $i \in\{1, \ldots, r\}$ we have
\begin{equation*}
\label{lambdabase2}
\Lambda_i \colon \!\!=\Lambda_{\{i\}} = \ell_i \cdot 1_{\cF_i \times (\prod_{j=1, j\ne i} ^r\cO_j)\times (\prod_{j=r+1}^m \cF_j) \times \cF_0}.
\end{equation*}

\begin{lemma}
\label{lemma:freermod5}
Modulo $I(\cT)\sC_c^{\flat}(S_1, S_2, \bC_p)$ we have
\begin{equation*}
\label{lambdabase4}
(\sum_{i=1}^m\, \ell_i)^r 1_{\cF} \equiv   r! \det\left(\begin{matrix} \Lambda_1 + \ell_1(t_1)\Lambda_{\emptyset} & \ell_1(t_2)\Lambda_{\emptyset} & \ldots &  \ell_1(t_r)\Lambda_{\emptyset}\\
\ell_2(t_1)\Lambda_{\emptyset} &\Lambda_2 + \ell_2(t_2)\Lambda_{\emptyset} & \ldots &  \ell_2(t_r)\Lambda_{\emptyset}\\
\vdots & \vdots & \ddots & \vdots\\
\ell_r(t_1)\Lambda_{\emptyset} & \ell_r(t_2)\Lambda_{\emptyset} & \ldots & \Lambda_r + \ell_r(t_r)\Lambda_{\emptyset}
\end{matrix}\right).
\end{equation*}
\end{lemma}

{\em Proof.} For $\Xi \subseteq \{1, \ldots, r\}$ we denote by $M(\Xi)\subseteq \Maps(\Xi, \{1, \ldots, m\})$ the set of maps $f: \Xi\to \{1, \ldots, m\}$ with $f(S) \not\subseteq S$ for all $S\subseteq \Xi$, $S\ne \emptyset$. 

Let $\un\in \bN_0^m$ with $\sharp(\Xi) + |\un| =r$. Firstly, we show that modulo $I(\cT)\sC_c^{\flat}(S_1,\break S_2, \bC_p)$ we have
\begin{equation}
\label{lambdabase3}
\lambda(\Xi, \un) \quad \equiv\quad \un!\,\, \sum_{(\Upsilon, f)} \,\, (-1)^{|\un(f)|}\left( \prod_{i\in \Upsilon}\, \ell_{f(i)}(t_i)\right)\, \Lambda_{\Xi^c-\Upsilon}
\end{equation}
where the sum runs through all pairs $(\Upsilon,f)$ with $\Upsilon\subseteq\Xi^c$ and $f\in M(\Upsilon)$ such that $\un(f) + \un(\Xi^c-\Upsilon) = \un$. 

In the case $n_i>0$ for all $i\in \Xi^c$ we show that both sides of \eqref{lambdabase3} are $=\Lambda_{\Xi^c}$. Since $|\un| \ge \sharp(\Xi^c) = r- \sharp(\Xi) = |\un|$ we have $\un = \un(\Xi^c)$ hence $\lambda(\Xi, \un) = \Lambda_{\Xi^c}$. On the other hand if $\Upsilon\subseteq \Xi^c$ and $f: \Upsilon\to \{1, \ldots, m\}$ is a map with $\un(f) + \un(\Xi^c- \Upsilon) = \un$, then $\un(f) = \un(\Xi^c) - \un(\Xi^c-\Upsilon) = \un(\Upsilon)$ and therefore $f(\Upsilon) = \Upsilon$, hence $\Upsilon =\emptyset$. Consequently, the left hand side of \eqref{lambdabase3} consists of only one summand $\Lambda_{\Xi^c}$. 

Now assume $n_j=0$ for some $j\in \Xi^c$. By \ref{lemma:freermod4}, \eqref{lambdared}, \eqref{lambdared2} and \eqref{lambdared3} we get
\begin{eqnarray*}
\lambda(\Xi,\un) & \equiv &  - ((1-t_j)\prod_{i=1}^m\, \ell_i^{n_i})\cdot 1_{\cF_{\Xi\cup\{j\}}}\\
 & \equiv  & - \sum_{\un'< \un, |\un'| = |\un|-1}\, a_{\un'} \,\, \lambda(\Xi\cup\{j\},\un')\\
 & = & - \sum_{(f, \un')}\, n_{f(j)}\, \ell_{f(j)}(t_j) \,\,\lambda(\Xi\cup\{j\},\un')
\end{eqnarray*}
where the last sum runs through all pairs $(f,\un')$ with $f\in M(\{j\})$ and $\un(f) + \un' = \un$. Now \eqref{lambdabase3} can be easily deduced by induction on $\sharp(\Xi^c)$.

By \eqref{lambda} and \eqref{lambdabase3} (for $\Xi = \emptyset$) we have
\begin{equation}
\label{lambdabase4a}
(\sum_{i=1}^m\, \ell_i)^r 1_{\cF}  \,\, \equiv  \,\, r! \, \sum_{\Xi\subseteq \{1, \ldots, r\} } (-1)^{|\un(\Xi^c)|}\left( \sum_{f\in M(\Xi^c)}\, \prod_{i\in \Xi^c}\, \ell_{f(i)}(t_i)\right) \Lambda_{\Xi}.
\end{equation}
On the other hand
\[
\det\left(\begin{matrix} \Lambda_1 + \ell_1(t_1)\Lambda_{\emptyset} & \ell_1(t_2)\Lambda_{\emptyset} & \ldots &  \ell_1(t_r)\Lambda_{\emptyset}\\
\ell_2(t_1)\Lambda_{\emptyset} &\Lambda_2 + \ell_2(t_2)\Lambda_{\emptyset} & \ldots &  \ell_2(t_r)\Lambda_{\emptyset}\\
\vdots & \vdots & \ddots & \vdots\\
\ell_r(t_1)\Lambda_{\emptyset} & \ell_r(t_2)\Lambda_{\emptyset} & \ldots & \Lambda_r + \ell_r(t_r)\Lambda_{\emptyset}
\end{matrix}\right) = \sum_{\Xi\subseteq \{1, \ldots, r\} } \,\,a_{\Xi} \, \Lambda_{\Xi}
\]
where \(
a_{\Xi} = \det\left( \ell_i(t_j) \right)_{i,j\in \Xi^c} =  \det\left(\begin{matrix} \ell_{i_1}(t_{i_1}) & \ldots &  \ell_{i_1}(t_{i_k})\\
\vdots & &  \vdots\\
\ell_{i_k}(t_{i_1}) & \ldots &  \ell_{i_k}(t_{i_k})
\end{matrix}\right)
\) if $\Xi^c = \{i_1 < \ldots < i_k\}$.

Note that $\sum_{j=1}^m \, \ell_j(t_i) = \log_p(\Norm_{K/\bQ}(t_i)) = 0$ since $\Norm_{K/\bQ}(t_i)$ is a power of $p$ (for all $i\in \{1, \ldots, m\}$). The assertion follows from \eqref{lambdabase4a} and the following result about determinants.\enddemo

\begin{lemma}
\label{lemma:determinant}
Let $k\le m$ be positive integers and let $\left( a_{ij} \right)_{i=1, \ldots, k, j=1, \ldots, m}$ be a $k\times m$--matrix with entries in a commutative ring such that $\sum_{j=1}^m a_{ij} =0$ for all $i=1, \ldots, k$. Then,
\[
\det\left( a_{ij} \right)_{i,j=1, \ldots, k} = (-1)^k\, \sum_{f}\, \prod_{i=1}^k\, 
a_{i f(i)}
\]
where the sum runs through all maps $f: \{1, \ldots,k\} \to \{1, \ldots,m\}$ with $f(S) \not\subseteq S$ for all $S\subseteq \{1, \ldots,k\}$, $S\ne \emptyset$.
\end{lemma}

{\em Proof.}\footnote{Due to V.\ Paskunas} By replacing $a_{ii}$ with $-\sum_{j\ne i} a_{ij}$ in  $\sum_{\sigma\in S_k} \sign(\sigma) \prod_{i=1}^m a_{i \sigma(i)}$ and expanding the sum we get
\begin{equation*}
\det\left( a_{ij} \right)_{i,j=1, \ldots, k} \,\,\, = \,\,\,\sum_{(\Xi, \sigma, g)}\, (-1)^{k-\sharp(\Xi)} \sign(\sigma) \prod_{i\in \Xi} a_{i\sigma(i)} \prod_{i\in \Xi^c} a_{i g(i)}
\end{equation*}
where the sum ranges over all triples $(\Xi, \sigma, g)$ with $\Xi\subseteq \{1, \ldots, k\}$, $\sigma$ is a permutation of $\Xi$ without fixed points and $g$ is a map $\Xi^c \colon \!\! = \{1, \ldots, k\}-\Xi \to \{1, \ldots, m\}$ without fixed points.

Let $f: \{1, \ldots,k\} \to \{1, \ldots,m\}$ be a map without fixed points and let $\Xi\subseteq \{1, \ldots, k\}$ be the largest subset with $f(\Xi) = \Xi$. If we decompose the permutation $f|_{\Xi}$ into disjoint cycles $\sigma_1\cdots\sigma_t$ of length $l_1, \ldots, l_t$, then it is easy to see that the coefficient of $\prod_{i=1}^k\, 
a_{i f(i)}$ is $(-1)^k \prod_{j=1}^t (1 + \sign(\sigma_j)(-1)^{l_j})$. Thus it is $=0$ except when $\Xi=\emptyset$.\enddemo

{\em Proof Prop.\ \ref{prop:thetalog} (b)} We first show
\begin{equation}
\label{keyequal}
[\ell^r 1_{\cF}] \, = \, (-1)^{\binom{r}{2}} (c_{\ell_1} \cup \ldots \cup c_{\ell_r}) \cap \varrho_{\cT}
\end{equation}
in $H_0(\cT,\sC_c^{\flat}(S_1, S_2, \bC_p))$. As before if $\cT'$ is a subgroup of finite index of $\cT$ the injectivity of $\res: H_0(\cT, \sC_c^{\flat}(S_1, S_2, \bC_p))\to H_0(\cT', \sC_c^{\flat}(S_1, S_2, \bC_p))$ together with $\res(\varrho_{\cT}) = \varrho_{\cT'}$ implies that in order to prove \eqref{keyequal} we may shrink $\cT$ so we can assume that $\cT$ is of the form \eqref{tdecomposition}. 

By \eqref{ellvanish3} we have modulo $I(\cT_p)\sC_c^{\flat}(S_1, S_2, \bC_p)$
\[
\ell^r 1_{\cF} \,=\, \sum_{j=0}^r \binom{r}{j}\ell_p^j(\ell^P)^{r-j}1_{\cF} \,\equiv\, \ell_p^r 1_{\cF} \,=\, (\sum_{i=1}^m \ell_i)^r 1_{\cF}
\]
so we may replace $\ell^r 1_{\cF}$ by $(\sum_{i=1}^m \ell_i)^r 1_{\cF}$ on the left hand side of \eqref{keyequal}.

We will use the notation of the proof of part (a). Furthermore we put $\cT_1 \colon \!\! = \langle t_1,\ldots, t_r\rangle, \cT_2 \colon \!\! =\langle t_{r+1},\ldots, t_m\rangle \times \cT_0$ and $\cF_2 \colon \!\! = \prod_{i=r+1}^m \cF_i \times \cF_0$ so that $\cT= \cT_1\times \cT_2$, $\cF = \prod_{i=1}^r \cF_i\times \cF_2$ and $\cF_2$ is a $\cT_1$-stable fundamental domain for the action of $\cT_2$ on $\bI^{S_1,\infty}/ U^{p,\infty}$. 

Since $\varrho_1\in H_r(\cT_1, \bZ)\cong H_r(\bR^r/\cT_1, \bZ)$ can be represented by the $r$-cycle $\sum_{\sigma\in S_r} \,\sign(\sigma)\, [t_{\sigma(1)}| \ldots| t_{\sigma(r)}]$ (compare the proof of Prop. \ref{prop:thetaord}) together with Remark \ref{remarks:varrho} (b) we conclude that 
\begin{equation*}
\label{keyequal2}
\sum_{\sigma\in S_r} \,\sign(\sigma)\, z_{\ell_1}(t_{\sigma(1)}) \otimes t_{\sigma(1)} z_{\ell_2}(t_{\sigma(2)})\otimes \ldots \otimes t_{\sigma(1)}\ldots t_{\sigma(r-1)}z_{\ell_r}(t_{\sigma(r)})\otimes 1_{\cF_2}
\end{equation*}
is a representative of $(-1)^{\binom{r}{2}}(c_{\ell_1} \cup \ldots \cup c_{\ell_r}) \cap \varrho\in H_0(\cT,\sC_c^{\flat}(S_1, \bC_p))$. We have
\[
z_{\ell_i}(t_j) \equiv \delta_{ij} \ell_i 1_{\cF_i} + \ell_i(t_j)1_{\cO_i} \quad \mbox{and} \quad t_k z_{\ell_i}(t_j) \equiv  z_{\ell_i}(t_j)
\]
(modulo $K \cdot 1_{\cF_i}$) for all $i,j,k \in \{1, \ldots, r\}$ with $j\ne k$. Hence by \ref{lemma:freermod4} we obtain 
\begin{eqnarray*}
& z_{\ell_1}(t_{\sigma(1)}) \otimes t_{\sigma(1)} z_{\ell_2}(t_{\sigma(2)})\otimes \ldots \otimes t_{\sigma(1)}\ldots t_{\sigma(r-1)}z_{\ell_r}(t_{\sigma(r)})\otimes 1_{\cF_2}\\
& \equiv \,\,
(\delta_{1\sigma(1)} \ell_1 1_{\cF_1} + \ell_1(t_{\sigma(1)})1_{\cO_1}) \otimes \ldots \otimes(\delta_{r\sigma(r)} \ell_r 1_{\cF_r} + \ell_r(t_{\sigma(r)})1_{\cO_r})\otimes 1_{\cF_2}
\end{eqnarray*}
(modulo $I(\cT)$) for all $\sigma\in S_r$. Now \eqref{keyequal} follows from \ref{lemma:freermod5}.

To finish the proof of Prop.\ \ref{prop:thetalog} (b) consider the commutative diagram 
\begin{equation*}
\label{cap3}
{\footnotesize
\xymatrix@+1pc{H^r(F^*_+, C^{\flat}_c(F_{S_1}))\!\times \! H_r(\cT,C_c(\bI^{S_1,\infty}/ U^{p,\infty}))^{E_+} \ar[r]^{\hspace{2.3cm}\cap\circ (\res\times \id)}\ar[d]^{\cap\circ (\id\times \eqref{coreta}))} & H_0(\cT,\sC_c^{\flat}(S_1, S_2))^{E_+}\ar[d]^{\cap\eta}\\
H_{d}(F^*_+,\sC_c^{\flat}(S_1,S_2)) \ar[r]^{\coinf} & H_d(E_+,\! H_0(\cT\! ,\sC_c^{\flat}(S_1, S_2)))}
}
\end{equation*}
where $C^{\flat}_c(F_{S_1})= C^{\flat}_c(F_{S_1},\bC_p)$, $C_c(\bI^{S_1,\infty}/ U^{p,\infty}))=C_c(\bI^{S_1,\infty}/ U^{p,\infty},\bC_p)$ etc.\ and the maps $\res$ and $\coinf$ denote the restriction and coinflation with respect to $\cT \leq F^*_+$ (recall that the latter is an isomorphism by Prop.\ \ref{prop:freermod2}). By \eqref{keyequal} the image of the pair $(c_{\ell_1} \cup \ldots \cup c_{\ell_r}, \varrho_{\cT})$ under the composition of the upper horizontal map, the right vertical map and the inverse of $\coinf$ is $(-1)^{\binom{r}{2}}\partial((\log_p \circ \cN)^r)$. On the other hand its image under the first vertical map is $(c_{\ell_1} \cup \ldots \cup c_{\ell_r}) \cap \vartheta$.
\enddemo

\begin{remark}
\label{remark:sigmafix}
\rm For $\mu\in \{\pm 1\}$ let $e(\mu) \in \bZ/2\bZ$ be given by $\mu = (-1)^{e(\mu)}$. Let $\Sigma \colon \!\! = \{\pm 1\}^{d+1}$. We write elements of $\Sigma$ in the form $\umu = (\mu_0, \ldots, \mu_d)$. Define the pairing $\langle \,\,\, , \,\,\, \rangle: \Sigma\times \Sigma \lra \{\pm 1\}$ by $\langle \umu, \unu\rangle = (-1)^{\sum_i\, e(\mu_i) e(\nu_i)}$. Let $\bfC$ be a field of characteristic zero. For a $\bfC[\Sigma]$-module $V$ and $\umu\in \Sigma$ we put $V_{\umu} = \{v\in V\mid \unu v = \langle\unu, \umu\rangle \, v \,\,\,\forall \, \unu \in \Sigma\}$ so that $V = \bigoplus_{\umu\in \Sigma} \, V_{\umu}$. For $v\in V$ we denote by $v_{\umu}\in V_{\umu}$ its projection to $V_{\umu}$. If $\umu = (+1, \ldots, +1)$ we shall also write $v_+$ instead of $v_{\umu}$. 

We identify $F^*/F^*_+$ with $\Sigma$ via the isomorphism $F^*/F^*_+ = F_{\infty}^*/U_{\infty} \cong \prod_{i=0}^d \,\bR^*/\bR^*_+ \cong \Sigma$. Hence for any $F^*$-module $M$ we obtain an action of $\Sigma$ on $H_q(F^*_+, M)$. Note that $\vartheta$ is $\Sigma$-invariant (since $\Sigma$ acts trivially on $H_{d+r}(F^*_+,C_c(\bI^{S_1,\infty}/ U^{S_1,\infty}, \bZ))$) as well as the cohomology classes defined in Def.\ \ref{definition:logclass}. Consequently, the cap-product  $(c_{\ell_{\fp_1}} \cup \ldots \cup c_{\ell_{\fp_r}}) \cap \vartheta$ for $\ell_{\fp} = \ord_{\fp}$ or $\ell_{\fp} = \log_p\circ N_{F_{\fp}/\bQ_p}$ are $\Sigma$-invariant.
\end{remark}

\subsection{$p$-adic $L$-functions attached to cohomology classes}
\label{subsection:abstrezc}

Let $S_1, S_2$ be arbitrary (possibly empty) disjoint subsets of $S_p$. For a ring $R$ and an $R$-module $M$ put
\[
\sD(S_1, S_2, M)\! = \!\Dist(F_{S_1}\times F_{S_2}^*\times \bI^{p,\infty}/U^{p,\infty}\!, M) \! = \! \Hom_R(\sC_c^0(S_1, S_2, R), M).
\]
It is easy to see that the functor $M \mapsto \sD(S_1, S_2, M)$ is exact. Let
\begin{equation}
\label{distr2}
\langle\,\,\, ,\,\,\, \rangle: \, \sD(S_1, S_2, M)\times \sC_c^0(S_1, S_2,R) \lra M
\end{equation}
be the canonical (evaluation) pairing. 
Also for $\fp\in S_1$ by \eqref{fubini} we obtain a pairing 
\begin{equation*}
\label{equation:distr2a}
\sD(S_1, S_2, M)\times C^0_c(F_{\fp},R) \lra \sD(S_1-\{\fp\}, S_2, M).
\end{equation*}
If $H$ is a subgroup of $\bI^{\infty}$ and $M$ is an $R[H]$-module we define a $H$-action on $\sD(S_1, S_2, M)$ by requiring that $\langle x\lambda, x f\rangle = x\langle \phi,  f\rangle$ for all $x\in H$, $f\in \sC^0_c(S_1, S_2, R)$ and $\lambda\in \sD(S_1, S_2, M)$. 

If $K$ is a $p$-adic field and $V$ a $K$-Banach space then we denote the subspace of measures of $\sD(S_1, S_2, V)$ by 
\[
\sD^b(S_1, S_2, V) \,=\, \Meas(F_{S_1}\times F_{S_2}^*\times \bI^{p,\infty}/U^{p,\infty}, V).
\]
The pairing \eqref{distr2} when restricted to $\sD^b(S_1, S_2, V)$ extends canonically to a pairing
\begin{equation}
\label{meas2}
\langle\,\,\, ,\,\,\, \rangle: \, \sD^b(S_1, S_2, V)\times \sCd(S_1, S_2, K) \lra V.
\end{equation}
Also for $\fp\in S_1$ by \eqref{fubini2} we obtain a pairing
\begin{equation}
\label{meas2a}
\sD^b(S_1, S_2, V)\times \Cd(F_{\fp}, K) \lra \sD^b(S_1-\{\fp\}, S_2, V).
\end{equation}

Assume now that $S_1 \cup S_2 = S_p$, $H= F_+^*$ and that $F_+^*$ acts trivially on $M$. Again, we order the places above $p$ so that $S_1 = \{\fp_1, \ldots, \fp_r\}$ and $S_2=\{\fp_{r+1}, \ldots, \fp_{m}\}$. The pairing \eqref{distr2} yields the bilinear map
\begin{equation}
\label{cap1}
\cap: H^{d-1}(F^*_+,\sD(S_1, S_2, M))\times H_{d-1}(F^*_+, \sC_c^0(S_1, S_2, R)) \to H_0(F^*_+, M) = M.
\end{equation}
For $\kappa\in H^d(F^*_+, \sD(S_1, S_2, M))$ define $\mu_{\kappa}\in \Dist(\cG_p, M)$ by
\begin{equation}
\label{cap2a}
\int_{\cG_p} f(\gamma) \,\mu_{\kappa}(d\gamma) \,\,= \,\,\kappa \cap \partial(f) 
\end{equation}
for all $f\in C^0(\cG_p, R)$. 

Suppose now that $R=K$ is a $p$-adic field and $M=V$ a finite dimensional $K$-vector space and let $\kappa\in H^d(F^*_+, \sD^b(S_1, S_2, V))$. By abuse of notation we denote its image under $H^d(F^*_+, \sD^b(S_1, S_2, V)) \to H^d(F^*_+, \sD(S_1, S_2, V))$ also by $\kappa$. It is then easy to see that $\mu_{\kappa}$ is actually a measure. Thus we obtain a map 
\begin{equation}
\label{cohmeas}
H^d(F^*_+, \sD^b(S_1, S_2, V))\lra \Meas(\cG_p, V), \,\, \kappa \mapsto \mu_{\kappa}.
\end{equation}
The integral $C(\cG_p, K) \to V, f\mapsto \int_{\cG_p} f \mu_{\kappa}$ admits also a description as a cap-product. More precisely \eqref{cap2a} holds more generally for all $f\in C(\cG_p, K)$ (where $\partial$ denotes now the map \eqref{recipr2} for $R=K$ equipped with the topology induced by $|\,\,\, |_p$ and the cap-product is induced by \eqref{meas2}).

Recall that $\cN: \cG_p \to \bZ_p^*$ is defined by $\gamma \zeta = \zeta^{\cN(\gamma)}$ for all $p$-power roots of unity $\zeta$. For $s\in \bZ_p$ and $\gamma\in \cG_p$ we put $\langle \gamma\rangle^s \colon \!\! = \exp_p(s \log_p(\cN(\gamma)))$. 

\begin{definition}
\label{definition:lpfunc}
Let $K$ be a $p$-adic field and $V$ a finite-dimensional $K$-vector space. We define the $p$-adic $L$-function of $\kappa\in H^d(F^*_+, \sD^b(S_1, S_2, V))$ by
\[
L_p(s, \kappa) \colon \!\! = \int_{\cG_p}\, \langle \gamma\rangle^s \mu_{\kappa}(d\gamma).
\]
\end{definition}

The main result of this chapter is the following

\begin{theorem}
\label{theorem:abstrezc} 
Let $r \colon \!\! = \sharp(S_1)$ and let $\kappa\in H^d(F^*_+, \sD^b(S_1, S_2, V))$.

\noi (a) $L_p(s, \kappa)$ is a locally analytic $V$-valued function on $\bZ_p$. We have 
\[
\ord_{s=0} L_p(s, \kappa) \ge r
\]
(b) For $\fp\in S_1$ put $\ell_{\fp} \colon \!\! = \log_p \circ \Norm_{F_{\fp}/\bQ_p}: F_{\fp}^* \to K$. For the $r$-th derivative of $L_p(s, \kappa)$ at $s=0$ we have 
\[
L^{(r)}_p(0, \kappa)\,\, =\,\, (-1)^{\binom{r}{2}}\,\,r! \,\, (\kappa \cup c_{\ell_{\fp_1}} \cup \ldots \cup c_{\ell_{\fp_{r}}}) \cap \vartheta.
\]
Here the cup-product is induced by \eqref{meas2a} and the cap-product by \eqref{meas2}.
\end{theorem}

{\em Proof.} We have
\[
L^{(k)}_p(0, \kappa) = \int_{\cG_p}\, (\log_p \circ \cN)^k \mu_{\kappa}(d\gamma) = \kappa\cap \partial((\log_p \circ \cN)^k)
\]
for all $k\in \bN$. Hence the assertion follows from Prop.\ \ref{prop:thetalog}.\enddemo

\begin{remark}
\label{remark:sigmaaction}
\rm The group $\Sigma\cong F^*/F^*_+$ acts on $H^q(F^*_+, \sD(S_1, S_2,  M))$  or $H^q(F^*_+, \sD^b(S_1, S_2, V))$. Let $\cG_p^+$ be the Galois group of the maximal abelian extension of $F$ which is unramified outside $p$. By class field theory we have an exact sequence $F^*/F^*_+ = F_{\infty}^*/U_{\infty} \to \cG_p \to \cG^+_p\to 1$, which yields an action of $\Sigma$ on $\cG_p$. It is easy to see that \eqref{cohmeas} is $\Sigma$-equivariant. The fact that $\gamma \mapsto \langle \gamma\rangle^s$ factors through $\cG_p \to \cG^+_p$ implies that $L_p(s, \kappa) = L_p(s, \kappa_+)$ for all $\kappa\in H^d(F^*_+, \sD^b(S_1, S_2, V))$. Also by Remark \ref{remark:sigmafix} we have 
$\kappa \cap (c_{\ell_{\fp_1}} \cup \ldots \cup c_{\ell_{\fp_r}}) \cap \vartheta = \kappa_+ \cap (c_{\ell_{\fp_1}} \cup \ldots \cup c_{\ell_{\fp_r}}) \cap \vartheta$.
\end{remark}

\subsection{Integral cohomology classes}
\label{subsection:intcohom} 

For a given cohomology class $\kappa\in H^d(F^*_+, \sD(S_1, S_2, \bC))$
we will introduce a condition -- called {\it integral} -- which
guarantees that $\mu_{\kappa}$ is a $p$-adic measure (in the sense of section \ref{subsection:padicmeas}) and which allows
us to apply Theorem \ref{theorem:abstrezc}. To begin with we define the module of {\it periods} of $\kappa$.

\begin{definition}
\label{definition:period}
Let $\kappa\in H^d(F^*_+, \sD(S_1,  S_2, \bC))$ and let $R$ be a subring of $\bC$. The image of 
\begin{equation*}
\label{equation:cap5}
H_d(F^*_+, \sC_c^0(S_1, S_2, R)) \lra H_0(F^*_+, \bC) = \bC,\,\,\, x \longmapsto\kappa\cap x
\end{equation*}
will be denoted by $L_{\kappa,R}$. If $R\subset \barQ$ then it is called the $R$-module of periods of $\kappa$.
\end{definition}

\begin{lemma}
\label{lemma:period}
Let $R \subseteq \barQ$ be a Dedekind ring. 

\noi (a) If $R'\supseteq R$ is a subring of $\bC$ then $L_{\kappa,R'} = R'L_{\kappa,R}$.

\noi (b) If $\kappa\ne 0$ then $L_{\kappa,R}\ne 0$.
\end{lemma}

{\em Proof.} (a) Since $\sC_c^0(S_1, S_2, R')= \sC_c^0(S_1, S_2, R)\otimes R'$ and $R'$ is flat $R$-algebra we have $H_d(F^*_+, \sC_c^0(S_1, S_2, R))\otimes R' = H_d(F^*_+, \sC_c^0(S_1, S_2, R'))$.

\noi (b) By (a) it is enough to show $L_{\kappa,\bC}\ne 0$. This follows from the fact that the pairing \eqref{cap1} (for $M=R = \bC$) is nondegenerated.\enddemo

\begin{definition}
\label{definition:integral}
A cohomology class $\kappa\in H^d(F^*_+, \sD(S_1, S_2,  \bC))$,
$\kappa\ne 0$ is said to be integral (or better $p$-integral) if there
exists a Dedekind ring $R \subset \barcO$ such that $\kappa$ lies in
the image of \(H^d(F^*_+, \sD(S_1, S_2, R))\otimes_R \bC \to H^d(F^*_+,
\sD(S_1, S_2,  \bC))\). If in addition there exists a finitely generated
$R$-submodule $M \subseteq H^d(F^*_+, \sD(S_1, S_2, R))$ of rank $\le 1$ (i.e.\ $\rank_R M/M_{\tor} \le1$) such that $\kappa$ lies in the image of $M\otimes_R \bC \to H^d(F^*_+, \sD(S_1, S_2,  \bC))$ then $\kappa$ is called integral of rank $\le 1$.
\end{definition}

\begin{prop}
\label{prop:integral} 
Let $\kappa\in H^d(F^*_+, \sD(S_1, S_2,  \bC))$, $\kappa\ne 0$. The
following conditions are equivalent.

\noi (i) $\kappa$ is integral (resp.\ integral of rank $\le 1$).

\noi (ii) There exists a Dedekind ring $R\subseteq \barcO$ such that
$L_{\kappa, R}$ is a finitely generated $R$-module (resp.\ $L_{\kappa,
  R}$ is $=0$ or an invertible $R$-module).

\noi (iii) There exists a Dedekind ring $R\subseteq \barcO$, a finitely generated $R$-module $M$ (resp.\ an invertible $R$-module $M$ of rank 1) and an $R$-linear map $f: M\to\bC$ such that $\kappa$ lies in the image of the induced map \(f_*: H^d(F^*_+, \sD(S_1, S_2, M)) \break \to H^d(F^*_+, \sD(S_1, S_2,  \bC))\).
\end{prop}

{\em Proof.} We consider only the case of arbitrary rank and will leave the necessary modifications to the rank $\le 1$ case to the reader. 

\noi (i) $\Rightarrow$ (ii) Let $R$ be as in Definition \ref{definition:integral}. If we write $\kappa$ in the form $\kappa = \sum_{i=1}^n\,\Omega_i\kappa_i$ with $\kappa_i\in \Image(H^d(F^*_+, \sD(S_1, S_2, R)) \to H^d(F^*_+, \sD(S_1, S_2,  \bC))$ and $\Omega_i\in \bC$ then $L_{\kappa,R}\subseteq R\Omega_1 + \ldots + R\Omega_n$. 

\noi (ii) $\Rightarrow$ (iii) Consider the diagram
\begin{equation*}
\label{equation:homcohom}
\begin{CD}
H^d(F^*_+, \sD(S_1, S_2, L_{\kappa,R})) @>>> \Hom_R(H_d(F^*_+,\sC_c^0(S_1, S_2,  R)), L_{\kappa,R})\\
@VVV@VVV\\
H^d(F^*_+, \sD(S_1, S_2, \bC)) @>>> \Hom_R(H_d(F^*_+,\sC_c^0(S_1, S_2,  R)), \bC)) 
\end{CD}
\end{equation*}
where the horizontal maps are induced by the cap-product and the vertical maps by the inclusion $L_{\kappa,R} \hookrightarrow \bC$. By the universal coefficient theorem the lower horizontal map is an isomorphism and the kernel and cokernel of the upper horizontal map are $R$-torsion. Hence some multiple $a\cdot \kappa$ with $a\in R$, $a\ne 0$ is contained in the image of the left vertical map. Define $f: L_{\kappa,R} \to \bC, \Omega \mapsto a^{-1}\Omega$.

\noi (iii) $\Rightarrow$ (i) We may assume that $M = R^n$ (for example replace $M$ by $f(M)$ and $f$ by the inclusion and then enlarge $M$ if necessary). Let $\Omega_1, \ldots, \Omega_n\in \bC^n$ be the images of the standard basis under $f$. It follows 
\[
\kappa\in \Image(f_*) = \sum_{i=1}^n\, \Omega_i \cdot \Image(H^d(F^*_+, \sD(S_1, S_2, R)) \to H^d(F^*_+, \sD(S_1, S_2,  \bC))).
\] 
\enddemo

\begin{corollary}
\label{corollary:arithmetmeas} 
Assume $\kappa\in H^d(F^*_+, \sD(S_1, S_2,  \bC))$ is integral and let
$R\subseteq \barcO$ be as in Definition \ref{definition:integral}. Then,

\noi (a) $\mu_{\kappa}$ is a $p$-adic measure. 

\noi (b) The map $H^d(F^*_+, \sD(S_1, S_2, L_{\kappa,R}))\otimes \barQ\to H^d(F^*_+, \sD(S_1, S_2,  \bC))$ is injective and $\kappa$ lies in its image.
\end{corollary}

{\em Proof.} (a) The image of $C^0(\cG_p, \overline{\cO}) \lra \bC, f\mapsto \int\, f\,\, d\mu_{\kappa}$ is contained in $L_{\kappa,\barcO}$. The assertion follows from Proposition \ref{prop:integral}.

\noi (b)  follows from the proof of (ii) $\Rightarrow$ (iii) above.
\enddemo

Let $\kappa$, $R$ be as above. By abuse of notation we define the $p$-adic $L$-function of $\kappa$ by $L_p(s, \kappa) \colon \!\! = \int_{\cG_p}\, \langle \gamma\rangle^s \mu_{\kappa}(d\gamma)$. Because of (b) we can view $\kappa$ as an element of $H^d(F^*_+, \sD(S_1, S_2, L_{\kappa, R}))\otimes_R \barQ$. Put $V_{\kappa} =   L_{\kappa,R}\otimes_R\bC_p$ and let $\wkappa$ denote the image of $\kappa$ under the homomorphism
\begin{equation}
\label{equation:kapmeas}
H^d(F^*_+, \sD(S_1, S_2, L_{\kappa, R}))\otimes_R\bC_p  \lra H^d(F^*_+, \sD^b(S_1, S_2,V_{\kappa}))
\end{equation}
induced by the obvious map $\sD(S_1, S_2, L_{\kappa,R})\to \sD(S_1, S_2, L_{\kappa,R})\otimes_R\bC_p\to\sD^b(S_1, S_2, V_{\kappa})$. It is easy to see that the $\wkappa$ does not depend on the choice of $R$. Since $L_p(s, \kappa) = L_p(s, \wkappa)$ we can apply Theorem \ref{theorem:abstrezc}.

\begin{corollary}
\label{corollary:abstrezc2} 
Assume $\kappa\in H^d(F^*_+, \sD(S_1, S_2,  \bC))$ is integral. For $\fp\in S_1$ put $\ell_{\fp} \colon \!\! = \log_p \circ \Norm_{F_{\fp}/\bQ_p}$. Then,

\noi (a) $\ord_{s=0} L_p(s, \kappa) \ge r$.

\noi (b) \(L^{(r)}_p(0, \kappa)\,\, =\,\, (-1)^{\binom{r}{2}}\,\,r! \,\, (\wkappa_+ \cup c_{\ell_{\fp_1}} \cup \ldots \cup c_{\ell_{\fp_{r}}}) \cap \vartheta.\)
\end{corollary}

\subsection{Another construction of distributions on $\cG_p$}
\label{subsection:globalmeas}

Let $\cA(\bG_m)$ be the space of smooth $\bC$-valued functions on $\bI/F^*$ which are rapidly decreasing as $|x| \to \infty$ or $|x| \to 0$ (i.e.\ for $f\in \cA(\bG_m)$ and $N>0$ there exists a constant $C>0$ such that $|f(x)| < |x|^{-N}$ for $|x|>C$ and $|f(x)| < |x|^{N}$ for $|x|< C^{-1}$).

Let $S_1, S_2$ be disjoint subsets of $S_p$ with $S_1\cup S_2= S_p$. We consider maps $\phi: \sC^0_c(S_1, S_2, \bZ)\times F_{\infty} \to \bC$ with the following properties
\begin{itemize}
\label{itemize:adelmeas}
\item[(i)] For $x_{\infty}\in F_{\infty}$ the map $\phi(\wcdot, x_{\infty}): \sC^0_c(S_1, S_2, \bZ)\to \bC, f\mapsto \phi(f, x_{\infty})$ lies in $\sD(S_1, S_2, \bC)$.
\item[(ii)]
For all $f\in \sC^0_c(S_1, S_2, \bZ)$ the function 
\[
\bI^{\infty}\times F_{\infty}\to \bC, \quad x=(x^{\infty}, x_{\infty})\mapsto \phi(x^{\infty} f, x_{\infty})
\]
lies in $\cA(\bGm)$. In particular we have $\phi(\xi f, \xi x_{\infty}) = \phi(f, x_{\infty})$
 for all $\xi \in F^*$.
\end{itemize}
We denote the space of all $\phi$ satisfying (i),(ii) by $\cD(\bG_m,S_1)$. Note that the map
\[
C_c^0(F_{S_1} \times F_{S_2}^*,\bZ)\times C_c^0(\bI^{p,\infty}/^{p,\infty} ,\bZ)\lra  \sC^0_c(S_1, S_2, \bZ),\,\,\, (f,g)\mapsto f\otimes g
\]
induces an isomorphism $C_c^0(F_{S_1} \times F_{S_2}^*,\bZ)\times C_c^0(\bI^{p,\infty}/^{p,\infty} ,\bZ)\cong  \sC^0_c(S_1, S_2, \bZ)$ and that any element of $C_c^0(\bI^{p,\infty}/^{p,\infty} ,\bZ)$ can be written as a finite sum of the characteristic functions of elements of $\bI^{p,\infty}/^{p,\infty}$. Hence we can (and will) view an element $\phi\in\cD(\bG_m,S_1)$ also as a map
\begin{equation}
\label{phiold}
\phi: C_c^0(F_{S_1} \!\times F_{S_2}^*,\bZ)\times \bI^p/U^{p, \infty}\to \bC,\,\, (f, x^p) \mapsto \phi(f, x^p).
\end{equation}
For $f\in C_c^0(F_{S_1} \times F_{S_2}^*,\bZ)$ we denote by $\phi_f\in \cA(\bGm)$ the map 
\[
\bI = F_p\times \bI^p\to \bC, \quad \phi_f(x_p,x^p) \colon \!\! = \phi(x_p f, x^p)
\]
In particular for a compact open subset $U$ of $F_{S_1} \times F_{S_2}^*$ we define $\phi_U\in \cA(\bGm)$ as $\phi_U = \phi_{1_U}$. 

Given $\phi\in \cD(\bG_m, S_1)$, $f\in C^0(\bI/F^*, \bC)$ and $s\in \bC$ we define the integral $\int_{\bI/F^*} f(x)|x|^s \phi(dx^{\infty}, x_{\infty}) d^{\times}x_{\infty}$ as follows. By Lemma \ref{lemma:locallycomp} below, there exists an open subgroup $U$ of $U_p$ such that $f(x_pu, x^p) = f(x_p, x^p)$ for all $(x_p, x^p)\in F_p^*\times \bI^p$ and $u\in U$. We set
\begin{equation}
\label{equation:intdistr1}
\int_{\bI/F^*} f(x)|x|^s \phi(dx^{\infty}, x_{\infty}) d^{\times}x_{\infty}\, = \,[U_p:U] \int_{\bI/F^*} f(x) |x|^s\phi_U(x) d^{\times}x.
\end{equation}
It is easy to see that the integral does not depend on the choice of $U$. Moreover since $\phi_U$ is rapidly decreasing it is holomorphic in $s$. Hence there exists a unique distribution $\mu= \mu_{\phi}$ on $\cG_p$ such that
\begin{equation*}
\label{equation:intdistr2}
\int_{\cG_p} f(\gamma) \mu_{\phi}(d\gamma) = \int_{\bI/F^*} f(\rho(x)) \phi(dx_p, x^p) d^{\times}x^p
\end{equation*}
for all $f\in C^0(\cG_p, \bC)$ (here $\rho: \bI/F^*\to \cG_p$ denotes the reciprocity map).
\begin{lemma}
\label{lemma:locallycomp} 
Let $X$ be a set and $f: \bI/F^* \to X$ be a locally constant map. Then there exists an open subgroup $U$ of $\bI$, such that $f$ factors through $\bI/F^*U$.
\end{lemma}

{\em Proof.} Since $U_{\infty} = \prod_{v\in S_{\infty}} \bR^*_+$ is connected, $f$ factors as
$\bI/F^*\to \bI/F^*U_{\infty}\break\stackrel{\bar{f}}{\lra} X$ and since $\bI/F^*U_{\infty}$ is profinite, $\bar{f}$ factors through a finite quotient of $\bI_F/F^*U_{\infty}$.\enddemo

We will construct now a cohomology class $\kappa = \kappa_{\phi} \in H^d(F^*_+, \sD(S_1, S_2, \bC))$ such that $\mu_{\phi} = \mu_{\kappa}$. Put $S_{\infty}^0 = S_{\infty}-\{\infty_0\} = \{\infty_1, \ldots, \infty_d\}$ and $U_{\infty}^0 = \prod_{w\in S_{\infty}^0}\bR^*_+$. We write elements of $F_{\infty} = F_{\infty_0}^* \times F_{S_0}$ as pairs $(x_0, x^{0})$. For $\phi\in \cD(\bG_m, S_1)$ we denote the function
\[
\sC^0_c(S_1, S_2, \bZ)\times F_{S_0} \lra \bC, \,\, (f, x^{0}) \mapsto \int_0^{\infty}\phi(f, x_0, x^0) d^{\times} x_0
\]
by $\int_0^{\infty} \,\phi \,\, d^{\times} x_0$. It is easy to see that we have $(\int_0^{\infty} \,\phi \,\, d^{\times} x_0)(\xi f, \xi x^{0})= (\int_0^{\infty} \,\phi\,\, d^{\times} x_0)(f, x^{0})$ for all $\xi \in F^*_+$. Therefore we obtain a homomorphism
\begin{equation}
\label{eichlershimura1}
\cD(\bG_m, S_1) \lra  H^0(F^*_+,\sD(S_1, S_2, C^{\infty}(U_{\infty}^0))),\,\,\phi \mapsto \int_0^{\infty} \,\phi\,\, d^{\times}x_0.
\end{equation}
Here $C^{\infty}(U_{\infty}^0)$ denotes the ring of smooth $\bC$-valued functions on $U_{\infty}^0$ (the homeomorphism $\lambda: U_{\infty}^0\to \bR^{S_{\infty}^0} = \bR^d, (x_v)\mapsto (\log(x_v))$ provides $U_{\infty}^0$ with the structure of a real manifold). Note that $U_{\infty}^0$ carries the canonical $d$-form $d^{\times}x_1 \ldots d^{\times}x_d = \prod_{v\in S_{\infty}^0} d^{\times}x_v$ so we obtain a map
\begin{equation}
\label{eichlershimura3}
C^{\infty}(U_{\infty}^0) \lra \Omega^d(U_{\infty}^0,\bC),\,\, f\mapsto f(x_1, \ldots, x_d) \,\, d^{\times}x_1 \ldots d^{\times}x_d
\end{equation} 

Define
\begin{equation}
\label{eichlershimura2}
\cD(\bG_m, S_1)\to H^d(F^*_+, \sD(S_1, S_2,  \bC)),\,\, \phi\mapsto \kappa_{\phi}
\end{equation}
as follows. Put $C^{\bcdot} \, \colon \!\! \,= \,\,\sD(S_1, S_2, \Omega^{\bcdot}(U_{\infty}^0,\bC))$. Since $U_{\infty}^0 \approx \bR^d$, the complex $C^{\bcdot}$ is a resolution of $\sD(S_1, S_2, \bC)$ and we have $C^i= 0$ if $i>d$.  The map \eqref{eichlershimura2} is the composite of \eqref{eichlershimura1} with the map
\begin{equation}
\label{eqn:edgecohomss}
H^0(F^*_+,\sD(S_1, S_2, C^{\infty}(U_{\infty}^0)))\to H^0(F^*_+, C^d) \to H^d(F^*_+,\sD(S_1, S_2, \bC))
\end{equation}
where the first arrow is induced by \eqref{eichlershimura3} 
while the second is an edge morphism of the spectral sequence 
\[
E^{pq}_1= H^q(F^*_+, C^p) \Longrightarrow E^{p+q}= H^{p+q}(F^*_+,C^{\bcdot}) = H^{p+q}(F^*_+,\sD(S_1, S_2, \bC)).
\]

\begin{prop}
\label{prop:meascoh}
For $\phi\in \cD(\bG_m, S_1)$ and $\kappa = \kappa_{\phi}\in
H^d(F^*_+, \sD(S_1, S_2,  \bC))$ we have $\mu_{\phi} = \mu_{\kappa}$.
\end{prop}

{\em Proof.} Define a pairing
\begin{equation*}
\label{equation:period1}
\langle\,\, , \,\,\rangle: \cD(\bG_m, S_1) \times C^0(\cG_p, \bC) \lra \bC
\end {equation*}
as the composite of the product of \eqref{eichlershimura1} and \eqref{recipr1}with the map
\begin{align}
\label{align:eichlershimura4}
H^0(F^*_+,\sD(S_1, S_2, C^{\infty}(U_{\infty}^0))) \times H_0(F^*_+/E_+,H^0(E_+, \sC))\\
\stackrel{\cap}{\lra}\,\,\,  H_0(F^*_+/E_+, H^0(E_+, C^{\infty}(U_{\infty}^0))) \, \lra \,\,H_0(F^*_+/E_+, \bC) \cong \bC\nonumber
\end{align}
(where $\sC\colon \!\! = \sC^0_c(S_1, S_2, \bZ)$). Here the first map is induced by \eqref{distr2} and the second by 
\begin{equation}
\label{eichlershimura5}
H^0(E_+, C^{\infty}(U_{\infty}^0))\to \bC, \,f\mapsto \int_{U_{\infty}^0/E_+} f(x_1, \ldots, x_d) \,d^{\times}x_1 \ldots d^{\times}x_d.\end {equation}
A simple computation shows that 
\begin{equation*}
\label{equation:distrpair}
\langle \phi , f \rangle = \int_{\cG_p} f(\gamma) \mu_{\phi}(d\gamma).
\end{equation*}
for all $f\in C^0(\cG_p, \bC)$. Thus, we need to show  
\(
\kappa_{\phi}\cap \partial(f) = \langle \phi , f \rangle
\), i.e.\ we have to show that the diagram
\begin{equation}
\label{equation:cap2b}
{\small\xymatrix@+0.1pc{
H^0(F^*_+,\sD(S_1, S_2, C^{\infty}(U_{\infty}^0))) \times H_0(F^*_+/E_+, H^0(E_+, \sC)) \ar[dd]^{(\ref{eqn:edgecohomss})\times \ep}  \ar[dr]_{(\ref{align:eichlershimura4})}\\
& \bC\\
H^d(F^*_+,\sD(S_1, S_2, \bC)) \times H_d(F^*_+, \sC) \ar[ur]^{\cap}
}}
\end{equation}
commutes. For that consider the commutative diagram
\begin{equation*}
\label{equation:cap3b}
{\scriptsize
\xymatrix@-1.0pc{
H^0(F^*_+,\sD(S_1, S_2, C^{\infty}(U_{\infty}^0)) \times H_0(F^*_+/E_+,
H^0(E_+, \sC)) \ar[d]^{\id \times \eta}\ar[r] & H_0(F^*_+/E_+, H^0(E_+, C^{\infty}(U_{\infty}^0))) \ar[d]^{\eta}\\
H^0(F^*_+,\sD(S_1, S_2, C^{\infty}(U_{\infty}^0))) \times H_0(F^*_+/E_+,
H_d(E_+, \sC)) \ar[d]^{3 \times \id}\ar[r] & H_0(F^*_+/E_+, H_d(E_+, C^{\infty}(U_{\infty}^0))) \ar[d]^{4}\\
H^0(F^*_+,\sD(S_1, S_2, \Omega^d(U_{\infty}^0))) \times H_0(F^*_+/E_+,
H_d(E_+, \sC)) \ar[d]^{\id \times 5} \ar[r]& H_0(F^*_+/E_+, H_d(E_+,\Omega^d(U_{\infty}^0))) \ar[d]^{6}\\
H^0(F^*_+,\sD(S_1, S_2, \Omega^d(U_{\infty}^0))) \times H_d(F^*_+,
\sC)) \ar[d]^{7 \times \id} \ar[r] & H_d(F^*_+, \Omega^d(U_{\infty}^0)) \ar[d]^{8}\\
H^d(F^*_+,\sD(S_1, S_2, \bC)) \times H_d(F^*_+, \sC)\ar[r] & H_0(F^*_+,\bC) =\bC\\
}}
\end{equation*}
Here the horizontal maps are cap-products induced by the pairings \eqref{distr2}. The maps 3 and 4 are induced by the map \eqref{eichlershimura3}, the maps 5 and 6 are edge morphisms in a Hochschild-Serre spectral sequence and 7 and 8 are edge morphisms of a $E^1$- (resp.\ $E_1$-) hyper(co-)homology spectral sequence for the resolution $0 \to \bC \to \Omega^0(U_{\infty}^0)\to \Omega^1(U_{\infty}^0) \to \ldots$

The commutativity of (\ref{equation:cap2b}) follows once we have shown that the composition of the right column of vertical maps is induced by the map \eqref{eichlershimura5}. However this can be easily deduced from the commutativity of the obvious diagram
\begin{equation*}
\label{equation:cap4}
\begin{CD}
H^0(E_+, C^{\infty}(U_{\infty}^0)) @>>> H^0(E_+,\Omega^d(U_{\infty}^0)) @>>> H^d(E_+, \bC)\\
@VV \cap \eta V@VV \cap \eta V@VV \cap \eta V\\
H_d(E_+, C^{\infty}(U_{\infty}^0)) @>>> H_d(E_+, \Omega^d(U_{\infty}^0)) @>>> H_0(E_+, \bC)
\end{CD}
\end{equation*}
and the fact that the trace map $H_{\DR}^d(M) \to \bC, [\omega]
\mapsto \int_M \omega$ for a $d$-dimen\-sional oriented manifold $M$ corresponds under the canonical isomorphism $H_{\DR}^d(M)\cong H_{\sing}^d(M)$ to the map $x\mapsto x\cap\eta_M$ where $\eta_M$ denotes the fundamental class of $M$.\enddemo

For $\phi\in \cD(\bG_m, S_1)$ put $\phi_0 = \phi_{U_0}$ where $U_0 = \prod_{\fp\in S_1} \cO_{\fp}\times  \prod_{\fp\in S_2} \cO_{\fp}^*$. 

\begin{corollary}
\label{corollary:Lwert0} 
For $\fp\in S_1$ let $c_{\fp} = c_{\ord_{\fp}}\in H^1(F^*_+, C_c^0(F_{\fp}, \bZ))$ be the cohomology class of the 1-cocyle \eqref{cycle} for $\ell = \ord_{\fp}$. Then we have
\[
\int_{\bI/F^*}\phi_0(x)\, d^{\times}x \,\, =\,\, (-1)^{\binom{r}{2}}\,\, (\kappa_+ \cup c_{\fp_1}\cup \ldots \cup c_{\fp_{r}}) \cap \vartheta.
\]
Here the cup-product is induced by \eqref{meas2a} and the cap-product by \eqref{meas2}. 
\end{corollary}

{\em Proof.} We use the notation of Prop.\ \ref{prop:thetaord}. Note that $\cF_1\subseteq \bI^{p,\infty}/U^{p,\infty}$ is a finite set and $\cX = U_0\times \cF_1$ so we have $1_{\cX} = \sum_{x\in \cF} x 1_{\cX_0}$ (where $\cF \colon \!\!= \{1\}\times \cF_1,\cX_0 \colon \!\!= U_0\times \{1\}\subseteq F_p^*\times \bI^{p,\infty}/U^{p,\infty} = \bI^{\infty}/U^{p,\infty}$). Because of the commutativity of (\ref{equation:cap2b}) and \ref{prop:thetaord} it is enough to prove that the pair $(\int_0^{\infty} \,\phi\,\, d^{\times}x_0, [1_{\cX}])$ is mapped to $\int_{\bI/F^*}\phi_0(x) d^{\times}x$ under the pairing (\ref{align:eichlershimura4}). In fact by definition of (\ref{align:eichlershimura4}) the pair is mapped to
\begin{eqnarray*}
&& \int_{U_{\infty}/E_+}\, \phi(1_{\cX}, x_0, \ldots, x_d) \,d^{\times}x_0\ldots d^{\times}x_d\\
&& \sum_{x\in \cF}\, \int_{U_{\infty}/E_+}\, \phi(x 1_{\cX_0}, x_0, \ldots, x_d) \,d^{\times}x_0\ldots d^{\times}x_d\\ 
&& = \,\, \int_{\bR_+\times \cE\times \cF}\, \phi_0(x) \,d^{\times}x \,\, = \,\, \int_{\bR_+\times \cE\times U_p\times \cF}\, \phi_0(x)\, d^{\times}x \\
&& = \,\, \int_{\bI/F^*} \phi_0(x)\, d^{\times}x
\end{eqnarray*}
 where $\cE\subseteq U^0_{\infty}$ is a fundamental domain for the action $E_+$. \enddemo

\section{$p$-adic $L$-functions of Hilbert modular forms}
\label{section:pmeashmf}

\subsection{$p$-ordinary cuspidal automorphic representations}
Let $\pi = \otimes_v \, \pi_v$ be a unitary cuspidal automorphic
representation of $G(\bA)$.  Thus
$\pi$ is an irreducible direct summand of the right regular
representation of $G(\bA)$ in $L^2_{\disc}(G(F)\backslash G(\bA))$. If the archimedian local representations $\pi_v$ are
discrete series then a $p$-adic $L$-function $L_p(s, \pi)$ for $\pi$
can be defined. The first construction under certain restrictions on $\pi$ 
is due to Manin \cite{manin}; a construction in most generality (based on earlier work of Panchishkin \cite{panchishkin}) is due to Dabrowski \cite{dabrowski}; see (\cite{dabrowski}, sect.\ 12) for further references.

In this section we shall give another definition of $L_p(s, \pi)$ well-suited for the proof of the weak exceptional zero conjecture. We assume that $\pi$ has parallel weight $(2,\ldots,2)$ and is $p$-ordinary. The first condition means that $\pi_v = \cD(2)$ is the discrete series representation of $G(\bR)$ of lowest weight 2 for all $v\in S_{\infty}$ and the second that $\pi_{\fp}$ is ordinary for all $\fp\in S_p$ in the sense of section \ref{subsection:ordinary}.
We shall attach an element $\phi_{\pi}\in \cD(\bG_m, S_1)$ to $\pi$, show that the corresponding cohomology class $\kappa_{\pi} = \kappa_{\phi_\pi} $ is integral and define $L_p(s, \pi)$ as the $p$-adic $L$-series associated to $\kappa_{\pi}$ (here $S_1$ denotes the set of $\fp\in S_p$ with $\pi_{\fp} = \St$). 

We introduce some notation. Firstly, we denote by $\fA_0(G, \uzwei)$
the set of all unitary cuspidal automorphic representation $\pi$ of
$G(\bA)$ of parallel weight $(2,\ldots,2)$. For each $\fp\in S_p$ we
fix an ordinary parameter $\alpha_{\fp}\in {\barcO}^*$ and put $a_{\fp} =
\alpha_{\fp} + \Norm(\fp)/\alpha_{\fp}$. As before we let $m =
\sharp(S_p)$ and $r$ be the number of $\fp\in S_p$ with $\alpha_{\fp} =1$. 
We choose an ordering $\fp_1, \ldots, \fp_m$ of the
places above $p$ so that $\alpha_{\fp_1} = \ldots = \alpha_{\fp_r} = 1$.  
We write $F_i$, $\alpha_i$ and $a_i$
instead of $F_{\fp_i}$, $\alpha_{\fp_i}$ and $a_{\fp_i}$  and put 
$\ualpha = (\alpha_1, \ldots, \alpha_m)$ and $\ua = (a_1, \ldots, a_m)$. Moreover we denote by $\fA_0(G,
\uzwei, \ualpha)$ the subset of $\pi\in \fA_0(G, \uzwei)$ such that $\pi_{\fp_i} = \pi_{\alpha_i}$ for $i= 1,\ldots, m$. 

For $\pi\in \fA_0(G, \uzwei, \ualpha)$ and a finite set of places $S$ of $F$ we put $\pi_S = \otimes_{v\in S} \, \pi_v$ and $\pi^S = \otimes_{v\not\in S} \, \pi_v$. We also write $\pi_{\infty}, \pi_{p, \infty}, \pi^{\infty}$ etc.\ for $\pi_{S_{\infty}}, \pi_{S_p \cup S_{\infty}}, \pi^{S_{\infty}}$ etc. For each finite place $\fq$ of $F$ we denote by $\ff(\pi_{\fq})$ the conductor of $\pi_{\fq}$ and we set $\ff(\pi) \colon \!\! = \prod_{\fq} \, \ff(\pi_{\fq})$. Thus the multiplicity $\ord_{\fp}(\ff(\pi))$ of $\fp\in S_p$ in $\ff(\pi)$ is $=1$ if $\alpha_{\fp} = \pm 1$ and $=0$ otherwise.

\subsection{Adelic Hilbert modular forms}

In section \ref{subsection:pmeashmf} we shall define a certain element $\phi_{\pi}\in \cD(\bG_m, S_1)$ for $\pi\in \fA_0(G, \uzwei)$. Firstly, we need to recall the notion of a Hilbert modular cusp form of parallel weight $(2,\ldots,2)$ (in the adelic setting). It is a function $\Phi: G(\bA) \to \bC$ with the following properties
\begin{itemize}
\label{itemize:cuspform}
\item[(i)]
For $\gamma \in G(F)$ we have $\Phi(\gamma g) = \Phi(g)$. 
\item[(ii)]
For any $g\in G(\bA)$, $k_{\infty}\in K_{\infty}^+$ we have $\Phi(g k_{\infty}) = j(k_{\infty}, \ui)^{-2} \Phi(g)$.
\item[(iii)]
For any $g \in G(\bA^{\infty})$, $\uz\in \bH^{d+1}$ define $\bff_{\Phi}(\uz, g) \colon \!\! = j(g_{\infty}, \ui)^{2}\Phi(g_{\infty}, g)$ where $g_{\infty}\in G(F_{\infty})^{+}$ is such that $g_{\infty} \ui = \uz$ (by (ii) $\bff_{\Phi}(\uz,g)$ is well defined). Then $\uz \mapsto \bff_{\Phi}(\uz, g)$ is a holomorphic function on $\bH^{d+1}$.
\item[(iv)] There exists a compact open subgroup $K$ of $G(\bA^{\infty})$ such that $\Phi(g k) \break = \Phi(g)$ for all $k\in K$ and $g\in G(\bA)$.
\item[(v)] (Cuspidality) For any $g\in G(\bA)$ we have
\[
\int_{\bA/F}\, \Phi\left(\left(\begin{matrix} 1 & x\\ 0 & 1\end{matrix}\right) g\right) dx = 0.
\]
\end{itemize}

We denote by $\cA_0(G, \hol, \uzwei)$ the space of all functions $\Phi$ satisfying (i) -- (v) above.  It is a left $G(\bA^{\infty})$-module via the right action on $G(\bA)$. For a compact open subgroup $K\subseteq G(\bA^{\infty})$ we set $S_2(G,K) = \cA_0(G, \hol, \uzwei)^K$. If $K = K_0(\fn)$ for an ideal $\fn$ of $\cO_F$, we write $S_2(G, \fn)$ instead of $S_2(G, K_0(\fn))$.

Let $\Phi\in \cA_0(G, \hol, \uzwei)$ and let $\bff_{\Phi}$ be as in (iii) above.
We define 
\begin{equation}
\label{intz0}
\int_{\sigma_0(Q)}^{\sigma_0(P)}\, \bff_{\Phi}(z_0, z_1, \ldots, z_d, g) dz_0 
\end{equation}
by integrating the function $z_0\mapsto \bff_{\Phi}(z_0, z_1, \ldots, z_d, g)$ along the geodesic in $\bH$ from $\sigma_0(Q)$ to $\sigma_0(P)$. Using the well-known fact that $\bff_{\Phi}(\uz, g)$ for fixed $g\in G(\bA^{\infty})$ rapidly decreases at the set of cusps $\bP^1(F)$ of $\bH^{d+1}$ it is easy to see that \eqref{intz0} is well-defined and that
\[
\int_{\sigma_0(Q)}^{\sigma_0(P)}\, \bff_{\Phi}(z_0, \ldots ) dz_0 + \int_{\sigma_0(P)}^{\sigma_0(R)}\, \bff_{\Phi}(z_0, \ldots ) dz_0 = \int_{\sigma_0(Q)}^{\sigma_0(R)}\, \bff_{\Phi}(z_0, \ldots ) dz_0
\]
for all $P,Q, R\in \bP^1(F)$. 

Let $\Div(\bP^1(F))$ be the free abelian group over $\bP^1(F)$ and let $\cM= \Div_0(\bP^1(F))$ be the subgroup of elements $\sum_{j=1}^r \, m_i P_i\in \Div(\bP^1(F))$ with $\deg(\sum_{j=1}^r \, m_i P_i) = \sum_{j=1}^r \, m_i = 0$. The natural $G(F)$-action on $\bP^1(F)$ induces a $G(F)$-action on $\cM$. For $g\in G(\bA^{\infty})$ we obtain a homomorphism 
\begin{equation*}
\label{intdiv1}
\cM \lra \cO_{\hol}(\bH^d),\,\,\, m\mapsto \int_{\sigma_0(m)}\, \bff_{\Phi}(z_0, z_1, \ldots, z_d, g) dz_0
\end{equation*}
 which coincides with \eqref{intz0} for $m= P- Q$. For $m\in \cM$ we define
\begin{equation*}
\label{intdiv2}
\int_m \omega_{\Phi}(g) \,=\, \left(\int_{\sigma_0(m)}\, \bff_{\Phi}(z_0, z_1, \ldots, z_d, g) dz_0\right) dz_1 \ldots dz_d \in \Omega_{\hol}^d(\bH^d).
\end{equation*}
We let the group $G(F)^+$ act on $\bH^d$ via the embedding $G(F)^+\to (G(\bR)^+)^d, \break \gamma \mapsto (\sigma_1(\gamma), \ldots, \sigma_d(\gamma))$. A simple computation using (i) shows
\begin{equation}
\label{intginv}
\gamma^*\left(\int_{\gamma m} \omega_{\Phi}(\gamma g)\right) = \int_m \omega_{\Phi}(g)
\end{equation}
for all $\gamma \in G(F)^+$, $g\in G(\bA^{\infty})$ and $m\in
\cM$. The integral $\int_m \omega_{\Phi}(g)$ will be used 
in the construction of the {\it Eichler-Shimura map}
\eqref{eichshiml2} in section \ref{subsection:eichshim}.

\begin{definition}
\label{definition:auttree}
(a) We denote by $\cA_0(G, \hol, \uzwei, \ua)$ the $\bC$-vector space of maps $\Phi: G(\bA^p)\to \sB^{\ua}(F_p, \bC)$ such that 
\begin{itemize}
\label{itemize:cuspform1}
\item[(i)]
There exists a compact open subgroup $K$ of $G(\bA^{p,\infty})$ such that $\Phi(g k) \break = \Phi(g)$ for all $k\in K$ and $g\in G(\bA^p)$.
\item[(ii)]
For $\psi\in \sB_{\ua}(F_p,\bC)$ the map 
\[
\quad\langle \Phi, \psi\rangle: G(\bA) = G(F_p) \times G(\bA^p)\mapsto \bC, (g_p, g^p) \mapsto \langle g_p\cdot\psi, \Phi(g^p) \rangle
\]
lies in $\cA_0(G, \hol, \uzwei)$. 
\end{itemize}
$\cA_0(G, \hol, \uzwei, \ua)$ is a left $G(\bA^{p,\infty})$-module via the right action on $G(\bA^{p,\infty})$.

\noi (b) For a compact open subgroup $K\subseteq G(\bA^{p,\infty})$ we set $S_2(G, K,\ua) = \cA_0(G, \hol, \uzwei, \ua)^K$. If $K = K_0(\fm)^p$ where $\fm$ is an ideal of $\cO_F$ which is relatively prime to $p\cO_F$, we shall write $S_2(G, \fm, \ua)$ instead of $S_2(G, K_0(\fm)^p,\ua)$.
\end{definition}

\begin{remarks}
\label{remarks:adeliccuspforms}\rm
(a) Let $\pi\in \fA_0(G, \uzwei, \ualpha)$. It is easy to see that 
\begin{eqnarray*}
\cA_0(G, \hol, \uzwei) & \cong & \Hom_{G(F_{\infty})}\left(\pi_{\infty} \,, \,\,L_0^2(G(F)\backslash G(\bA))\right), \\
\cA_0(G, \hol, \uzwei, \ua) & \cong & \Hom_{G(F_p\times F_{\infty})}\left( 
\pi_{p, \infty}\, , \,\,L_0^2(G(F)\backslash G(\bA))\right)
\end{eqnarray*}
as representations of $G(\bA^{\infty})$ and $G(\bA^{p, \infty})$
respectively.

\noi (b) Assume that $F$ has narrow class number 1. Let
$\fn$ be a non-zero ideal of $\cO_F$ and let $\Gamma_0(\fn)$ be the
subgroup of matrices $A\in \SL_2(\cO)$ which are upper triangular
modulo $\fn$. Then $S_2(G, \fn)$ can be identified with the space
$S_2(\Gamma_0(\fn))$ of usual Hilbert modular cusp forms of parallel
weight $(2, \ldots, 2)$ on $\Gamma_0(\fn)$. Moreover if the ideal
$\fm$ of $\cO_F$ is relatively prime to $p\cO_F$ and if $\fn$ is the
product of $\fm$ with all $\fp\in S_p$ with $\alpha_{\fp} = \pm 1$
then $S_2(G, \fm, \ua)$ can be identified with the subspace of
$f\in S_2(\Gamma_0(\fn))$ which satisfy (i) $T_{\fp} f = a_{\fp} f$
for all $\fp\in S_p$ with $\alpha_{\fp} \ne \pm 1$, (ii) $f$ is
$\fp$-new and $U_{\fp} f = \alpha_{\fp} f$ for all $\fp\in S_p$ with
$\alpha_{\fp} = \pm 1$. Here, if $\fp \nmid \fn$ (resp.\ $\fp \mid
\fn$) $T_{\fp}$ (resp.\ $U_{\fp}$) denotes the Hecke operator at
$\fp$. 
\end{remarks}

\subsection{Hecke Algebra} 

We recall here a few facts about the Hecke algebra of $G(\bA^S)$ (a reference for what follows is e.g.\ (\cite{bump}, 3.4 and 4.2) or (\cite{bushhen}, 1.2--4)). Fix a finite set of places $S$ of $F$ containing $S_{\infty}$. Let $dg$ denote the Haar measure on $G(\bA^S)$ normalized such that $\int_K dg =1$ for $K = \prod_{v\not\in S}\, G(\cO_v)$. For a field $\Omega$ of characteristic zero we denote by $\cH_{\Omega}^S = \cH_{G(\bA^S)}$ the Hecke algebra of $G(\bA^S)$ with coefficients in $\Omega$, i.e.\ it is the convolution ring of locally constant compactly supported $\Omega$-valued functions on $G(\bA^S)$ (see \cite{bump}, p.\ 309). If $K\subseteq G(\bA^S)$ is any compact open subgroup then we let $\cH^S_{K,\Omega}$ be the subspace of $K$-biinvariant functions in $\cH_R^S$. Let $\fq\in \bfP_F-S$ and assume $K = K' \times G(\cO_{\fq})$ for some compact open subgroup $K'$ of $G(\bA^{S'})$ where $S' = S\cup\{\fq\}$. Then $\cH^S_{K,\Omega}$ is isomorphic to the tensor product of $\cH^{S'}_{K',\Omega}$ and the Hecke algebra $\cH_{\Omega}(G(F_{\fq}), G(\cO_{\fq}))$. For the latter we have $\cH_{\Omega}(G(F_{\fq}), G(\cO_{\fq}))\cong \Omega[T_{\fp}]$  (see \cite{bump}, 4.6.5) so in this case $\cH^S_{K,\Omega} = \cH^{S'}_{K',\Omega}[T_{\fq}]$.

Recall that the concepts ``smooth $\cH_{\Omega}^S$-module'' and ``smooth $\Omega$-repre\-sentation of $G(\bA^S)$'' are interchangeable. In the following we view a smooth $\cH_{\Omega}^S$-module often as a smooth $G(\bA^S)$-module and vice versa. A sequence of smooth $\cH_{\Omega}^S$-modules $V_1 \to V_2 \to V_3$ is exact if and only if $V_1^K \to V_2^K \to V_3^K$ is exact for all compact open subgroups $K$ of $G(\bA^S)$. We call a smooth representation $V$ of $G(\bA^S)$ semisimple if it is isomorphic to a direct sum of smooth irreducible representations of $G(\bA^S)$. A smooth representation $V$ of $G(\bA^S)$ is irreducible if and only if $V^K$ is either zero or a simple $\cH^S_{K,\Omega}$-module for all $K$. More generally it is easy to see that a smooth representation $V$ of $G(\bA^S)$ is semisimple if and only if $V^K$ is a semisimple $\cH^S_{K,\Omega}$-module for all $K$. 

Let $V$ and $W$ be smooth $\cH_{\Omega}^S$-modules and assume that $V$ is irreducible and $W$ is semisimple. Let $K_0$ be a compact open subgroup of $G(\bA^S)$ such that $V^{K_0} \ne 0$. Then the canonical map $\Hom_{G(\bA^S)}(V, W) \to \Hom_{\cH^S_{K_0}}(V^{K_0}, \break W^{K_0})$ is an isomorphism. For that it is enough to assume that $W$ is irreducible in which case the assertion follows from (\cite{bushhen}, Prop.\ on p.\ 38).

For $\pi\in \fA_0(G, \uzwei)$ with the complex representation $\pi^S$ of $G(\bA^S)$ is an example of a smooth irreducible representation. It is also known that $\pi^S$ can be defined over a finite extension of $\bQ$. More precisely there exists a smallest finite extension $\Omega = \Omega_{\pi}\subseteq \barQ$ of $\bQ$ (the {\it field of definition of $\pi$}) and a smooth irreducible $\Omega$-representation $G(\bA^S)\to \GL(V)$ such that $\pi^S \cong V\otimes_{\Omega} \bC$ as $G(\bA^S)$-representations. By abuse of notation we also write $\pi^S$ (resp.\ $\pi^{S,K}$) for the $\cH_{\Omega}^S$-module $V$ (resp.\ for the $\cH^S_{K}$-module $V^K$). For a field $\bfC$ containing $\Omega_{\pi}$, a compact open subgroup $K$ of $G(\bA^S)$ and a smooth semisimple $\bfC$-representation $W$ of $G(\bA^S)$ we write $W_{\pi}$ for $\Hom_{G(\bA^S)}(\pi^S, W) = \Hom_{G(\bA^S)}(V\otimes_{\Omega} \bfC, W)$ and $W_{\pi}^K$ for $\Hom_{\cH^S_{K}}(V^K, W^K)$. Also if $f: W' \to W$ is a homomorphism of smooth semisimple $G(\bA^S)$-representations we denote the induced homomorphism $W'_{\pi} \to W_{\pi}$ of $\bfC$-vector spaces by $f_{\pi}$. We have $W_{\pi}^K = W_{\pi}$ if $K\subseteq K_0(\ff(\pi))^S$ and $W_{\pi}^K = 0$ otherwise. If $K = K_0(\ff(\pi))^S$, then $\pi^{S,K}$ is onedimensional as a $\Omega_{\pi}$-vector space \cite{casselman}. Thus the action of $\cH^S_{K}$ is given by a homomorphism $\lambda_{\pi}: \cH^S_{K}\to \Omega_{\pi}$ (it is known that $\lambda_{\pi}(T_{\fq})$ lies in the ring of integers of $\Omega_{\pi}$). In this case we have $W_{\pi} = \{w\in W^{K}\mid\, t w = \lambda_{\pi}(t) w \,\, \forall t\in \cH^S_{K}\}$. We also remark that if $W'\to W \to W''$ is an exact sequence of smooth semisimple representations of $G(\bA^S)$ then $W'_{\pi}\to W_{\pi} \to W''_{\pi}$ is exact as well.

Finally, a representation $W$ of $G(\bA^S)$ will be said to be of
{\it automorphic type} if $W$ is smooth and semisimple and the only
irreducible subrepresentations of $W$ are either the onedimensional representations or the representations $\pi^S$ for $\pi\in \fA_0(G, \uzwei)$. By strong multiplicity one, if $W$ is of automorphic type then $W_{\pi}$ is independent of the set $S$ in the following sense. Let $S' \supset S$ and $K = \prod_{v\in S'-S} \, K(\ff(\pi))_v$. Then we have $\Hom_{G(\bA^{S'})}(\pi^{S'}, W^{K})= \Hom_{G(\bA^S)}(\pi^S, W)$ (this fact will be used in the proof of Prop.\ \ref{prop:shimuraisomorphism} below).

\subsection{Cohomology of $\GL_2(F)$}
\label{subsection:finiteautomorph}

In this section we introduce and study the cohomology groups of certain $G(F)^+$-modules $\sA(\ua, \cM, \bC)$ on which the group $G(\bA^{p,\infty})$ acts smoothly and so that each $\pi\in \fA_0(G, \uzwei, \ua)$ occurs with multiplicity $2^{d+1}$ in $H^d(G(F)^+, \sA(\ua, \cM, \bC))$ (see Prop.\ \ref{prop:shimuraisomorphism} below). 

Let $H\subset G(F)$ be any subgroup and let $M$ be a left $G(F)$-module. Let $R$ is a ring and $N$ an $R[H]$-module. For a finite subset $S$ of $\bfP_F^{\infty}$ we denote by $\sA(S, M; N)$ the $R$-module of maps $\Phi: G(\bA^{S, \infty})\times M \to N$ such that there exists a compact open subgroup $K$ of $G(\bA^{S, \infty})$ with $\Phi(g k, m) = \Phi(g, m)$ for all $k\in K$, $g \in  G(\bA^{S, \infty})$ and $m\in M$. 

We have commuting $G(\bA^{S, \infty})$- and $H$-actions on $\sA(S, M; N)$; the first is induced by right multiplication on $G(\bA^{S, \infty})$ and the second is given by $(\gamma \wcdot \Phi)(g, m) =\gamma \Phi(\gamma^{-1} g, \gamma^{-1}m)$. For a compact open
subgroup $K\subseteq G(\bA^{S, \infty})$ we set $\sA(K, S, M; N) = \sA(S, M; N)^K$. If $K = K_0(\fm)^S$ for an ideal $\fm$ of $\cO_F$ not divisible by any $\fq\in S$ we write $\sA(\fm, S , M; N)$ for $\sA(K_0(\fm)^S, S , M; N)$. If $S = \emptyset$ we will often drop $S$ from the notation, i.e.\ we write $\sA(M; N)$, $\sA(K, M; N)$ etc.\  for $\sA(S, M; N)$, $\sA(K, S, M; N)$ etc.

In contrast to our previous notation in this section we denote by $S_1$, $S_2$ subsets of $S_p$ with  $S_1 \subseteq S_2\subseteq S_p$. We define the $G(\bA^{S_2, \infty})$-$H$-module $\sA(\ua_{S_1}, S_2, M; N)$ by
\begin{equation*}
\label{aftree}
\sA(\ua_{S_1}, S_2, M; N) =  \sA(S_2, M; \sB^{\ua}(F_{S_1}, N))
\end{equation*}
Also for a compact open subgroup $K\subseteq G(\bA^{S_2, \infty})$ we put 
$\sA(K, \ua_{S_1}, S_2, M;\break N) = \sA(\ua_{S_1}, S_2, M; N)^K$ and if $K = K_0(\fm)^{S_2}$ for an ideal $\fm$ of $\cO_F$ not divisible by any $\fp\in S_2$ we set $\sA(\fm, \ua_{S_1}, S_2, M; N)$ for $\sA(K_0(\fm)^p, \ua_{S_1}, S_2,\break M; N)$. If $S= S_1 = S_2$ (resp.\ $S_1 = S_2 = S_p$) (we deal mostly with the latter case) we shall drop $S_2$ (resp.\ $S_1$ and $S_2$) again from the notation, i.e.\ we put
\begin{eqnarray*}
\sA(\ua_S, M; N) & = & \sA(\ua_S, S, M; N) \,\, =\,\, \sA(S, M; \sB^{\ua}(F_S, N)),\\
\sA(\ua, M; N) & = & \sA(\ua_{S_p}, S_p, M; N) \,\, =\,\, \sA(S_p, M; \sB^{\ua}(F_p, N)).
\end{eqnarray*}
So $\sA(\ua, M; N)$ can be identified with the $R$-module of maps $\Phi: G(\bA^{p,\infty})\times M \times \sB_{\ua}(F_p, R) \to N$ which are invariant under some compact open subgroup of $G(\bA^{p,\infty})$. 

The pairing \eqref{upairdual} induces a pairing 
\begin{equation*}
\label{eqn:phichar2}
\langle \,\,\, , \,\,\,\rangle: \sA(\ua_{S_1}, S_2, M; N) \times \sB_{\ua}(F_{S_1}, R) \lra \sA(S_2, M; N)
\end{equation*}
hence a homomorphism $\sA(\ua_{S_1}, S_2, M; N) \to \sB^{\ua}(F_{S_1}, \sA(S_2, M; N))$ which becomes an isomorphism when restricting it to $K$-invariant elements  
\begin{equation*}
\label{eqn:vertauschac1}
\sA(K, \ua_{S_1}, S_2, M; N) \stackrel{\cong}{\lra} \sB^{\ua}(F_{S_1}, \sA(K, S_2, M; N))
\end{equation*}
(for any compact open subgroup $K\subseteq G(\bA^{S_2, \infty})$). Similarly, for any $\fp\in S_1$ and $S_0 \colon \!\! = S_1-\{\fp\}$ we have an isomorphism
\begin{equation}
\label{vertauschac2}
\sA(K, \ua_{S_1}, S_2, M; N) \stackrel{\cong}{\lra} \sB^{a_{\fp}}(F_{\fp}, \sA(K, \ua_{S_0}, S_2, M; N))
\end{equation}

\begin{remark}
\label{remark:afcofree}
\rm Assume that  $M$ is free as an abelian group. For a compact open subgroup $K\subseteq G(\bA^{S_2, \infty})$ we have 
\begin{eqnarray*}
\label{eqn:rmodstr}
\sA(K, \ua_{S_1}, S_2, M; R) & \cong & \Coind^{G(\bA^{S_2, \infty})}_K \Hom_{\bZ}(M, \sB^{\ua}(F_{S_1}, R))\\
& \cong & \Coind^{G(\bA^{S_2, \infty})}_K \Hom_R(M\otimes_{\bZ} \sB_{\ua}(F_{S_1}, R), R)
\end{eqnarray*}
Since $\sB_{\ua}(F_{S_1}, R)$ is a free $R$-module we see that the $R$-module $\sA(K, \ua_{S_1}, S_2,\break  M; R)$ is isomorphic to a product of copies of $R$. 

Assume now that $R= \bC_p$, let $\cO=\cO_{\bC_p}$ so that $L=\sA(K, \ua_{S_0}, S_2, M; \cO)$ is a complete lattice in $V=\sA(K, \ua_{S_0}, S_2, M; \cO)\otimes_{\cO}\bC_p$ (see section \ref{section:distmeas}). Let $\fp\in S_1$ such that $\alpha_{\fp} =1$ and put $S_0= S_1- \{\fp\}$. By \eqref{vertauschac2} the vector space $\sA(K, \ua_{S_1}, S_2, M; \cO)\otimes_{\cO}\bC_p$ can be identified with $\Meas(F_{\fp}, V)$ and so the evaluation map \eqref{distr1} extends to a pairing $\Meas(F_{\fp}, V)\times \Cd(F_{\fp}, \bC_p) \lra V$, i.e.\ we have a canonical bilinear map
\begin{equation*}
\label{vertauschac4}
\sA(K, \ua_{S_1}, S_2, M; \cO)\otimes_{\cO}\bC_p \times \St(F_{\fp}, \bC_p) \to \sA(K, \ua_{S_0}, S_2, M; \cO)\otimes_{\cO}\bC_p.
\end{equation*}
It will be used for the construction of the $\cL$-invariants $\cL_{\fp}(\pi)$ in the next chapter.
\end{remark}

The next result follows immediately from \eqref{heckeinduced3} and \eqref{heckeinduced4}.

\begin{lemma}
\label{lemma:ordinaryaut}
Let $S\subseteq S_p$, let $\fp\in S$ and let $S_0 \colon \!\! = S-\{\fp\}$. Let $K\subseteq G(\bA^{S, \infty})$ be a compact open subgroup and put $K_0 = K \times G(\cO_{\fp})$, $K_1 = K \times K_0(\fp)_{\fp}$.

\noi (a) If $\alpha_{\fp} \ne \pm 1$ then the following sequence is exact
\begin{equation*}
\label{tree1}
\xymatrix@-1.0pc{
0 \ar[r] & \sA(K, \ua_{S},\! M; N) \ar[r] & \sA(K_0, \ua_{S_0},\! M; N) \ar[rr]^{(T_{\fp} - a_{\fp})} & &\sA(K_0, \ua_{S_0},\! M; N) \ar[r] & 0 }
\end{equation*}
\noi (b) If $\alpha_{\fp} = \pm 1$ then there exists a short exact sequence
\begin{equation*}
\label{tree2}
\xymatrix@-1.2pc{
0 \ar[r] & \sA(K, \ua_{S}, \! M;N) \ar[r] & \sA(K_1, \ua_{S_0}, \! M; N)^{W=\mp1} \ar[r] &\sA(K_0, \ua_{S_0},\! M; N) \ar[r] & 0}
\end{equation*}
where $W = W_{\fp}$ is a certain involution acting on $\sA(K_1, \ua_{S_0}, M; N)$.
 \end{lemma}

\begin{remark}
\label{remark:atkinlehnerinvolutions}
\rm The involution in part (b) above induces an involution -- also denoted by $W_{\fp}$ -- on the cohomology groups $H^{\bu}(G(F)^+, \sA(K_1, \ua_{S_0}, M; N))$. In particular 
if $\fn$ is an ideal of $\cO_F$ such that $\fp$ devides $\fn$ exactly once we have an involution $W_{\fp}$ acting on $H^{d+1}(G(F)^+, \sA(\fn, R))$. This is the analogue of the Atkin-Lehner involution. As in the classical case we have $W_{\fp} = -U_{\fp}$ on $H^{d+1}(G(F)^+, \sA(\fn, R))$.
\end{remark}

\begin{prop}
\label{prop:basechange}
Let $S_1\subseteq S_2 \subseteq S_p$ and let $K$ be a compact open subgroup of $G(\bA^{S_2, \infty})$. 

\noi (a) Let $N$ be a flat $R$-module (equipped with the trivial $G(F)$-action). Then the canonical map 
\begin{equation*}
\label{basechange1}
H^q(G(F)^+\!\!, \sA(K, \ua_{S_1}, S_2, \cM; R))\otimes_R N \to  H^q(G(F)^+\!\! \sA(K, \ua_{S_1}, S_2, \cM; N))
\end{equation*}
is an isomorphism for all $q\ge 0$.

\noi (b) If $R$ is noetherian then $H^q(G(F)^+, \sA(K, \ua_{S_1}, S_2, \cM; R))$ is a finitely generated $R$-module.
\end{prop}

{\em Proof.}  (a) The sequence $0\to \cM \to \Div(\bP^1(F))\cong \Ind^{G(F)}_{B(F)}\bZ \to \bZ\to 0$ yields a short exact sequence
\begin{eqnarray}
\label{div}
&& 0 \lra \sA(K, \ua_{S_1}, S_2; N) \lra \Coind^{G(F)^+}_{B(F)^+}\sA(K, \ua_{S_1}, S_2; N)\\
&& \hspace{3cm}\lra \sA(K, \ua_{S_1}, S_2, \cM; N) \lra 0\nonumber
\end{eqnarray}
(where $\sA(K, \ua_{S_1}, S_2; N) \colon \!\! = \sA(K, \ua_{S_1}, S_2, \bZ; N)$). By considering the associated long exact cohomology sequences it is enough to prove the assertion when we replace the coefficients $\sA(K, \ua_{S_1}, S_2, \cM;\cdot)$
with $\sA(K, \ua_{S_1}, S_2;\cdot)$ or $\Coind^{G(F)^+}_{B(F)^+}\sA(K, \ua_{S_1}, S_2;\cdot)$ (here we use the flatness assumption). Furthermore using Lemma \ref{lemma:ordinaryaut} it is enough to consider the case $S_1= \emptyset$, $S=S_2$. Since $\sA(K, S, \bZ; N) \cong \Coind^{G(\bA^{S, \infty})}_K N$ it suffices to show that
\begin{eqnarray*}
\label{basechange2}
H^q(G(F)^+, \Coind^{G(\bA^{S, \infty})}_K R)\otimes_R N\lra H^q(G(F)^+, \Coind^{G(\bA^{S, \infty})}_K N)\\
\label{basechange3}
H^q(B(F)^+, \Coind^{G(\bA^{S, \infty})}_K R)\otimes_R N \lra H^q(B(F)^+, \Coind^{G(\bA^{S, \infty})}_K N)
\end{eqnarray*}
are isomorphism for all $q\ge 0$ and all $R$-modules $N$. For that it is enough to prove that the functors $N \mapsto H^q(G(F)^+, \Coind^{G(\bA^{S, \infty})}_K N)$ and $N \mapsto H^q(B(F)^+, \Coind^{G(\bA^{S, \infty})}_K N )$ commute with direct limits (since any module is the direct limit of free modules of finite rank). For $g\in G(\bA^{S, \infty})$ put $\Gamma_g = G(F)^{+}\cap g K g^{-1}$. By the strong approximation theorem there are only finitely many double cosets $G(F)^+ g K$ in $G(\bA^{S, \infty})$. If $g_1, \ldots, g_n\in G(\bA^{S, \infty})$ is a system of representatives then 
\[
H^q(G(F)^+, \Coind^{G(\bA^{S, \infty})}_K N)\,\,\, = \,\,\,\bigoplus_{i=1}^n \, H^q(\Gamma_{g_i}, N).
\]
Since the group $\Gamma_g$ is $S$-arithmetic, hence of type (VFL), the functor $N \mapsto H^q(\Gamma_g, N)$ commutes with direct limits (see \cite{serre}, p.\ 101). The same proof works for $H^q(B(F)^+, \Coind^{G(\bA^{S, \infty})}_K N)$ as well. Indeed using the Iwasawa decomposition $G(\bA^{S, \infty}) = B(\bA^{S, \infty})\prod_{v\nmid \infty} G(\cO_v)$ one can easily see that $B(F)^+\backslash G(\bA^{S, \infty})/K$ is finite.

\noi (b) can be deduced using similar arguments and the fact that the groups $\Gamma_g$ are $S$-arithmetic and (\cite{serre}, remarque on p.\ 101).\enddemo

Let $S_1 \subset S_2 \subseteq S_p$ be as before, let $R$ be a ring and let $M$ be a left $G(F)$-module. We define 
\begin{equation*}
\label{cohom}
H^q_*(G(F)^+, \sA(\ua_{S_1}, S_2, M; R)) = \underset{\lra}{\lim}  \,\,H^q(G(F)^+, \sA(K, \ua_{S_1}, S_2, M; R))
\end{equation*}
where $K$ runs through all compact open subgroups of $G(\bA^{S_2, \infty})$. By \ref{prop:basechange} we have

\begin{corollary}
\label{corollary:basechange}
Let $S\subseteq S_p$ and let $R\to R'$ be a flat ring homomorphism. Then the canonical map 
\begin{equation*}
\label{basechange}
H^q_*(G(F)^+, \sA(\ua_{S}, \cM; R))\otimes_R R' \lra H^q_*(G(F)^+, \sA(\ua_{S}, \cM; R'))
\end{equation*}
is an isomorphism for all $q\ge 0$.
\end{corollary}

If $R= \bfC$ is a field of characteristic zero then $H^q_*(G(F)^+, \sA(\ua_{S_1}, S_2, M;\break \bfC))$ is a smooth $G(\bA^{S_2, \infty})$-module and it is easy to see that we have 
\[
H^q_*(G(F)^+, \sA(\ua_{S_1}, S_2, M; \bfC))^K = H^q(G(F)^+, \sA(K, \ua_{S_1}, S_2, M; \bfC)).
\]
We identify $G(F)/G(F)^+$ with the group $\Sigma = \{\pm 1\}^{d+1}$
via the isomophism $G(F)/G(F)^+\stackrel{\det}{\lra} F^*/F^*_+\cong
\Sigma$ (compare Remark \ref{remark:sigmaaction}). Hence $\Sigma$ acts
on $H^q_*(G(F)^+, \sA(\ua_{S_1}, S_2, M; \bfC))$ and $H^q(G(F)^+, \sA(K, \ua_{S_1}, S_2, M; \bfC))$ by conjugation.

\begin{prop}
\label{prop:shimuraisomorphism}
Let $S\subseteq S_p$ and let $\bfC$ be a field containing the field of definition of $\pi\in \fA_0(G, \uzwei, \ualpha)$.

\noi (a) The $G(\bA^{S, \infty})$-representation $H^q_*(G(F)^+, \sA(\ua_{S}, \cM; \bfC))$ is of automorphic type for all $q\in \bZ$.

\noi (b) Let $\umu\in \Sigma$ and $q\in \{0,1, \ldots, d\}$. Then
\[ 
H^q_*(G(F)^+, \sA(\ua_{S}, \cM; \bfC))_{\pi, \umu} \,\,\, = \,\,\, \left \{\begin{array}{ll} 
\bfC & \mbox{if $q=d$;}\\
0 & \mbox{if $q\le d-1$.} 
\end{array}\right.
\]
\end{prop}

{\em Proof.} Firstly, we assume $S= \emptyset$. Consider the long exact sequence\begin{eqnarray*}
\label{divcoh}
&\ldots \lra H^q_*(G(F)^+, \sA(\bfC))\lra  H^q_*(B(F)^+, \sA(\bfC)) \lra \hspace{2cm}\\ 
& \lra H^q_*(G(F)^+,\sA(\cM; \bfC)) \lra H^{q+1}_*(G(F)^+, \sA(\bfC)) \ldots
\end{eqnarray*}
associated to \eqref{div} (the second group is defined
similarly as the direct limit of the groups $H^q(B(F)^+, \sA(K,
\bfC))$). The action of $G(\bA^{\infty})$ on $H^q_*(G(F)^+,\break
\sA(\bfC))$ and $H^q_*(B(F)^+,\sA(\bfC))$ has been determined in
\cite{harder}. As a $G(\bA^{\infty})$-module $H^q_*(G(F)^+,\sA(\bfC))$ is a direct sum of onedimensional representations except for $q=d+1$ in which
case there exists a $G(\bA^{\infty})$-stable decomposition 
\[
H^{d+1}_*(G(F)^+,\sA(\bfC)) \,\, = \,\, H^{d+1}_{\cusp} \oplus H^{d+1}_{\res} \oplus H^{d+1}_{\Eis}.
\]
Again, the action of $G(\bA^{\infty})$ on the second and third summand is
direct sum of onedimensional representations. On the first factor it is of automorphic type and we have
(\cite{harder}, 3.6.2.2)
\[
H^{d+1}_{\cusp}(G(F)^+,\sA(\bfC))_{\pi,\umu} \,\,\, = \,\,\, \bfC. 
\]
Using \ref{lemma:ordinaryaut}, (a) can now be easily deduced from the case $S=\emptyset$. For (b) we may pass to the $K_0(\fm)^S$-invariant part $H^q(G(F)^+, \sA(\fm, \ua_{S}, \cM; \bfC))_{\umu}$ where $\fm$ denotes the maximal prime-to-$S$ divisor of $\ff(\pi)$. By keeping in mind that the Hecke operator $T_{\fp}$ (resp.\ $U_{\fp}= -W_{\fp}$) acts by multiplication with $a_{\fp}$ (resp.\ $\pm 1$) on $H^d_*(G(F)^+,\sA(\ua_{S'},\cM; \bfC))_{\pi}$ for $\fp\in S$ with $\alpha_{\fp} \ne \pm 1$ (resp.\ $\alpha_{\fp} = \pm 1$) and $S'\subseteq S-\{\fp\}$, (b) can also be deduced from the case $S = \emptyset$ using Lemma \ref{lemma:ordinaryaut}. \enddemo

\subsection{Eichler-Shimura map}
\label{subsection:eichshim}

For the rest of this chapter we change the notation slightly again and denote by $S_1=\{\fp_1, \ldots, \fp_r\}$ the subset  of $\fp\in S_p$ such that $\alpha_{\fp} = 1$ (i.e.\ $\pi_{\fp} = \St$) and put $S_2 = S_p-S_1$ (thus with our previous notation we have $S_1=\{\fp_1, \ldots, \fp_r\}$ and $S_2=\{\fp_{r+1}, \ldots, \fp_m\}$). 

Let $K\subseteq G(\bA^{p, \infty})$ be a compact open subgroup. Our aim is to define a $\cH^{p,\infty}_K$-equivariant homomorphism (Eichler-Shimura map)
\begin{equation}
\label{eichshiml2} S_2(G, K, \ua)\to H^d(G(F)^+, \sA(K, \ua, \cM; \bC)),\,\,\, \Phi\mapsto \kappa_{\Phi}.
\end{equation}
Its definition is similar to \eqref{eichlershimura2} (the role of the manifold $U_{\infty}$ with its natural $F^*_+$-action is replaced by $\bH^{d+1}$ with $G(F)^+$ acting on it).
Firstly, define
\begin{equation}
\label{eichshiml3}
I_0: S_2(G, K, \ua) \lra H^0(G(F)^+, \sA(K, \ua, \cM; \Omega_{\hol}^d(\bH^d)))
\end{equation}
by 
\[
\langle I_0(\Phi), \psi\rangle(g, m) = \int_m \, \omega_{\langle \Phi, \psi\rangle}(1,g)
\] 
for $\psi\in C_{\ua}(F_p,\bC)$, $g\in G(\bA^{p, \infty})$ and $m\in \cM$ (here $(1,g)$ denotes the element $(1_{G(F_p)},g)\in G(F_p)\times G(\bA^{p,\infty}) = G(\bA^{\infty})$). That the image of $I_0$ is $G(F)^+$-invariant follows from \eqref{intginv} by a tedious but straightforward computation. Note that the complex $C^{\bcdot} \, \colon \!\! \,= \,\,\sA(K,\ua, \cM; \Omega_{\hol}^{\bcdot}(\bH^d))$ is a resolution of $\sA(K,\ua, \cM;\bC)$ and we have $C^q=0$ for $q>d$. We define \eqref{eichshiml2} to be the composite of \eqref{eichshiml3} with the edge morphism 
\[
H^0(G(F)^+, C^d)\lra H^d(G(F)^+,C^{\bcdot})\cong H^d(G(F)^+, \sA(K, \ua, \cM;\bC))
\]
of the spectral sequence $E^{pq}_1= H^q(G(F)^+, C^p) \!\Rightarrow \! E^{p+q}= H^{p+q}(G(F)^+,C^{\bcdot})$.

Next we define two maps 
\begin{eqnarray}
\label{deltaglob1}
&&\Delta^{\ualpha}:  S_2(G, \fm, \ua) \,\,\,  \lra \,\,\, \cD(\bG_m, S_1)\\
\label{deltaglob2}
&&\Delta^{\ualpha}: \sA(\fm, \ua_{S_0}, S_p, \cM; N)\,\,\,  \lra \,\,\, \sD(S_0\cap S_1, S_0\cap S_2, N)
\end{eqnarray}
(for $S_0$ any subset of $S_p$) which are global analogues of the map $\delta^{\ualpha}$ defined in section \ref{subsection:diststeinberg}. The first is given by
\begin{eqnarray*}
\label{deltaglob1a}
\Delta^{\ualpha}(\Phi)(f, x^p) & = & 
\delta^{\ualpha}\left(\Phi\left(\begin{matrix} x^p & 0\\ 0 & 1\end{matrix}\right) \right)(f) \,\,\, =\,\,\, \langle\left(\Phi\left(\begin{matrix} x^p & 0\\ 0 & 1\end{matrix}\right) \right), \delta_{\ualpha}(f) \rangle\end{eqnarray*}
for $f\in \sC^0(S_1, S_2, \bZ)$ and $x^p \in \bI^p$ and the second by 
\begin{eqnarray*}
\Delta^{\ualpha}(\Phi)(f, x^{p, \infty}) & = &\delta^{\ualpha}\left(\Phi\left(\left(\begin{matrix} x^{p, \infty} & 0\\ 0 & 1\end{matrix}\right) , \infty - 0\right)\right)(f)\\
& = &  \langle \left(\Phi\left(\left(\begin{matrix} x^{p, \infty} & 0\\ 0 & 1\end{matrix}\right) , \infty - 0\right)\right), \delta_{\ualpha}f\rangle
\end{eqnarray*}
for $f\in C_c^0(F_{S_0\cap S_1} \times F_{S_0\cap S_2}^*,\bZ)$ and $x^{p, \infty}\in \bI^{p, \infty}$ (as usual $\infty, 0$ denote the points  $[1:0]$ and $[0:1]$ of $\bP^1(F)$ so $\infty - 0\in \cM$). 

One checks easily that \eqref{deltaglob2} is $T(F)$-equivariant (we let $T(F)$ act on  $\sD(S_0\cap S_1, S_0\cap S_2, N)$ via the identification $T = \bG_m$). Hence the maps $T(F)^+ \to F^*_+$, $\Delta^{\ualpha}: \sA(\ua_{S_0}, \cM; N)\to \sD(S_0\cap S_1, S_0\cap S_2, N)$ induce a $\Sigma$-equivariant homomorphism
\begin{equation}
\label{deltacohom}
 H^q(G(F)^+, \sA(\fm, \ua_{S_1}, S_p, \cM; N)) \lra H^q(F^*_+,\sD(S_0\cap S_1, S_0\cap S_2, N)).
\end{equation}
The proof of the following lemma is straightforward and will be left to the reader.

\begin{lemma}
\label{lemma:eichlershimcomp}
The following diagram commutes
\begin{equation*}
\label{eichshimcomp}
\begin{CD}
S_2(G, \fm, \ua) @> \eqref{eichshiml2} >> H^d(G(F)^+, \sA(\fm, \ua, \cM; \bC))\\
@VV \Delta^{\ualpha} V @VV \eqref{deltacohom} V\\
\cD(\bG_m,S_1) @> \eqref{eichlershimura2} >>  H^d(F^*_+,\sD(S_1, S_2, \bC))
\end{CD}
\end{equation*}
\end{lemma}

\subsection{$p$-adic measures attached to Hilbert modular forms}
\label{subsection:pmeashmf}

Let $\pi\in \fA_0(G, \uzwei, \ualpha)$. As in section \ref{subsection:localdistr} let $\mu_{\alpha_i}$ denote the distribution $\chi_{\alpha_i}(x)\psi_{\fp_i}(x) dx$ on $F_i$ (resp.\ $F_i^*$) if $i\le r$ (resp.\ $i> r$) and let $\mu_{\pi_p} = \mu_{\alpha_1}\times \ldots \times \mu_{\alpha_m}$ be the product distribution on $F_{S_1} \times F_{S_2}^*$.

For $v\in \bfP_F$ let $\cW_v$ denote the Whittaker model of $\pi_v$ and let $\cW= \cW(\pi)$ be the global Whittaker model. We can choose $W_v\in \cW_v$ such that the local zeta function 
\[
\zeta(s,W_v, \chi_v) = \int_{F_v^*} \, W_v\left(\begin{matrix} x & 0\\ 0 & 1\end{matrix}\right) \, \chi_v(x) \, |x|^{s-\frac{1}{2}} d^{\times} x
\]
is equal to the local $L$-factor $L(s, \pi_v\otimes \chi_v)$ for any unramified quasicharacter $\chi_v: F^*_v\to \bC^*$ and $\Real(s) \gg 0$. In fact if $v$ is finite we can (and will) take $W_v$ to be $K_0(\ff(\pi_v))_v$-invariant. It is then uniquely determined (\cite{casselman}, Theorem 1). If $\pi_v$ is spherical (i.e.\ $\ord_v(\ff(\pi)) = 0$), then $W_v= W_{v,0}$ is the unique $G(\cO_v)$-invariant element of $\cW_v$ with $W_{v,0}(g) = 1$ for all $g\in G(\cO_v)$. Put $W^p(g) = \prod_{v\nmid p}\, W_v(g_v)$ for $g = (g_v)\in G(\bA^p)$.

We define $\phi= \phi_{\pi}: C^0_c(F_{S_1} \times F_{S_2}^*, \bZ)\times \bI^p/U^{p,\infty} \to \bC$ in $\cD(\bG_m, S_1)$ by 
\begin{equation*}
\label{phipi}
\phi(f, x^p) \colon \!\! = \sum_{\zeta\in F^*}\,\, \mu_{\pi_p}(\zeta f)\,\, W^p\left(\begin{matrix} \zeta x^p & 0\\ 0 & 1\end{matrix}\right).
\end{equation*}
That $\phi_{\pi}$ is well-defined follows from \ref{prop:whittdelta} (b). In fact for $f\in C^0_c(F_{S_1} \times F_{S_2}^*, \bZ)$ there exists an element $W_f$ of the Whittaker model $\cW_p$ of $\pi_p = \otimes_{v\in S_p} \pi_v$ such that
\[
\phi_f(x) \colon \!\! = \phi(x_p f, x^p) = \sum_{\zeta\in F^*}\,\, W\left(\begin{matrix} \zeta x & 0\\ 0 & 1\end{matrix}\right)
\]
where $W(g) \colon \!\! = W_f(g_p) W^p(g^p)$ for $g= (g_p, g^p) \in G(\bA)$. It follows $\phi_{\pi}\in \cD(\bG_m, S_1)$.

Let $\mu_{\pi} \colon \!\! = \mu_{\phi_{\pi}}$ be the corresponding distribution on $\cG_p$ and $\kappa_{\pi} \colon \!\! = \kappa_{\phi_{\pi}} \in H^d(F^*_+, \sD(S_1, S_2, \bC))$. We have 

\begin{prop}
\label{prop:interpol}
(a) (Interpolation property) Let $\chi: \cG_p\to \bC^*$ be a character of finite order with conductor $\ff(\chi)$. Then
\[
\int_{\cG_p} \,\chi(\gamma) \mu_{\pi}(d\gamma) \,\, = \,\, \tau(\chi) \,\,\prod_{\fp\in S_p} e(\alpha_{\fp},\chi_{\fp}) \,\cdot \,L(\einhalb, \pi\otimes \chi)
\]
where
\[
e(\alpha_{\fp},\chi_{\fp}) = \left \{\begin{array}{ll} 
(1- \alpha_{\fp}\chi(\varpi)^{-1}) & \mbox{if $\ord_{\fp}(\ff(\chi)) =0, \alpha_{\fp} = \pm 1$;}\\
(1-\frac{\chi(\varpi)}{\alpha_{\fp}})(1-\frac{1}{\alpha_{\fp}\chi(\varpi)}) & \mbox{if $\ord_{\fp}(\ff(\chi)) =0, \alpha_{\fp} \ne \pm 1$;}
\\
\alpha_{\fp}^{-\ord_{\fp}(\ff(\chi))} & \mbox{if $\ord_{\fp}(\ff(\chi)) >0$.} 
\end{array}\right.
\]
(b)  Let $U_0 = \prod_{\fp\in S_1} \cO_{\fp}\times  \prod_{\fp\in S_2} \cO_{\fp}^*$ and put $\phi_0 \colon \!\! = (\phi_{\pi})_{U_0}$. Then,
\[
\int_{\bI/F^*}\,\phi_0(x)\, d^{\times}x \,\, =\,\,\prod_{\fp\in S_2} e(\alpha_{\fp},1) \,\cdot \,L(\einhalb, \pi)
\]
(c) $\kappa_{\pi}$ is integral (compare Def.\ \ref{definition:integral}). For $\umu\in \Sigma$ let $\kappa_{\pi, \umu}$ denote the projection of $\kappa_{\pi}$ onto $H^d(F^*_+, \sD(S_1, S_2, \bC))_{\umu}$. Then $\kappa_{\pi, \umu}$ is integral of rank $\le 1$.
\end{prop}

{\em Proof.} (a) We view $\chi$ as a character of $\bI/F^*$ and choose an open subgroup $U$ of $U_p$ which lies in the kernel of $\chi_p = \chi|_{F_p^*}$. Let $W_U\colon = W_{1_U}\in \cW_p$ be as above and $W(g) \colon \!\! = W_U(g_p) W^p(g^p)\in \cW$. It suffices to show
\begin{eqnarray*}
\label{interpol}
&&[U_p:U]\int_{\bI/F^*} \chi(x) |x|^s \phi_U(x)\, d^{\times} x \\
&& \hspace{2cm} = \,\, \Norm(\ff(\chi))^s\tau(\chi ) \,\,\prod_{\fp\in S_p} e(\alpha_{\fp},\chi_{\fp}|\wcdot|_{\fp}^s) \,\cdot \,L(s + \einhalb, \pi\otimes \chi)
\end{eqnarray*}
for all $s\in \bC$ with $\Real(s) >-1$. Since both sides are holomorphic functions it is enough to prove this for $\Real(s) \gg 0$. Using \ref{prop:whittdelta} (c) and \ref{prop:localmeasurespherical} we obtain 
\begin{eqnarray*}
&  [U_p:U]\int_{\bI/F^*} \chi(x) |x|^s \phi_U(x)\, d^{\times}x =  [U_p:U]\int_{\bI} \chi(x) |x|^s W\left(\begin{matrix} x & 0\\ 0 & 1\end{matrix}\right)\,d^{\times}x\\
&  = \,\,[U_p:U]\int_{F^*_p} \chi_p(x) |x|^s  W_U\left(\begin{matrix} x & 0\\ 0 & 1\end{matrix}\right) d^{\times}x  \int_{\bI^p} \chi^p(y) |y|^s W^p\left(\begin{matrix} y & 0\\ 0 & 1\end{matrix}\right) d^{\times}y\\
&  = \,\, \prod_{\fp\in S_p}\int_{F^*_{\fp}} \chi_{\fp}(x) |x|_{\fp}^s \,\,\mu_{\alpha_{\fp}}(dx)\,\, \cdot L_{S_p}(s + \einhalb, \pi\otimes \chi)\hspace{3cm}\\
& = \,\,\Norm(\ff(\chi))^s\tau(\chi ) \,\,\prod_{\fp\in S_p} e(\alpha_{\fp},\chi_{\fp}|\wcdot|_{\fp}^s) \,\cdot \,L(s + \einhalb, \pi\otimes \chi).
\hspace{2.4cm}\end{eqnarray*}
(b) Again it suffices to show that
\begin{equation*}
\label{Lwerteinhalb}
\int_{\bI/F^*}\,|x|^s\phi_0(x) d^{\times}x \,\, =\,\,\prod_{\fp\in S_2} e(\alpha_{\fp},|\wcdot|_{\fp}^s) \,\cdot \,L(s+\einhalb, \pi)
\end{equation*}
for $\Real(s) \gg 0$. A similar computation as above yields
\begin{eqnarray*}
&\int_{\bI/F^*}\,|x|^s\phi_0(x) d^{\times}x  \,\, =\,\,\hspace{5cm} \\
&\prod_{\fp\in S_1}\!\int_{F^*_{\fp}} |x|_{\fp}^s \mu_1(x\cO_{\fp})  d^{\times}x \wcdot \prod_{\fp\in S_2}\! \int_{F^*_{\fp}} |x|_{\fp}^s \mu_{\alpha_{\fp}}(x\cO_{\fp}^*)  d^{\times}x \wcdot L_{S_p}(s + \einhalb, \pi).
\end{eqnarray*}
For $\fp\in S_1$ we have $\mu_1(x\cO_{\fp})  = \int_{F_{\fp}} 1_{x\cO_{\fp}}(y)\psi_{\fp}(y)  dy = |x|_{\fp}1_{\cO_{\fp}}(x)$ hence 
\begin{eqnarray*}
\label{Lwerteinhalb1}
\int_{F^*_{\fp}} |x|_{\fp}^s \mu_1(x\cO_{\fp})  d^{\times}x & = &\int_{F^*_{\fp}} |x|_{\fp}^{s+1} 1_{\cO_{\fp}}(x) d^{\times}x\\
 &=& (1-\Norm(\fp)^{-(s+1)})^{-1} \,\,= \,\, L(s+\einhalb, \pi_{\fp}).
\end{eqnarray*}
On the other hand by Prop.\ \ref{prop:localmeasurespherical} we get for $\fp\in S_2$ 
\begin{eqnarray*}
\label{Lwerteinhalb2}
\int_{F^*_{\fp}} |x|_{\fp}^s \,\,\mu_{\alpha_{\fp}}(x\cO_{\fp}^*)  d^{\times}x = \int_{F^*_{\fp}} |x|_{\fp}^s \,\,\mu_{\alpha_{\fp}}(dx) \,\,= \,\,  e(\alpha_{\fp},|\wcdot|_{\fp}^s) \,\cdot \,L(s+\einhalb,\pi_{\alpha}).
\end{eqnarray*}
The assertion follows.

\noi (c) Let $\lambda_{\ua}\in \sB^{\ua}(F_p,\barQ)$ be the image of $\otimes_{i=1}^m \, \lambda_{\alpha_i}$ under \eqref{prodsteinbergdual} and define $\Phi_{\pi}\in \cA_0(G, \hol, \uzwei, \ua)$ by 
\begin{equation*}
\label{interpol1}
\langle \psi, \Phi_{\pi}(g^p) \rangle \,\, = \,\, \sum_{\zeta\in F^*}\,\, \langle \left(\begin{matrix} \zeta & 0\\ 0 & 1\end{matrix}\right) \psi, \lambda_{\ua}\rangle \,  W^p\left( \left(\begin{matrix} \zeta  & 0\\ 0 & 1\end{matrix}\right) g^p\right)
\end{equation*}
for $g^p\in G(\bA^p)$ and $\psi\in \sB_{\ua}(F_p,\bC)$. To see that $\Phi_{\pi}$ satisfies property (ii) of Def.\ \ref{definition:auttree}, let $\psi\in \sB_{\ua}(F_p,\bC)$ and define $W_{\psi}\in \cW(\pi)$ by $W_{\psi}(g_p, g^p) \colon \!\! = \langle g_p \psi, \lambda_{\ua}\rangle \, W^p(g_p)$. Then 
\[
\langle \psi, \Phi_{\pi} \rangle(g) \,\, = \,\, \sum_{\zeta\in F^*}\,\, W_{\psi}\left( \left(\begin{matrix} \zeta & 0\\ 0 & 1\end{matrix}\right)  g \right)\in \cA_0(G, \hol, \uzwei).
\] 
Let $\fm$ be the maximal prime-to-$p$ divisor of $\ff(\pi)$. Because
$W^p$ is $K_0(\fm)$-invariant we have $\Phi_{\pi}\in S_2(G, \fm,
\ua)$. Since $\Delta^{\ualpha}(\Phi_{\pi}) = \phi_{\pi}$ by
\ref{prop:whittdelta} (a), we can apply Lemma
\ref{lemma:eichlershimcomp} to conclude that $\kappa_{\pi}$ lies in
the image of \eqref{deltacohom}. That $\kappa_{\pi}$ is
integral now follows from the fact that $R \mapsto H^d(G(F)^+,
\sA(\fm, \ua, \cM; R))$ commutes with flat base change (Prop.\
\ref{prop:basechange}). The second assertion is a consequence of
\ref{prop:basechange} and \ref{prop:shimuraisomorphism}.
\enddemo

Let $\pi\in \fA_0(G, \uzwei, \ualpha)$. By \ref{corollary:arithmetmeas}
and \ref{prop:meascoh} the distribution $\mu_{\pi}$ is a $p$-adic
measure. We define the {\it $p$-adic $L$-function of $\pi$} by 
\begin{equation*}
\label{plfunc}
L_p(s, \pi) \colon \!\! = L_p(s, \kappa_{\pi}) = L_p(s \kappa_{\pi, +})=\int_{\cG_p}\, \langle \gamma\rangle^s \mu_{\pi}(d\gamma) \hspace{0.5cm} \mbox{for $s\in \bZ_p$.}
\end{equation*}
It is a locally analytic function with values in the one-dimensional vector space $L_{\kappa_{\pi,+}}\otimes_{\barcO}\,\bC_p$ (compare Remark \ref{remark:sigmaaction}).

\section{Exceptional zero conjecture}
\label{section:ezc}

\subsection{Automorphic $\cL$-invariants}
\label{subsection:linv}

We keep the notation and assumptions from the end of last section. In this section we define for each $\fp\in S_1$ a certain number $\cL_{\fp}(\pi)\in \bC_p$, the {\it $\cL$-invariant of $\pi$ at $\fp$}. It has the important property that it does not change under suitable quadratic twists (see Lemma \ref{lemma:twistinv} below).

Let $\cO$ denote the valuation ring of $\bC_p$. Fix $\fp\in S_1$ and put $S_0= S_p-\{\fp\}$. Recall that by Remark \ref{remark:afcofree}) there exists a canonical pairing
\begin{eqnarray}
\label{pairkey}
&& \sA(K, \ua, \cM; \cO)\otimes_{\cO} \bC_p \times \St(F_{\fp}, \bC_p)\\
&& \lra\,\, \sA(K, \ua_{S_0}, S_p, \cM; \cO)\otimes_{\cO} \bC_p\subseteq \sA(K, \ua_{S_0}, S_p, \cM; \bC_p).\nonumber
\end{eqnarray}
It induces a cup-product pairing on $G(F)^+$-cohomology. Together with \ref{prop:basechange} this yields a pairing 
\begin{eqnarray}
\label{eqn:pairkeycohom}
&& \cup: H^p(G(F)^+, \sA(K, \ua, \cM;\bC_p)) \times H^q(G(F)^+, \St(F_{\fp}, \bC_p))\\
&& \hspace{2cm}\lra H^{p+q}(G(F)^+,\sA(K, \ua_{S_0}, S_p, \cM; \bC_p))\nonumber
\end{eqnarray}
Hence by passing to the direct limit over all $K$ we obtain a homomorphism of smooth $G(\bA^{p, \infty})$-representations 
\begin{equation*}
\label{eqn:pairkeycohomlim}
\wcdot\cup b: H^p_*(G(F)^+, \sA(\ua, \cM;\bC_p)) \to  H^{p+q}_*(G(F)^+,\sA(\ua_{S_0}, S_p, \cM; \bC_p))
\end{equation*}
for all $b\in H^q(G(F)^+, \St(F_{\fp}, \bC_p))$.

\begin{remarks}\rm 
\noi (a) Note that in \eqref{pairkey} we cannot replace
$\sA(\fm, \ua, \cM; \cO)\otimes_{\cO} \bC_p$ by $\sA(\fm,
\ua, \cM; \bC_p)$. Therefore the compatibility of $R \mapsto
H^{q_1}(G(F)^+, \break\sA(\fm, \ua, \cM; R))$ with flat base change is crucial in the definition of (\ref{eqn:pairkeycohom}).

\noi (b) Let $\fm$ be an ideal of $\cO_F$ which is relatively prime to $p\cO_F$.  We have a commutative diagram  
\begin{equation*}
\label{deltacup1}
\begin{CD}
\sA(\fm, \ua, \cM; \cO)\otimes_{\cO} \bC_p \times \St(F_{\fp}, \bC_p) @>>> \sA(\fm, \ua_{S_0}, S_p, \cM;\cO)\otimes_{\cO} \bC_p\\
@VV \Delta^{\ualpha}\times \delta_1^{-1} V @VV \Delta^{\ualpha} V\\
\sD^b(S_1, S_2, \bC_p) \times \Cd(F_{\fp}, \bC_p)   @>>> \sD^b(S_0\cap S_1, S_0\cap S_2, \bC_p)
\end{CD}
\end{equation*}
(the top horizontal arrow is the map \eqref{pairkey}). Hence for $b\in H^q(G(F)^+, \St(F_{\fp}, \break \bC_p))$ the diagram 
\begin{equation*}
\label{deltacup2}
\begin{CD}
H^p(G(F)^+ \! ,\sA(\fm, \ua, \! \cM;\bC_p)) @> \wcdot\cup b >>  H^{p+q}(G(F)^+ \! ,\sA(\fm, \ua_{S_0}, S_p, \! \cM; \bC_p))\\
@VV \eqref{deltacohom} V @VV \eqref{deltacohom} V\\
H^p(F^*_+,\sD^b(S_1, S_2, \bC_p)) @> \wcdot\cup\delta^*(b) >> H^{p+q}(F^*_+,\sD^b(S_0\cap S_1, S_0\cap S_2, \bC_p))
\end{CD}
\end{equation*}
commutes as well. Here $\delta^*: H^q(G(F), \St(F_{\fp}, \bC_p)) \to H^q(F^*, \Cd(F_{\fp}, \bC_p))$ is the canonical map induced by $\delta_1^{-1}$ (compare \ref{lemma:extsteinberg} (b)) and the first vertical map is given by 
\begin{eqnarray*}
\label{deltameas}
& H^p(G(F)^+,\sA(\fm, \ua, \cM; \bC_p)) \, \cong \, H^p(G(F)^+,\sA(\fm, \ua, \cM; \cO))\otimes_{\cO} \bC_p \\
& \stackrel{\eqref{deltacohom}}{\lra} H^p(F^*_+,\sD(S_1,
S_2, \cO))\otimes_{\cO} \bC_p  \lra  H^p(F^*_+,\sD^b(S_1, S_2, \bC_p))
\end{eqnarray*}
The definition of the second vertical map is analogous.
\end{remarks}

 The extension classes \eqref{steinext} associated to the
 homomorphisms $\ord_{\fp}: F_{\fp}^* \to \break\bC_p$ and $\ell_{\fp} =
 \log_p\circ \Norm_{F_{\fp}/\bQ_p}: F^* \to \bC_p$ define two
 cohomology classes $b_{\ord, \fp}$, $b_{\log, \fp}\in H^1(G(F)^+, \St(F_{\fp}, \bC_p))$. 

\begin{lemma}
\label{lemma:cupiso}
(a) The $G(\bA^{p, \infty})$-representation $H^q_*(G(F)^+,\sA(\ua_{S_0}, S_p, \cM; \break\bC_p))$ is of automorphic type for all $q\in \bZ$. 

\noi (b) For $\umu\in \Sigma$ the map $\wcdot\cup b_{\ord, \fp}$ induces an isomorphism 
\begin{equation}
\label{eqn:cupord}
H^d_*(G(F)^+, \sA(\ua, \cM;\bC_p))_{\pi, \umu} \lra H^{d+1}_*(G(F)^+,\sA(\ua_{S_0}, S_p, \cM; \bC_p))_{\pi, \umu}
\end{equation}
of onedimensional $\bC_p$-vector spaces. 
\end{lemma}

{\em Proof.} Let $K$ be a compact open subgroup of $G(\bA^{p, \infty})$  and put 
$K_0 = K \times G(\cO_{\fp})\subseteq G(\bA^{S_0, \infty})$. We define the $\cH^{S_p\cup S_{\infty}}_{K}$-module $\sQ_K$ by 
\begin{equation*}
\label{tree3}
\xymatrix@-0.9pc{
0 \ar[r] & \sQ_K \ar[r] & \sA(K_0, \ua_{S_0}, \cM;\bC_p) \ar[rr]^{T_{\fp} - (q+1)} & &\sA(K_0, \ua_{S_0}, \cM;\bC_p) \ar[r] & 0. }
\end{equation*}
By considering the corresponding long exact cohomology sequence we obtain that $H^{\bcdot}(G(F)^+, \sQ_K)$ is of automorphic type and $H^{\bcdot}(G(F)^+, \sQ_K)_{\pi} = 0$. For the latter note that $H^{\bcdot}(G(F)^+, \sA(K_0, \ua_{S_0}, \cM;\bC_p))_{\pi} = 0$ since $K_0 \not\subseteq K_0(\ff(\pi))^{S_0\cup S_{\infty}}$. By \eqref{heckeharmonic2} there exists an exact sequence 
\begin{equation*}
\label{tree4}
\xymatrix@-0.9pc{
0 \ar[r] & \sA(K, \ua_{S_0}, S_p, \cM; \bC_p) \ar[r] & \sQ_K \ar[r] &\sA(K, \ua, \cM;\bC_p)\ar[r] & 0. }
\end{equation*}
It follows from \ref{lemma:extsteinberg} (c) that (\ref{eqn:cupord}) is equal to the connecting homomorphism in the corresponding long exact sequence in degree $d$ (up to sign). Hence the assertions follow from Prop. \ref{prop:shimuraisomorphism}.\enddemo

\begin{definition}
\label{definition:cupiso}
For $\umu\in \Sigma$ there exists a unique $\cL_{\fp}(\pi, \umu)\in \bC_p$ -- called the $\cL$-invariant of $\pi$ at $\fp$ -- such that 
\[
(\wcdot\cup b_{\log, \fp})_{\pi,\umu}  \,\,\, = \,\,\, \cL_{\fp}(\pi, \umu) \,\,(\wcdot\cup b_{\ord, \fp})_{\pi,\umu}.
\]
 If $\umu = (1, \ldots, 1)$ then we write $\cL_{\fp}(\pi) = \cL_{\fp}(\pi, \umu)$.
\end{definition}

\begin{conjecture}
\label{conjecture:indmu}
$\cL_{\fp}(\pi, \umu)$ is independent of the choice of $\umu\in \Sigma$, i.e.\ we have $\cL_{\fp}(\pi, \umu) = \cL_{\fp}(\pi)$ for all $\umu\in \Sigma$.
\end{conjecture}

\begin{lemma}
\label{lemma:twistinv}
Let $\chi: \bI/F^* \to \{\pm 1\}$ be a quadratic character whose conductor is prime to $p$ 
and such that $\chi_{\fp} =1$. Then
\[
\cL_{\fp}(\pi, \umu) \,\,\, = \,\,\, \cL_{\fp}(\pi\otimes \chi, \sign(\chi)\umu)
\]
where $\sign(\chi) = \chi_{\infty}(-1,\ldots, -1)\in \Sigma$. 
\end{lemma}

{\em Proof.} For a smooth semisimple representation $V$ of $G(\bA^{p,\infty})$ we denote by $V_{\chi}$ the representation $V\otimes \det\circ \chi^{p,\infty}$. Note that $(V_{\chi})_{\pi} = V_{\pi\otimes \chi}$. Put $\uep \colon \!\! = \uep(\chi_p) = (\chi_{\fp_1}(\varpi_1), \ldots, \chi_{\fp_m}(\varpi_m))$ (compare section \ref{subsection:semiloc}).  We define a {\it twisting operator}
\begin{equation*}
\label{twist4}
\Tw_{\chi}: \sA(\ua, S_p, \cM; \bC_p) \lra \sA(\uep\ua, S_p, \cM; \bC_p)
\end{equation*}
by $\Tw_{\chi}(\Phi)(g, m) \,= \,\chi^{p,\infty}(\det(g)) \Tw_{\chi_p}(\Phi(g,m))$ for all $g\in G(\bA^{p,\infty})$, $m\in \cM$ and define $\Tw_{\chi}:\sA(\ua_{S_0}, S_p, \cM; \bC_p)\to \sA((\uep\ua)_{S_0}, S_p, \cM; \bC_p)$ analogously. Note that $\Tw_{\chi}$ is $G(F)^+$-linear and maps $\sA(K, \ua, S_p, \cM; \bC_p)$ onto $\sA(K, \uep\ua, S_p, \cM; \bC_p)$ as long as $\det(K)$ is contained in the kernel of $\chi^{p,\infty}$. Note also that $\Tw_{\chi}\circ \Tw_{\chi} = \id$. For $b\in H^q(G(F)^+, \St(F_{\fp}, \bC_p))$ we obtain a diagram
\begin{equation*}
\label{deltacup3}
\small{
\xymatrix@-1.0pc{
H^d_*(G(F)^+, \sA(\ua, \cM;\bC_p))_{\pi} \ar[rr]^{(\wcdot\cup
  b)_{\pi}} \ar[d]^{\cong} & &H^{d+1}_*(G(F)^+,\sA(\ua_{S_0}, S_p, \cM; \bC_p))_{\pi} \ar[d]^{\cong}\\
H^d_*(G(F)^+, \sA(\uep\ua, \cM;\bC_p))_{\pi\otimes\chi} 
\ar[rr]^{(\wcdot\cup b)_{\pi\otimes \chi}} & & H^{d+1}_*(G(F)^+,\sA((\uep\ua)_{S_0}, S_p, \cM; \bC_p))_{\pi\otimes\chi}
}}
\end{equation*}
where the vertical maps are induced by $\Tw_{\chi}$. The commutativity
of the diagram for $b= b_{\log, \fp}$ and $b= b_{\ord, \fp}$ implies the assertion.\enddemo

\begin{remarks} \rm (a) Conjecture \ref{conjecture:indmu} is known in
  the case $F=\bQ$ (\cite{bdi}, Theorem 6.8).

\noi (b) Let $D$ be a quaternion algebra over $F$ and let $G' =  D^*$
(viewed as an algebraic group). Let $\pi$ be an automorphic
representation of $G'(\bA)$ whose components $\pi_v$ are discrete
series for all $v\in S_{\infty}$. By a similar construction as above
one should be able to define an $\cL$-invariant $\cL_{\fp}(\pi)$
whenever $\fp$ does not divide the discriminant of $D$ and we have
$\pi_{\fp} \cong \St$. These $\cL$-invariants have been defined in the
case $F=\bQ$ and $D$ definite by Teitelbaum \cite{teitelbaum}, for $F=
\bQ$, $D = M_2(\bQ)$ by Darmon \cite{darmon} (for weight $2$) and
Orton \cite{orton} (higher even weight) and if $F$ has narrow class
number 1 by Greenberg \cite{greenberg} (in the case of parallel weight
$(2,\ldots,2)$) (it not difficult to see that for $D = M_2(F)$ the
$\cL$-invariant defined in \cite{greenberg} match with the one defined
above). An interesting and difficult problem is to show that
$\cL_{\fp}(\pi)$ is invariant under the Jacquet-Langlands
correspondence (this has been proved for $F=\bQ$ in certain cases in
\cite{bdi}, Cor.\ 6.9) or under base change. The author hopes to
return to the study of these $\cL$-invariants in the future.
\end{remarks}

\subsection{Main results}

Our first main result is

\begin{theorem}
\label{theorem:wezc} 
Let $\pi\in \fA_0(G, \uzwei, \ualpha)$. The vanishing order of $L_p(s, \pi)$ at $s=0$ is at least equal to $r$ (i.e.\ to the number of places $\fp$ of $F$ above $p$ with $\pi_{\fp} \cong \St$). Moreover we have 
\begin{equation}
\label{ezchmf}
L^{(r)}_p(0, \pi)\,\, =\,\, r! \,\, \prod_{\fp\in S_1}\cL_{\fp}(\pi) \,\cdot\,\prod_{\fp\in S_2} e(\alpha_{\fp},1) \,\cdot \,L(\einhalb, \pi)
\end{equation}
where \(e(\alpha_{\fp},1) = \left \{\begin{array}{ll} 2 & \mbox{if $\alpha_{\fp} = -1$;}\\
(1-\frac{1}{\alpha_{\fp}})^2 & \mbox{if $\alpha_{\fp} \ne \pm 1$.}\end{array}\right.
\)
\end{theorem}

{\em Proof.} The first statement follows from \ref{corollary:abstrezc2} (a). By \ref{corollary:Lwert0},
\ref{prop:interpol} (b) and \ref{corollary:abstrezc2} (b) we have
\begin{eqnarray*}
&\prod_{\fp\in S_2} e(\alpha_{\fp},1) \,\cdot \,L(\einhalb, \pi) \,\,\, = \,\,\, (-1)^{\binom{r}{2}}\,\, (\kappa_{\pi, +} \cup c_{\fp_1}\cup \ldots \cup c_{\fp_{r}}) \cap \vartheta \\
&L^{(r)}_p(0, \pi) \,\,\, = \,\,\, (-1)^{\binom{r}{2}}\,\,r! \,\, (\wkappa_{\pi, +} \cup c_{\ell_{\fp_1}} \cup \ldots \cup c_{\ell_{\fp_{r}}}) \cap \vartheta
\end{eqnarray*}
Thus it suffices to show $\wkappa_{\pi, +} \cup c_{\ell_{\fp}} =
\cL_{\fp}(\pi)\,\wkappa_{\pi, +} \cup c_{\fp}$ for $\fp\in S_1$. Let
$\fm$ be the maximal prime-to-$p$ divisor of $\ff(\pi)$. The proof of
\ref{prop:interpol} (c) and \ref{prop:shimuraisomorphism} shows that
there exists an element $\beta\in H^d(G(F)^+, \sA(\fm, \ua, \cM;
\barQ))_{\pi}$ which is mapped under \eqref{deltacohom} to
some non-zero multiple of $\kappa_{\pi}$. Also by Lemma
\ref{lemma:extsteinberg} (b) we have $2c_{\fp} = \delta^*(b_{\ord,
  \fp})$ and $2 c_{\ell_{\fp}} = \delta^*(b_{\log, \fp})$. Hence the
 assertion  follows from $\beta_+\cup b_{\ord, \fp} =
 \cL_{\fp}(\pi)\,\beta_+\cup b_{\log, \fp}$ (here we view $\beta$ as
an element of $H^d(G(F)^+, \sA(\fm, \ua, \cM; \bC_p))$).\enddemo

Now assume that $E/F$ is an elliptic curve which is {\it $p$-ordinary}, i.e.\ it has either good ordinary or multiplicative reduction at all places above $p$.  We also assume that $E$ is {\it modular} by which we mean that for some prime number $\ell$, the $\ell$-adic Tate module of $E$  is isomorphic (as a Galois representation) to the $\ell$-adic representation associated to some $\pi = \pi_E\in \fA_0(G, \uzwei)$ (this holds then for any $\ell$; compare \cite{wintenberger}). Then it is known (\cite{carayol}, \cite{taylor}) that the local $L$-factors of $E$ and $\pi$ all match up. In particular we have $\Lambda(E, \chi, s) = L(s- \einhalb,\pi\otimes \chi)$ for any character $\chi: \bI/F^*\to \bC^*$ (where $\Lambda(E, \chi, s)$ denotes the completed Hasse-Weil $L$-function) and the conductor of $E$ is $=\ff(\pi)$. Moreover $\pi$ is $p$-ordinary and $E$ has split multiplicative reduction at $\fp$ if and only if $\pi_{\fp} = \St$. Thus $\pi\in \fA_0(G, \uzwei, \ualpha)$ for some ordinary parameters $\ualpha = (\alpha_1, \ldots, \alpha_m)$ and if $S_1 = \{\fp_1, \ldots, \fp_r\}$ denotes the set of $\fp\in S_p$ where $E$ has split multiplicative reduction then $\alpha_1 = \ldots = \alpha_r=1$. For $\fp\in S_1$ recall that $\cL_{\fp}(E) = \frac{\log_p(\Norm_{F_{\fp}/\bQ_p}(q_{E/F_{\fp}}))}{\ord_{\fp}(q_{E/F_{\fp}})}$ where $q_{E/F_{\fp}}$ denotes the Tate period associated to $E/F_{\fp}$. The $p$-adic $L$-function of $E$ is defined as $L_p(E,s) \colon \!\! = L_p(s, \pi)$. Recall Hida's conjecture \cite{hida} 

\begin{conjecture}
\label{conjecture:hidaezc}
(a) $\ord_{s=0} L_p(E, s)\ge r$. 

\noi (b) $L^{(r)}_p(E, 0)\,\, =\,\, r! \,\, \prod_{\fp\in S_1}\cL_{\fp}(E) \,\cdot\,\prod_{\fp\in S_2} e(\alpha_{\fp},1) \,\cdot \, \Lambda(E, 1)$. 
\end{conjecture}
 
In the case $r=1$ this has been proved by Mok \cite{mok} (under the
additional assumption that $p$ is unramified in $F$ and $\ge 5$). We prove the first assertion unconditionally and deduce the second from Mok's result under some further (mild) restrictions using Theorem \ref{theorem:wezc} and a non-vanishing result for twisted $L$-values \cite{waldspurger}, \cite{friedhoff}. Let $w(\pi)$ denote the root number of $\pi$. We first show

\begin{prop}
\label{prop:complinv} 
Assume that $p\ge 5$ is unramified in $F$. If (i) $E$ has
multiplicative reduction at some place $\fq\nmid p$ or (ii) $r + w(\pi) \equiv 1 \mod 2$ then $\cL_{\fp}(E) = \cL_{\fp}(\pi)$ for all $\fp\in S_1$. 
\end{prop}

{\em Proof.} Suppose that there exists a quadratic character $\chi: \bI/F^* \to \{\pm 1\}$ whose conductor is relatively prime to $\ff(\pi)$ with the following properties
\begin{itemize}
\item{} $\chi_v = 1$ for all $v\in \{\fp\} \cup S_2 \cup S_{\infty}$ (where $S_2 = S_p-S_1$);
\item{} $\chi_{\fq}\ne 1$ for all $\fq\in S_1-\{\fp\}$;
\item{} $L(\einhalb, \pi\otimes \chi) \ne 0$.
\end{itemize}
Let $E_{\chi}$ denote the twist of $E$ by $\chi$. The first property
and \ref{lemma:twistinv} imply $\cL_{\fp}(E) = \cL_{\fp}(E_{\chi})$
and $\cL_{\fp}(\pi) = \cL_{\fp}(\pi\otimes \chi)$. It also follows
from the first two properties that $E_{\chi}$ is $p$-ordinary and
$\fp$ is the only place above $p$ where $E_{\chi}$ has good ordinary
reduction. Thus by (\cite{mok}, Theorem 1.1) Conj.\ \ref{conjecture:hidaezc} (b) holds for $E_{\chi}$. Since the formula coincides with \eqref{ezchmf} with the
possible exception of the $\cL$-invariants and because of the non-vanishing of $L(\einhalb, \pi\otimes \chi) $ we deduce $\cL_{\fp}(E) = \cL_{\fp}(E_{\chi}) = \cL_{\fp}(\pi\otimes \chi) = \cL_{\fp}(\pi)$.

Therefore it remains to show that under the assumptions (i) or (ii) there exists a quadratic character with the above properties. If $\chi$ is any quadratic character with $\ff(\chi)$ relatively prime to $\ff(\pi)$ then it is well known that $w(\pi \otimes \chi) = \prod_{v\in S_{\infty}} \chi_v(-1) \,\, \chi(\ff(\pi)) \,w(\pi)$. Hence under the assumptions (i) or (ii) it is clear that there exists $\chi$ so that at least the first two properties are satisfied and such that $w(\pi \otimes \chi)=1$. Then a theorem of Waldspurger (\cite{waldspurger}; see also \cite{friedhoff}, Thm.\ B) implies that we can choose $\chi$ so that $L(\einhalb, \pi\otimes \chi) \ne 0$ holds as well.\enddemo

\begin{theorem}
\label{theorem:ezc}
(a) $\ord_{s=0} L_p(E, s)\ge r$.

\noi (b) Assume $p\ge 5$ is unramified in $F$. If (i) $E$ has
multiplicative reduction at some place $\fq\nmid p$ or (ii) $r$ is odd or (iii) $w(\pi) = -1$ then Conj.\ \ref{conjecture:hidaezc} (b) holds.
\end{theorem}

{\em Proof.} (a) is part of Theorem \ref{theorem:wezc}. For (b) assume first $w(\pi) = -1$. Then both sides of the equation vanish (the left hand side by  \eqref{ezchmf}). If $w(\pi) = 1$ and (i) or (ii) hold then (b) follows from Thm.\ \ref{theorem:wezc} and Prop.\ \ref{prop:complinv}.\enddemo

\begin{remarks} \rm (a) In the case $F=\bQ$, Thm.\ \ref{theorem:wezc} is due to Darmon \cite{darmon}.

\noi (b) If Conjecture \ref{conjecture:indmu} holds then it is not necessary to assume (i) or (ii) in  Prop.\ \ref{prop:complinv} above. In fact, obviously, there exists a quadratic character $\chi$ with $\chi_v = 1 \ne \chi_{\fq}$ for all $v\in \{\fp\} \cup S_2$ and $\fq\in S_1-\{\fp\}$ and so that 
$w(\pi \otimes \chi) = 1$. Hence by Waldspurger's theorem we can choose $\chi$ which satisfies also $L(\einhalb, \pi\otimes \chi) \ne 0$. In the forthcoming work \cite{iovitaspiess} we shall give a different and unconditional proof of the equality $\cL_{\fp}(E) = \cL_{\fp}(\pi)$ (hence of Hida's conjecture) without using Mok's result.
\end{remarks}

\end{document}